\documentclass{amsart}
\usepackage{fix-cm}
\usepackage[english]{babel}
\usepackage{xcolor}
\usepackage{colortbl}
\usepackage{booktabs}
\usepackage{multirow}
\usepackage{graphicx}
\usepackage{amssymb}
\usepackage{amsmath}
\usepackage{amsfonts}
\usepackage{mathtools}
\usepackage{algpseudocode}
\usepackage{url}
\usepackage{algorithm}
\usepackage{tikz}
\usepackage{float}
\usepackage{bm}
\usepackage{multicol}
\usepackage[format=plain,labelfont={bf,sf},textfont=sf,compatibility=false]{caption}
\usepackage{subcaption}
\usepackage[title]{appendix}
\usepackage{textcomp}
\usepackage{manyfoot}
\usepackage{listings}
\usepackage[normalem]{ulem} % per \sout
\usepackage{xparse}  

\usetikzlibrary{positioning, arrows.meta, shapes}
\graphicspath{{images/}}
\definecolor{mplblue}{HTML}{1F77B4}

\newcommand{\vect}[1]{\bm{#1}}
\newcommand{\set}[1]{\mathbb{#1}}
\newcommand{\vectspace}[1]{\mathcal{#1}}
\newcommand{\norm}[1]{\parallel#1\parallel}

\newcommand{\revA}[2][]{%
  {\color{blue}#2}%
  \if\relax\detokenize{#1}\relax
  \else
    {\color{blue}\ \sout{#1}}%
  \fi
}
\newcommand{\revB}[2][]{%
  {\color{red}#2}%
  \if\relax\detokenize{#1}\relax
  \else
    {\color{red}\ \sout{#1}}%
  \fi
}

\newcommand{\revC}[2][]{%
  {\color{teal}#2}%
  \if\relax\detokenize{#1}\relax
  \else
    {\color{teal}\ \sout{#1}}%
  \fi
}

\definecolor{ssfem_green}{RGB}{92,130,67}
\definecolor{ssfem_blue}{RGB}{86,112,168}
\definecolor{ssfem_red}{RGB}{180,44,34}
\usepackage[foot]{amsaddr}
\usepackage{hyperref}
\usepackage[nameinlink,noabbrev]{cleveref}
\crefformat{subfigure}{#2(#1)#3}
\Crefformat{subfigure}{#2(#1)#3}
\usepackage[top=1in, bottom=1.25in, left=1.5in, right=1.5in]{geometry}

\title[Stochastic bifurcations in CFD]{A stochastic perturbation approach to nonlinear bifurcating problems}

\author{Isabella Carla Gonnella$^1$, Moaad Khamlich$^1$, Federico Pichi$^2$, Gianluigi Rozza$^1$}
\address{$^1$ mathLab, Mathematics Area, SISSA, via Bonomea 265, I-34136 Trieste, Italy}
\address{$^2$ MCSS, École Polytechnique Fédérale de Lausanne, Lausanne, 1015, Switzerland
}
\date{}

\begin{document} 

\begin{abstract}
Incorporating probabilistic terms in mathematical models is crucial for capturing and quantifying uncertainties in real-world systems, especially when the solution is not unique or exhibits sudden qualitative changes as parameters vary. 
However, stochastic models typically require large computational resources to produce meaningful statistics. 
In this work, we leverage the Polynomial Chaos (PC) expansion to propose a systematic approach for bifurcation detection in parametric systems of equations.
We show that the method, exploiting a perturbed version of the deterministic model, avoids repeated costly simulations across multiple parameter values and requires no prior information for initializing numerical solvers, while still providing accurate characterization of the bifurcation branches. 
We argue that the PC solutions of the perturbed model not only provide access to statistical information about the deterministic branches, but also approximate these branches in a meaningful sense.
Finally, we validate our claims by means of numerical tests on the pitchfork bifurcation, examining both its normal form and a classical realization in fluid-dynamics PDEs, namely the Coand\u{a} effect.


\medskip
\begin{center}
    \textbf{Code availability:} \url{https://github.com/ICGonnella/SSFEM-Coanda-Effect.git}
\end{center}
\end{abstract}

\maketitle

\section{Introduction}

Nonlinear parametric partial differential equations (PDEs) play a fundamental role in accurately modeling several physics problems, but their solutions may exhibit a non-trivial dependence on the system parameters.
This is indeed the case of bifurcating phenomena \cite{Seydel}, where the nonlinear terms in the equations give rise to a non-differentiable evolution of the solution manifold, and thus to the non-uniqueness of the solution with respect to the parameters. More specifically, the problem admits multiple co-existing states for the same value of a parameter $\mu$ (e.g.\ controlling the physical properties of the model), and an in-depth numerical investigation is needed to deal with such ill-posedness to fully describe the system in the simulations.
Studying a model's bifurcating behavior means evaluating the curve individuated in the state space by the solutions obtained solving the equations for different values of the parameters.  An intuitive way of describing the bifurcation is to represent it via the \textit{bifurcation diagram}, which consists in tracking a specific scalar quantity of interest, characterizing the different \textit{branches}
w.r.t.\ the corresponding parameter value, possibly including information regarding their stability properties.

The complexity of the topic, coupled with its wide range of potential applications, has stimulated the development of numerical methods to investigate such problems and enhance their exploitation. 
In particular, most numerical approaches to bifurcation analysis are commonly referred to as \textit{continuation methods}, and consist in the use of deterministic numerical solvers in an iterative pipeline, progressively exploring the parametric space. 
The idea behind this approach is to accurately ``continue" an already known solution branch in the yet unexplored regions, following its evolution. Nevertheless, the local nature of continuation methods creates a problem when interested in the behavior of the system far from the bifurcation point. Indeed, continuing all the branches for many parameter values easily becomes computationally unbearable, especially when employing \textit{full-order} methodologies, or when forced to take only small exploration steps.
Several efforts have been dedicated to reducing the computational complexity associated with the reconstruction of the bifurcation diagram and the detection of the bifurcation points, e.g.\ via the development of ad-hoc surrogate models \cite{PichiPhd}. Despite this, a complete and computationally affordable investigation of bifurcation problems is still far from being reached. 

% Bifurcation problems are especially significant in engineering contexts, as parameter fluctuations may cause the system to enter a bifurcating regime of the physical model. 
% There, even small parameter perturbations may trigger qualitatively different solutions, underscoring the importance of identifying these regimes for controlling.

{Bifurcations have been identified also in the context of stochastic PDEs \cite{Arnold1992, CrisovanModelOrderReduction2019, VenturiWeightedReducedOrder2019a, CarereWeightedPODreductionApproach2021}, which allow in principle a more realistic modeling of physical systems.} 
Indeed, by considering the parameters characterizing the physics of the system as a source of uncertainty (boundary conditions, physical parameters, and geometry), the stochastic PDEs inherently incorporate the characteristic randomness typical of real-world scenarios. 
The importance of such a research line is confirmed by the recent work \cite{KuehnUncertaintyQuantificationAnalysis2024}, where the authors investigate the statistical properties of the \textit{stochastic bifurcation} points for the Allen-Cahn equation.

{Although the focus of this work is on deterministic bifurcation detection, we pursue this goal by analyzing stochastic PDEs obtained by parameter randomization.
In particular, using stochastic numerical techniques, we show how the perturbed system characterizes the deterministic bifurcation branches of the original PDE.}
%{Thus, the core idea is to infer information about bifurcations in the deterministic PDE by numerically solving a correspondent stochastic PDE, obtained by modeling the bifurcation parameter as a random variable (or as a stochastic process in the non-scalar case).}

Specifically, we solve the stochastically-perturbed problem by approximating the unknowns with Polynomial Chaos (PC) expansion \cite{Sudret2008, Xiu2003}. 
PC has become increasingly prominent in stochastic computations, {since it allows for an efficient characterization of the solution probability distribution. 
We argue that PC solutions provide not only statistical information on the occurrence of stochastic bifurcations, as demonstrated in previous works \cite{Venturi2010, Breden2023}, but also meaningful approximations of the corresponding deterministic branches.

To this end, we present a systematic approach for inferring the bifurcation branches of a deterministic PDE directly from the PC solution of the associated perturbed PDE, revealing that possibility to simultaneously capture all deterministic branches within a single run of the numerical solver of the stochastic problem, without requiring any a-priori knowledge or ad-hoc initialization.} 

\begin{figure}
    \centering
    \includegraphics[width=\textwidth]{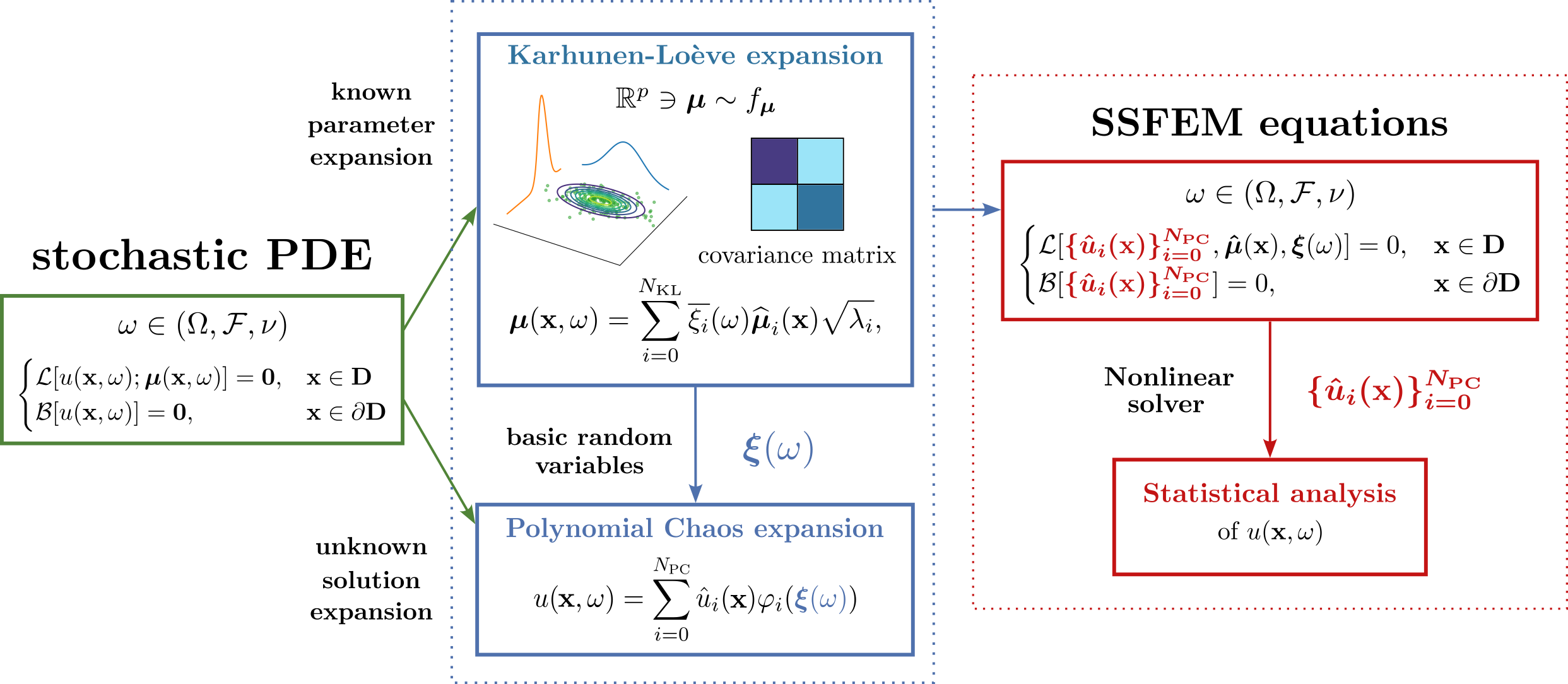}
    \caption{Statistical analysis of bifurcations via the stochastic-perturbation approach.}
    \label{fig:scheme}
\end{figure}

More precisely, to numerically solve the stochastically-perturbed PDE, we employ the Spectral Stochastic Finite Element Method (SSFEM) \cite{Stefanou2009}, which extends in a probabilistic sense the well-known Galerkin projection approach used in FEM, by choosing the PC expansion to approximate the stochastic unknowns. 
It works as follows: given a stochastic PDE, the known stochastic parameters are approximated using the Karhunen–Loève (K–L) expansion, which represents them in terms of a finite set of underlying seed random variables $\overline{\xi}$. 
These variables are then employed in the PC expansion of the unknown stochastic solution. The resulting formulation yields a deterministic system whose unknowns are the coefficients of the PC representation and which can be solved using standard numerical solvers.
A schematic overview of the pipeline is illustrated in Figure \ref{fig:scheme}.
{As we will show, the same workflow applies to bifurcating ordinary differential equations (ODEs), with the only difference that the coefficients of the PC expansion of the unknown no longer depend on the spatial variable.}

In this work, we {focus on pitchfork bifurcations, starting with their normal form, i.e., a simplified scalar ODE that reproduce the qualitative features of bifurcation behavior, and subsequently considering a more complex evidence of such phenomenon in computational fluid dynamics (CFD), namely the Coand\u{a} effect} \cite{QuainiSymmetryBreakingPreliminary2016}. 
We present computational evidence, by assessing the method capabilities first on the supercritical pitchfork normal form, and then demonstrating the effective recovery of the bifurcation diagram in a more challenging model: the Navier–Stokes equations.

We first introduce in Section \ref{section:stochastic_expansion} the two expansions used in the SSFEM setting, through which we can represent a generic second-order stochastic process in a finite-dimensional space: the K-L and the PC expansions. 
Moreover, we discuss non-intrusive sample-based and intrusive methods to retrieve PC coefficients in the bifurcating setting. 
Successively, we proceed with the description of the intrusive approach in subsection \ref{section:SSFEM}, {first applying it to the investigation of the bifurcation behavior of an ODE, namely, the supercritical pitchfork normal form in Section~\ref{sec:pitchfork}, and subsequently to the Navier-Stokes equations in the context of the Coand\u{a} effect in Section~\ref{section:testcase}.
} 
Successively, in Section \ref{sec:results} we show the numerical results for the Coand\u{a} effect, underlying the relation discovered between the PC solution and the deterministic bifurcation diagram. 
In particular, we investigate the relationship between the PC coefficients and the computational mesh, as well as the impact of varying the probability distribution assigned to the viscosity parameter on the resulting solutions.
Finally, we conclude by showing a reproduction of the deterministic bifurcation diagram, which we obtain only relying on the information collected through the stochastic perturbation approach.
\subsection{Related works} \label{sec:related_works}

{This section provides an overview of related literature on the numerical treatment of deterministic bifurcating systems and on the theory of stochastic bifurcations, including numerical methods for their solution. 
Indeed, despite not specifically focusing on stochastic bifurcations themselves, we leverage methodologies from such context to obtain information concerning deterministic bifurcations, creating a novel connection between the two settings.}
% We emphasize that our focus is not on stochastic bifurcations themselves, but on leveraging methodologies from that setting to extract information about deterministic bifurcations, which, to the best of our knowledge, has not been already investigated in previous works.}

%The study of bifurcating phenomena in nonlinear PDEs has been a subject of extensive research due to its relevance in various physical and engineering domains \cite{ambrosetti1995primer}.
For what concerns the numerical solution of deterministic bifurcating systems, traditional numerical methods \cite{CalozNumericalAnalysisNonlinear1997}, including the Finite Element Method and Spectral Element Method, has been widely employed to explore the coexisting solution branches and the critical regime in a neighborhood of the bifurcation points \cite{Cliffe2000, Pintore2020}. 
Additionally, to effectively explore the global behavior of the system, numerical methods such as continuation techniques \cite{allgower2003, PichiStrazzullo2022} and deflation methods \cite{FarrellDeflationTechniquesFinding2015} have been designed, respectively to reconstruct a specific solution branch and discover new admissible ones as parameters vary.
{Additional strategies for the analysis of bifurcation problems include the study of their normal forms \cite{BruntonDiscoveringGoverningEquations2016,KuznetsovElementsAppliedBifurcation2023} (see Section \ref{sec:pitchfork} for the example on the pitchfork case), or spectral submanifolds \cite{HallerNonlinearNormalModes2016,BettiniDatadrivenNonlinearModel2025} for reducing dynamical systems to low-dimensional smooth invariant manifolds.}

Nevertheless, the computational effort for a complete exploration of the solutions in the whole parameter space can easily become unbearable. 
Recent works focused on developing reduced-order models (ROMs) \cite{hesthaven2015certified,benner2017model} to enable many-query simulations, and thus efficiently recovering the full bifurcating pattern. 
These approaches can be divided into two categories: the ones exploiting standard projection-based techniques \cite{PichiPhd,Pitton_2017,Khamlich2022,Pintore2020}, and the data-driven ones, often based on machine-learning strategies \cite{BruntonDiscoveringGoverningEquations2016,PichiArtificialNeuralNetwork2023,PichiGraphConvolutionalAutoencoder2024,coscia2024generative}.

All the above works treat the bifurcation problem in the absence of randomness. Indeed, the deterministic assumption simplifies the analysis of the system's behavior, not requiring the additional computational burden typical of stochastic settings, but it is often not physically realistic. For this reason, the presence of stochasticity in the models cannot be disregarded for a comprehensive analysis of the topic, as bifurcating phenomena depend on various uncertainties due to input variability, parameter uncertainty, and numerical discretization \cite{uq-book}.

Towards a more realistic study of the problem, the concept of bifurcation has been expanded in the stochastic setting with the analysis of additive and multiplicative noise effects on the models, thus characterizing various types of \textit{stochastic bifurcations}
\cite{Arnold1992}. Indeed, while in the deterministic case the bifurcation corresponds to a qualitative change in the dynamics of the system as a response to a small variation of the parameters, in the stochastic case it can be observed as a sudden change in some invariant measure of the system's solution as a result of some imposed additive or multiplicative perturbation. Such changes can involve variations in the largest Lyapunov exponent (\textit{dynamical bifurcations}), as well as changes in the shape of the solution probability density function (\textit{phenomenological bifurcations}) \cite{arnold2010}.

Early investigations explored the dynamic response, stability, and bifurcating behavior of nonlinear dynamical systems under stochastic excitation. In \cite{SriNamachchivaya1990}, the author studied the impact of small stochastic perturbations on systems exhibiting codimension-one and -two bifurcations. The work exploited methods of stochastic averaging and stochastic normal forms to analyze the asymptotic behavior of nonlinear dynamical systems in the presence of noise. Furthermore, various works have been conducted for specific systems of equations, such as the stochastic reaction-diffusion equation perturbed by an infinite-dimensional Wiener process in \cite{Blumenthal2023}, and the stochastic regime-switching predator-prey system in \cite{Liu2024}. 
%Moreover, investigating the interplay between deterministic and stochastic bifurcations, the authors in \cite{Mirzakhalili2019} propose a method to derive probabilistic bifurcation diagrams for stochastic nonlinear dynamical systems using the Fokker-Planck equation. 

Going in the direction of identifying the stochastic perturbation applied to the system with the randomness associated with its parameters, in \cite{Venturi2010}, the authors explore a stochastic approach to the known deterministic bifurcations in Rayleigh–Bénard convection within a two-dimensional cavity. 
The study investigates the onset of instability and the existence of multiple stable states under specific ranges of Rayleigh number. 
Indeed, through stochastic simulations conducted around a point of deterministic bifurcation, the authors reveal the influence of also random initial flow states on convection patterns {through statistical considerations on the PC solution of the stochastically-perturbed PDE}. 
The importance of stochastic aspects has been highlighted via Monte Carlo simulations in \cite{Rubinstein2016}, investigating the probabilities of the solution branches through a stochastic sampling of the parameterized initial conditions. Moreover, they introduce a comparison between the accuracy of statistical moments computation of Monte Carlo sampling and the PC expansion of the solution, confirming the high reliability and computational convenience of the latter.
\section{Stochastic modeling in PDEs}\label{section:stochastic_expansion}

Given an established deterministic model for a physical system, we want to define a convenient setup for its stochastic counterpart, to study the effect of uncertainty in the input parameters on the solution itself.

As illustrated in Figure \ref{fig:scheme}, let us consider the stochastic PDE of the form: 
{
\begin{equation} \label{eq:stochasticPDE}
    \begin{cases}
    \vectspace{L}[u(\vect{x},\omega); \vect{\mu}(\omega)] = 0 &\quad \vect{x}\in\vectspace{D}, \quad \omega\in\Omega, \\
    \vectspace{B}[u(\vect{x},\omega)] = 0 &\quad \vect{x}\in\partial\vectspace{D},\quad \omega\in\Omega, 
\end{cases}
\end{equation}
}
where {$\vect{\mu}(\omega):\Omega\rightarrow\set{R}^p$ are the $p$ randomized parameters},  $\vectspace{L}$ and $\vectspace{B}$ incorporate respectively the physics and the boundary constraints, $u(\vect{x},\omega):\vectspace{D}\times\Omega\rightarrow\set{R}$ is the unknown solution of the system, $\vectspace{D}\subset\set{R}^n$ the spatial domain, and $(\Omega,\vectspace{F},\nu)$ a probability space, where $(\Omega, \mathcal{F})$ is a measurable space with $\mathcal{F}$ a $\sigma-$algebra defined on $\Omega$, and $\nu$ a probability measure on it.
We recall that a generic real stochastic process is defined as a family $\vect{\alpha}(\vect{x},\omega)=\{\vect{\alpha}(\vect{x},\cdot)\}_{\vect{x}\in \vectspace{D}}$ of random variables, which are in turn measurable functions of the type $\vect{\alpha}(\vect{x},\cdot):(\Omega,\vectspace{F},\nu)\rightarrow (\set{R}^d,\vectspace{B})$ for a generic $d\leq\infty$.

By considering the solution as a generic stochastic process $u(\vect{x},\omega)$, to enable numerical computations we need to reduce the inﬁnite-dimensional probability space to a ﬁnite-dimensional one. This is accomplished by exploiting a decomposition of the stochastic process via a finite set of random variables.

We are interested in obtaining meaningful expansions of the real stochastic processes having a finite second-order moment, thus belonging to the space:
$$
    \set{L}_2(\Omega, \vectspace{F}, \nu; \set{R}^d):= \left\{u(\vect{x},\omega):\vectspace{D}\times\Omega\rightarrow\set{R}^d\ \ \text{s.t.}\ \ \int_{\Omega} u(\vect{x},\omega) ^2d\nu(\omega)<\infty\right\}.
$$
{To this end, we first introduce the PC expansion, which constitute the basis for the generalization of FEM in the probabilistic context, and then present the complete SSFEM framework.}

\subsection{Polynomial Chaos expansion} \label{section:gPC_expansion}
The PC expansion decomposes a general stochastic
process $u(\vect{x}, \omega)\in\set{L}_2(\Omega, \vectspace{F}, \nu; \set{R}^d)$
thanks to the Cameron-Martin theorem \cite{CameronMartin1947}, which establishes the existence of a sequence of closed, pairwise orthogonal linear subspaces of $\set{L}_2(\Omega, \vectspace{F}, \nu; \set{R}^d)$, spanned by polynomial bases. 
In particular, the PC expansion expresses the process $u(\mathbf{x}, \omega)$ as the infinite sum:
\begin{equation}\label{eq:PC}
    u(\mathbf{x}, \omega) = \sum_{i=0}^{\infty}
    \widehat{u}_{i}(\mathbf{x}) \psi_i(\{\zeta_j(\omega)\}_{j=1}^{N_{\text{RV}}}), \quad \text{with} \quad 
    \int_{\set{R}^d} \psi_i(\vect{\zeta})\psi_j(\vect{\zeta})d\nu_{\vect{\zeta}} \ = \ \norm{\psi_i}^2\delta_{i,j},
\end{equation}
where $\widehat{u}_{i}(\mathbf{x})$ are the expansion coefficients associated with the stochastic polynomial basis defined as $\psi_i(\{\zeta_j(\omega)\}_{j=1}^{N_{\text{RV}}}):(\Omega,\vectspace{F},\nu)\rightarrow\set{R}$, which are orthogonal with respect to the Gaussian measure $\nu_{\vect{\zeta}}$ and depend on $N_{RV}$ independent standard normal random variables $\zeta_j(\omega)$.  
In the following, we will refer to such set of random variables, called also \textit{seeds} or \textit{basic random variables}, as $\vect{\zeta}=\{\zeta_j(\omega)\}_{j=1}^{N_\text{RV}}$. 
As an exaple, for a Gaussian measure $\nu_{\vect{\zeta}}$, the orthogonal basis $\bm{\psi}$ corresponds to the set of Hermite polynomials \cite{GhanemSpanos1990}.

For computational purposes, the infinite sum in \eqref{eq:PC} is in practice truncated to a finite number of terms $N_{PC}$. It is often useful to relate such a number to the desired maximum degree $M$ of the polynomials, and the number of variables $N_{RV}$ from which they depend as
$1+N_{PC} = \frac{(M+N_{RV})!}{M!N_{RV}!}.$

The generalization of the Polynomial Chaos (PC) expansion, originally developed for Gaussian measures, was introduced in \cite{XIU2002_gpc} under the name generalized Polynomial Chaos (gPC). 
In this framework, the expansion is constructed using different families of orthogonal polynomials,referred to as Wiener-Askey polynomials \cite{ernst2012},chosen according to the probability measure (including non-Gaussian ones) associated with the underlying random variables $\nu_{\vect{\zeta}}$. 
The correspondence between probability measures and polynomial families is summarized in Table~\ref{tab:askey_scheme}.

\newcommand{\BB}[1]{\cellcolor{blue!#1}}
\newcommand{\GG}[1]{\cellcolor{green!#1}}

\begin{table}[t]
\centering
\resizebox{\textwidth}{!}{
\begin{tabular}{c|ccc}
\rowcolor[HTML]{C1CDCD}
\toprule
\textbf{Type} & \textbf{Basic random variables} & \textbf{Wiener-Askey polynomials} & \textbf{Support} \\
\midrule
\multirow{4}{*}{\textbf{Continuous}} &\BB{20} Gaussian & \BB{20}Hermite &\BB{20}$(-\infty,\infty)$ \\
&\BB{10}Gamma &\BB{10}Laguerre & \BB{10}$[0,\infty)$ \\
&\BB{20}Beta & \BB{20}Jacobi & \BB{20}$[a,b]$ \\
& \BB{10}Uniform & \BB{10}Legendre & \BB{10}$[a,b]$ \\
\midrule
\multirow{3}{*}{\textbf{Discrete}} &\GG{20}Poisson &\GG{20}Charlier &\GG{20}$\{0,1,2,\dots\}$ \\
&\GG{10}Binomial &\GG{10}Krautchouk &\GG{10}$\{0,1,2,\dots,N\}$ \\
&\GG{20}Negative Binomial &\GG{20}Meixner &\GG{20}$\{0,1,2,\dots\}$ \\
\bottomrule
\end{tabular}}
\caption{Some of the correspondences between basic random variables, the Wiener--Askey orthogonal polynomials for gPC expansion, and their support.}
\label{tab:askey_scheme}
\end{table}

\subsubsection{Statistical moments computation} \label{sec_PC_moments}
A key advantage of the PC representation is that the solution coefficients $\widehat{\vect{u}}(\bm{x})$ allow for a straightforward computation of the first and second-order moments of the full solution $u(\bm{x},\omega)$.
Indeed, leveraging the orthogonality property of the polynomials $\{\psi_i\}_{i=1}^{N_\text{PC}}$ and recognizing that, since the polynomial of degree zero is constant, the following holds:
\begin{equation} \label{eq:GPC_mean}
    \mathbb{E}[\psi_i(\vect{\xi})] = \mathbb{E}[\psi_i(\vect{\xi})\psi_0(\vect{\xi})] = \delta_{i,0}, \quad \text{and} \quad 
 \set{E}\left[\sum_{i=0}^{N_\text{PC}}
     \widehat{u}_i(\vect{x})\psi_i(\vect{\xi})\right] = \widehat{u_0}(\vect{x}),
\end{equation}
which states that the mean of the solution is equivalent to the first coefficient of the expansion.

Furthermore, each element of the covariance matrix can be approximated by a weighted sum involving all the coefficients except the first one $\widehat{u_0}$, where the weights correspond to the polynomial stochastic norms $\set{E}[\vect{\psi}_i^2]$:
\begin{equation} \label{eq:GPC_covariance}
    \begin{split}
        \mathrm{Cov}\left[u(\vect{x},\omega), u(\vect{y},\omega) \right] &= \set{E}\left[\left(\sum_{i=0}^{N_\text{PC}}\widehat{u_i}(\vect{x})\psi_i(\vect{\xi})-\widehat{u_0}(\vect{x})\right)\left(\sum_{i=0}^{N_\text{PC}}\widehat{u_i}(\vect{y})\psi_i(\vect{\xi})-\widehat{u_0}(\vect{y})\right)\right] \\
        &=\sum_{i=1}^{N_\text{PC}} \widehat{u_i}(\vect{x})\widehat{u_i}(\vect{y})\set{E}[\vect{\psi}_i^2] .
    \end{split}
\end{equation}
The straightforward computation of the statistical moments from the PC surrogate model represents a key advantage. 
Indeed, both the PC statistical analysis and PDF reconstruction are facilitated, since the first can now be performed through a proper combination of the coefficients, and the latter by sampling the solution polynomial expansion, evaluating the polynomials at a sample of their basic random variables $\vect{\xi}$.

\subsubsection{Intrusive and non-intrusive computations of the coefficients}\label{section:intrusive_vs_nonintrusive}

%In the former sections, we have seen that modeling the stochastic components of the system as random variables and stochastic processes allow us to express them as series expansions. 
A crucial task in this setting is the determination of the PC coefficients. 
%Indeed, once found their representations, such expansions could also be used to infer relevant information, such as statistical moments.
Their computation can be carried out using either \textit{intrusive} or \textit{non-intrusive} methods.

While the former class of methodologies requires some knowledge about the PDE, approaches belonging to the latter class typically treat the model as a black-box, from which samples of the solution can be extracted and exploited to retrieve the coefficients (for instance by interpolation).
%In general, when well-established deterministic codes exist, sample-based non-intrusive stochastic simulations of the system can be performed, computing the solution for different parameter values \cite{babuvska2007stochastic}. 
Such methods include Monte Carlo simulation \cite{Rubinstein2016}, orthogonal spectral projection \cite{xiu2010}, and stochastic collocation \cite{babuvska2007stochastic}. 
In particular, Monte Carlo simulation gathers a set of solutions for a given set of input parameters and deduces statistical properties from the collected data. 
On the other hand, orthogonal spectral projection computes the coefficients through numerical integration of projection integrals on the polynomials, while stochastic collocation achieves this via interpolation of sampled deterministic solutions (we refer the reader to \cite{sullivan2015}).

{
In the specific context of bifurcation problems, non-intrusive, sample-based approaches may be unsuitable. Indeed, obtaining a comprehensive and representative set of solution samples for a given parameter set is non-trivial, as it would require a priori knowledge of the bifurcation structure (e.g., the number of branches) and multiple runs of the same deterministic model for identical parameter values, each exploiting a different ad hoc initial conditions. As an explicit example supporting this claim, we refer to Appendix \ref{sec:MC}, which presents a Monte Carlo non-intrusive approach applied to the perturbed bifurcating system and highlights the associated limitations.
}

Therefore, we proceed with the presentation of the SSFEM, an intrusive approach for the retrieval of the expansion coefficients that does not rely on the availability of solution samples of stochastic PDEs.

%While the SSFEM approach also requires initialization for solving the system of gPC coefficients, it explores the solution space differently than deterministic solvers. Rather than trying to find specific solutions, it aims to characterize the probability distribution of solutions under parameter perturbations. This probabilistic perspective can reveal multiple solution branches through the structure of the resulting distributions, without choosing ad-hoc initializations (see Section \ref{sec:pitchfork}).

\subsection{Spectral Stochastic Finite Element Method} \label{section:SSFEM}
{In this section, we describe the SSFEM \cite{Stefanou2009, GhanemSpanos1990, xiu2010} pipeline: an intrusive methodology for computing the coefficients of the PC expansion of a stochastic PDE solution.}

Let $\mathcal{L}$ denote the operator defined in \eqref{eq:stochasticPDE}. We consider the following abstract problem:
\begin{equation} \label{eq:general_prob}
\mathcal{L}[u(\bm{x},\omega);\bm{\mu}(\omega)]:=(\vect{D}(\vect{x}) + \vect{S}(\vect{x}, \vect{\mu}(\omega)))[u(\vect{x}, \omega)] = \vect{f}(\vect{x})\ \quad \text{for}\ \vect{x} \in \vectspace{D},\ \omega \in \Omega
\end{equation}
where $u(\vect{x}, \omega) : \vectspace{D} \times \Omega \rightarrow \mathbb{R}$ is the unknown function, and $\vect{\mu}(\omega) : \Omega \rightarrow \mathbb{R}^p$ represents the $p$ {random variables obtained from stochastic perturbation of the deterministic system's parameters}. The deterministic operator $\vect{D}$ and its stochastic counterpart $\vect{S}$ act on $u(\vect{x}, \omega)$ returning an object that belongs to $\set{R}^{out}$, while $\vect{f}(\vect{x}): \vectspace{D} \rightarrow \mathbb{R}^{out}$ represents the source term\footnote{In the following, we will omit the full notation with the dependence on ($\vect{x}, \omega$) for $u$, $\vect{\mu}$ and $\vect{f}$.}. With stochastic and deterministic operators, we intend operators that respectively depend or not from the stochastic variables we inserted in the system, namely $\vect{\mu}$.

\begin{figure}[tb]
    \centering
    \resizebox{\textwidth}{!}{%
    \begin{tikzpicture}[
        >=Stealth,
        node distance=0.5cm,
        box/.style={rectangle, draw, text width=2cm, align=center,
        font=\footnotesize, rounded corners, minimum height=1.3cm, line width=0.7mm},
        arrow/.style={->, line width=1pt},
        ]
        % Nodes
        \node[box, draw=ssfem_green] (pdf) {Known covariance of $\vect{\mu}$};
        \node[box, draw=ssfem_blue, right=of pdf] (KLPC) {K-L expansion of $\vect{\mu}$ and PC expansion of $u$};
        \node[box, draw=ssfem_red, right=of KLPC] (SSFEM) {Galerkin projection};
        \node[box, draw=ssfem_red, right=of SSFEM] (solve) {PC coefficients retrival};
        \node[box, draw=ssfem_red, right=of solve] (stats) {PDF reconstruction of $u$};
        % Arrows
        \draw [arrow, draw=ssfem_green] (pdf) -- (KLPC);
        \draw [arrow, draw=ssfem_blue] (KLPC) -- (SSFEM);
        \draw [arrow, draw=ssfem_red] (SSFEM) -- (solve);
        \draw [arrow, draw=ssfem_red] (solve) -- (stats);
    \end{tikzpicture}}
    \caption{Synthetic workflow diagram showing the process of computing the surrogate
    PDF of $u$ by solving Equation \eqref{eq:general_prob} using SSFEM pipeline.}
    \label{fig:workflow}
\end{figure}

The general workflow, schematically illustrated in Figure~\ref{fig:workflow}, begins by prescribing a covariance structure for $\bm{\mu}$. 
This covariance fully characterizes the probability distribution assigned to the parameters of the original deterministic system during the randomization step.
%Therefore, the stochastic parameters are expanded using the K-L decomposition presented in Appendix \ref{section:KL_expansion}. 
{In our setting, we employ a simplified Karhunen--Loève (K--L) expansion, whose general formulation tailored to stochastic processes with finite second-order moments is reported in Appendix \ref{section:KL_expansion}. 
Here, we consider a random vector $\bm{\mu}\in\mathbb{R}^p$ with no spatial dependence, consisting of $p$ random variables.
Thus, in our case the truncated K--L expansion reads as
\begin{align*}
\bm{\mu}(\omega)=\mathbb{E}[\bm{\mu}]+\sum_{i=1}^{N_{\mathrm{KL}}}\bar{\xi}_i(\omega)\,\widehat{\bm{\mu}}_i\,\sqrt{\lambda_i},
\end{align*}
where $N_{\mathrm{KL}}\le p$, and $\{\lambda_i,\widehat{\bm{\mu}}_i\}$ are the eigenpairs of the covariance matrix
\begin{align*}
(K_{\bm{\mu}})_{ij}:=\mathbb{E}\!\left[(\mu_i-\mathbb{E}[\mu_i])(\mu_j-\mathbb{E}[\mu_j])\right],
\end{align*}
and the random variables $\bar{\xi}_i$ satisfy
$\mathbb{E}[\bar{\xi}_i]=0$ and
$\mathbb{E}[\bar{\xi}_i\bar{\xi}_j]=\delta_{ij}$.
Not that, for $p=1$, the K--L expansion reduces to a simple normalization of the random
variable,
\begin{align*}
\mu(\omega)=\mathbb{E}[\mu]+\bar{\xi}(\omega)\sqrt{\mathbb{E}\!\left[(\mu-\mathbb{E}[\mu])^2\right]},
\end{align*}
with $\mathbb{E}[\bar{\xi}]=0$ and $\mathbb{E}[\bar{\xi}^2]=1$.
}

Since a straightforward K-L expansion of the solution $u$ is not feasible as its covariance function is not known a-priori, the finite-dimensional representation of $u$ is achieved via the truncated PC expansion:
\begin{equation}\label{eq:GPC}
    u(\vect{x},\omega) = \sum_{i=0}^{N_\text{PC}} \widehat{u}_{i}(\vect{x})\psi_i(\{\overline{\xi_l}(\omega)\}_{l=1}^{N_\text{KL}}),
\end{equation} 
where we identify the \textit{seed random variables} of the PC expansion with the uncorrelated normalized random variables\footnote{The index for KL random variables starts from $1$ due to the assumptions in \ref{section:KL_expansion}.} resulting from the K-L expansion of the parameters $\vect{\xi}(\omega):=\{\overline{\xi_j}(\omega)\}_{j=1}^{N_\text{KL}}$.

{The SSFEM pipeline then proceeds by approximating the spatial PC coefficients
$\widehat{\boldsymbol{u}}(\boldsymbol{x})$, along with any other spatially dependent
quantity, using a Finite Element basis. A Galerkin projection is subsequently
applied in both the spatial domain $\mathcal{D}$ and the probability space, with
the PC polynomials employed as test functions in the latter.

For clarity, we present here an illustrative example for the case $p=1$ and an affine
stochastic operator $\boldsymbol{S}(\mu,\boldsymbol{x})=\mu\,\boldsymbol{S}(\boldsymbol{x})$,
which is the setting adopted in all numerical experiments in this work, but does not exhaust the range of applications of SSFEM.}
In this illustrative case, substituting the PC expansion of $u$ and the K-L expansions of $\mu$, the problem can be recast as follows:
\begin{equation}\label{eq:general_prob_expansion}
    \left(\vect{D}(\vect{x})+\sum_{i=0}^{N_\text{KL}}\overline{\xi_i}(\omega)
    \widehat{\mu}_{i}(\vect{x})\sqrt{\lambda_i}\vect{S}(\vect{x})\right)\left[\sum_{j=0}^{N_\text{PC}}
\widehat{u}_{j}(\vect{x})\psi_j(\vect{\xi}(\omega))\right] =
\vect{f}(\vect{x}).
\end{equation}
We then exploit standard FEM to numerically approximate \eqref{eq:general_prob_expansion} thanks to a set of deterministic basis functions $\{\phi_k(\vect{x})\}_{k=1}^{N_\text{D}}$ such that
\begin{equation}
\widehat{u}_{j}(\vect{x})
=\sum_{k=1}^{N_D}\widehat{u}_{k,j}\phi_k(\vect{x}), \quad \forall j\in\{0,\dots,N_\text{PC}\}.
\end{equation}
Substituting the former expansion into Equation \eqref{eq:general_prob_expansion}, we obtain:
\begin{equation}\label{eq:general_prob_discretized}
    \left(\vect{D}(\vect{x})+\sum_{i=0}^{N_\text{KL}}\overline{\xi_i}(\omega)\widehat{\mu}_{i}(\vect{x})\sqrt{\lambda_i}\vect{S}(\vect{x})\right)\left[\sum_{j=0}^{N_\text{PC}} \sum_{k=1}^{N_D}\widehat{u}_{k,j}\phi_k(\vect{x})\psi_j(\vect{\xi}(\omega))\right] = \vect{f}(\vect{x}),
\end{equation}
that has to be solved for the unknown $\mathbf{\widehat{U}}={\{\widehat{u}_{k,j}\}_{k=1}^{N_{D}}}{\vphantom{=\{\widehat{u}_{k,j}\}}}_{j=0}^{N_\text{PC}} \in \mathbb{R}^{N_\text{D} \times(N_{\text{PC}}+1)}$.

Finally, to construct an algebraic system of equations for $\mathbf{\widehat{U}}$, one imposes the orthogonality of the residual in Equation \eqref{eq:general_prob_discretized} with respect to the deterministic and stochastic test functions, respectively  $\{\phi_n(\vect{x})\}_{n=1}^{N_\text{KL}}$ and the polynomials $\{\psi_m(\vect{\xi})\}_{m=0}^{N_\text{PC}}$.

For instance, in the case of linear operators $\bm{D}(\bm{x})$ and $\bm{S}(\bm{x})$ acting on $u(\bm{x},\omega)$, the discrete counterpart of the stochastic PDE reads as the following system of equations $\forall n\in\{1,\dots,N_D\}$, and $\forall m\in\{0,\dots,N_\text{PC}\}$:
\begin{equation} \label{eq:SSFEM_final_eq}
    \sum_{j=0}^{N_\text{PC}} \sum_{k=1}^{N_D}C_{j,m}\widehat{u}_{k,j}A_{k,n}+\sum_{i=0}^{N_\text{KL}}\sum_{j=0}^{N_\text{PC}}\sum_{k=1}^{N_D}E_{i,j,m}\widehat{u}_{k,j}B_{i,k,n}
    = F_n,
\end{equation}
having defined the tensors

    \begin{equation} \label{eq:SSFEM_mats_general}
    \begin{gathered}
            A_{k,n} = \int_{\vectspace{D}}\vect{D}(\vect{x})[\phi_k(\vect{x})]\phi_n(\vect{x})d\vect{x}, \quad
            B_{i,k,n} = \int_{\vectspace{D}}\widehat{\mu_i}(\vect{x})\vect{S}(\vect{x})[\phi_k(\vect{x})]\phi_n(\vect{x})d\vect{x}, \\
            F_n = \int_{\vectspace{D}}\vect{f}(\vect{x})\phi_n(\vect{x})d\vect{x}, \quad
            C_{j,m} = \int_{\set{R}^d}\psi_j(\vect{\xi})\psi_m(\vect{\xi})d\nu_{\vect{\xi}},  \\
            E_{i,j,m} = \int_{\set{R}^d}\overline{\xi_i}\sqrt{\lambda_i}\psi_j(\vect{\xi})\psi_m(\vect{\xi})d\nu_{\vect{\xi}}.  
    \end{gathered}
    \end{equation}

Therefore, the system can now be solved in function of the unknown $\widehat{\vect{U}}$ by means of a standard numerical solver. 
In the case of a generic nonlinear nature of the two operators, $\vect{D}(\bm{x})$ and $\vect{S}(\bm{x})$, a similar discrete system can be constructed and solved by a proper nonlinear solver (see Section \ref{section:testcase} for the Navier-Stokes Equations). 
%As will be underlined in the results section, the choice of the solver represents a crucial step for the methodology. 
%Indeed, when interested in simulations involving dense grids (i.e.\ high $N_{\text{D}}$), and high polynomial degrees (i.e.\ high $N_{\text{PC}}$), the convergence rate of the solver might be compromised by the number of the degrees of freedom of $\widehat{\vect{U}}$. 
%In such cases, it is convenient to perform the computations in parallel, splitting the grid between multiple processors and performing separately the residual and update step computation needed by the iterative solver.

{It is worth emphasizing that the same SSFEM pipeline, namely, approximating the unknown solution via a PC expansion built on the seed random variables arising from the K–L expansion of the randomized parameters, can also be applied to non-spatial systems, such as ordinary differential equations (ODEs).}
\section{Pitchfork bifurcation normal form}\label{sec:pitchfork}

{In this section, we apply the above intrusive pipeline to the normal form of a pitchfork bifurcation in the ODE setting. 
We exploit this toy problem to demonstrate how PC solutions of the randomized equation can provide insights into deterministic branches.}
For this preliminary analysis, we test two different distributions on the parameter, Gaussian and Uniform.

We consider the normal form of the supercritical pitchfork bifurcation given by $\dot{u} = \mu u - u^3$, whose equilibrium solutions satisfy $u(u^2-\mu)=0$. For $\mu>0$, three equilibria coexist: the trivial solution $u=0$ and the two diverging branches $u=\pm\sqrt{\mu}$. 

Applying the pipeline explained in section \ref{section:stochastic_expansion}, we obtain:
\begin{align}\label{eq:nomal_form_system}
\int_{\mathbb{R}}\sum_{i=0}^{N_{PC}}\hat{u}_i\phi_i(\bar{\xi})\Big(\Big(\sum_{j=0}^{N_{PC}}\hat{u}_j\phi_j(\bar{\xi})\Big)^2-\bar{\mu}-\sigma\bar{\xi}\Big)\phi_k(\bar{\xi})d\nu_{\bar{\xi}}=0,\quad\forall k\in\{0,\dots,N_{PC}\},
\end{align}
which can be solved using Newton method to find the scalar coefficients $\{\hat{u}_i\}_{i=0}^{N_{PC}}$.

\begin{figure}[t!]
    \captionsetup[sub][figure]{skip=-30pt,slc=off,margin={0pt,0pt}}
    
    \centering\subcaptionbox{\label{subfig:poly_big_perturbation_normal}}{\includegraphics[trim={0cm 0cm 0cm 0cm}, clip, width=0.47\textwidth]{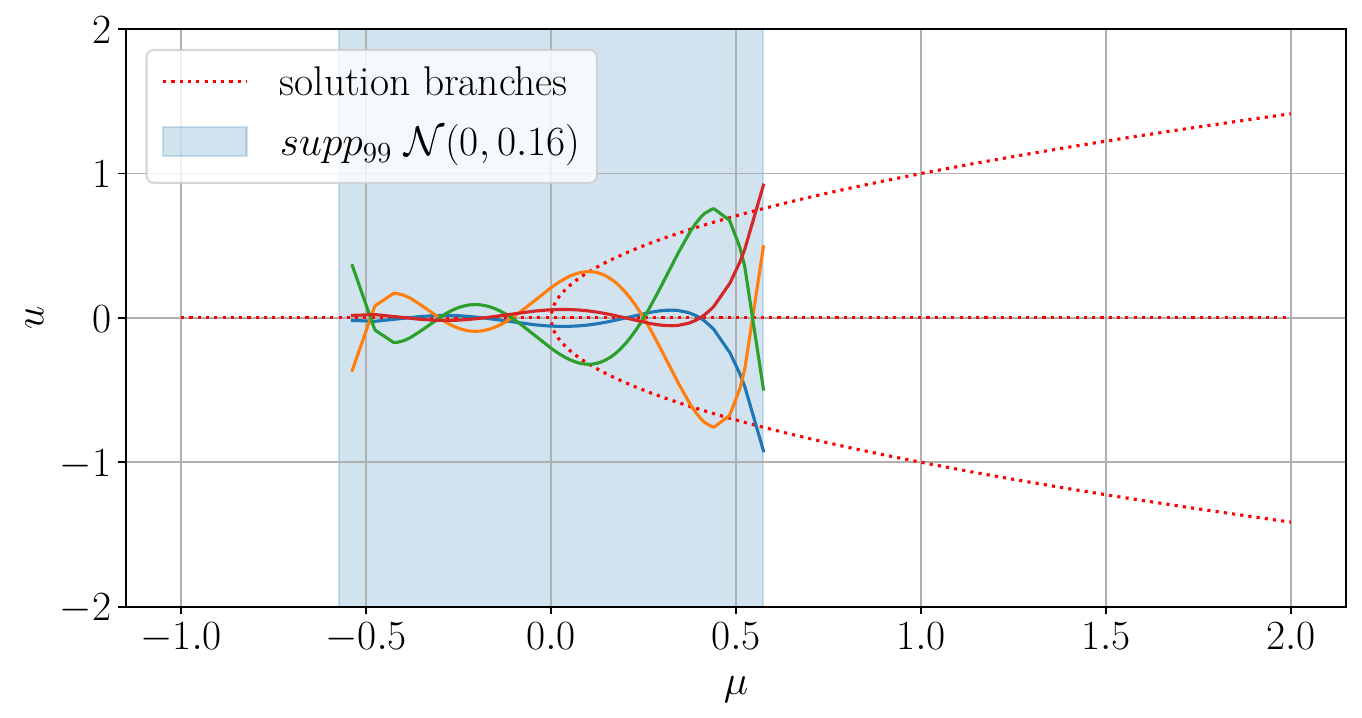}
    }
    \hfill
    \centering\subcaptionbox{\label{subfig:poly_big_perturbation_uniform}}{\includegraphics[trim={0cm 0cm 0cm 0cm}, clip, width=0.47\textwidth]{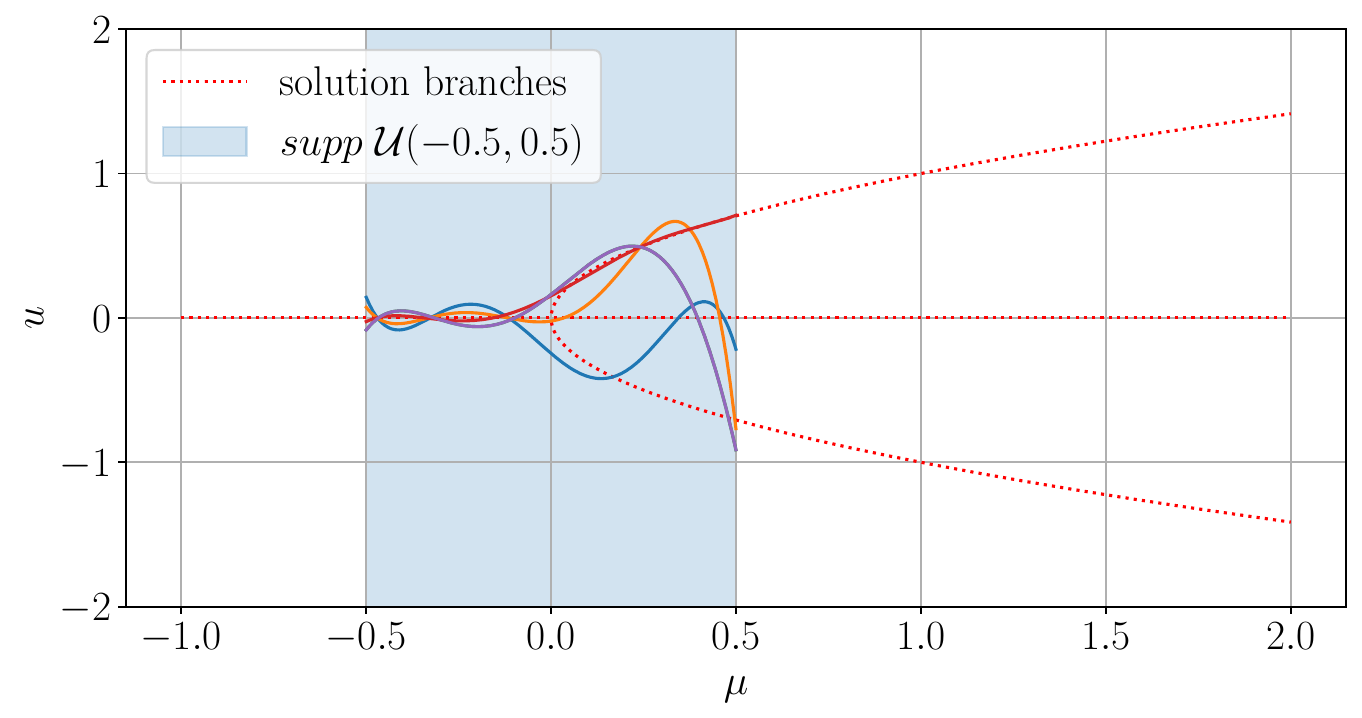}
    }
    \hfill
    \centering\subcaptionbox{\label{subfig:poly_little_perturbation_normal}}{\includegraphics[trim={0cm 0cm 0cm 0cm}, clip, width=0.47\textwidth]{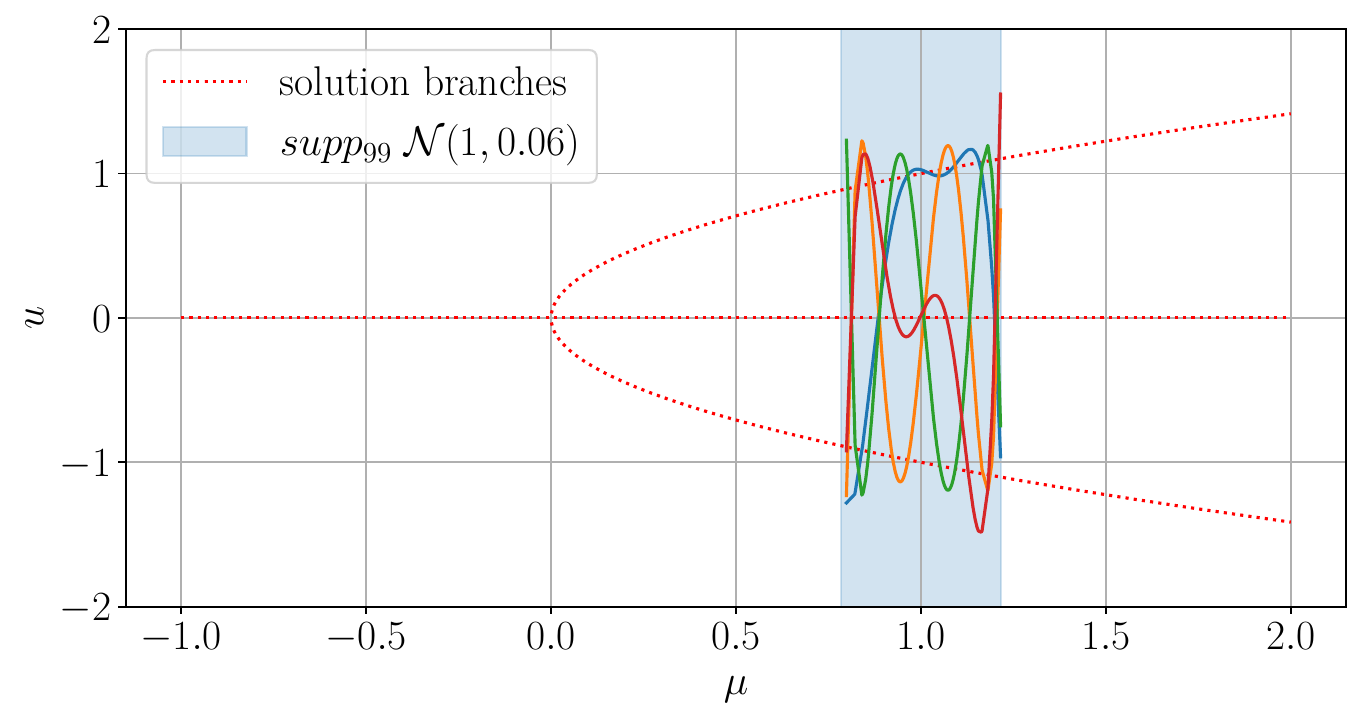}
    }
    \hfill
    \centering\subcaptionbox{\label{subfig:poly_little_perturbation_uniform}}{\includegraphics[trim={0cm 0cm 0cm 0cm}, clip, width=0.47\textwidth]{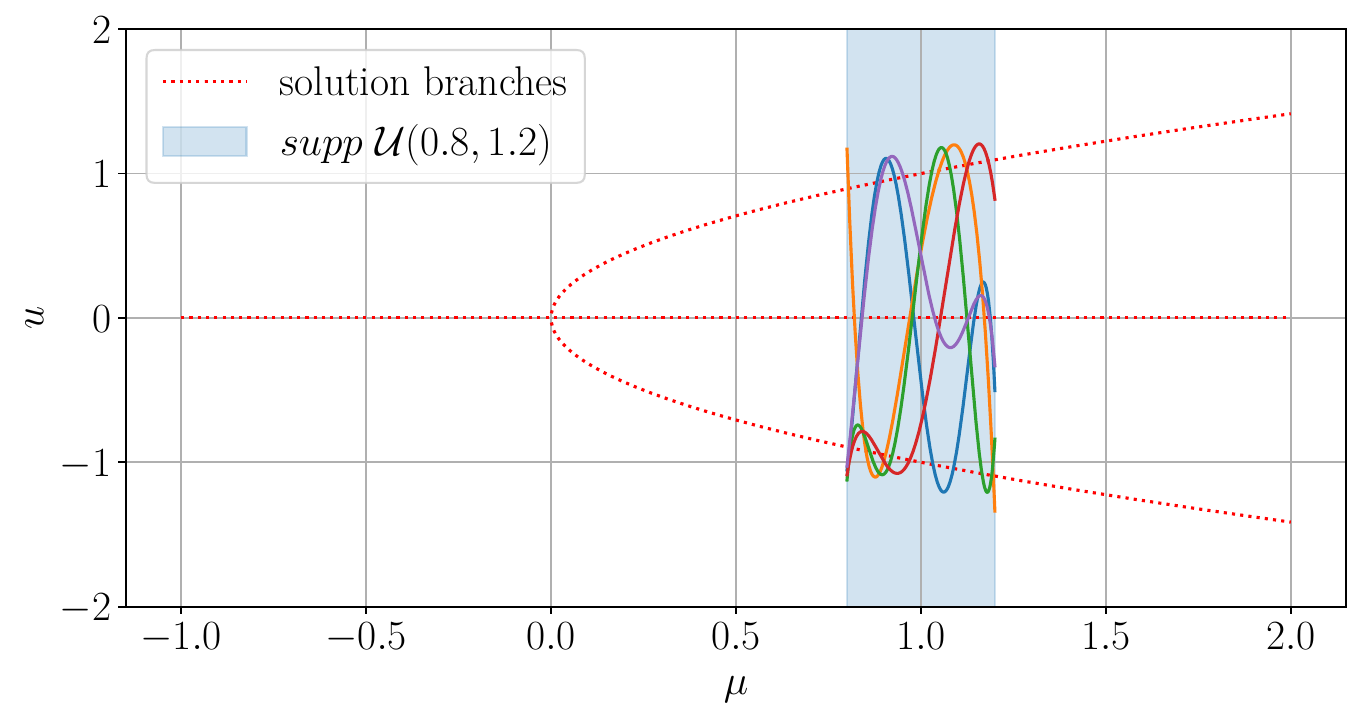}
    }
    \captionsetup{subrefformat=parens}
    
    \caption{Some admissible PC solutions with $N_{PC}=5$ for large and small perturbations of $\mu$, respectively top and bottom, considering a Gaussian and a Uniform distribution, respectively left and right. 
    The shaded area identifies the support of the distribution associated to $\mu$ ($99\%$ of it in the case of Gaussian distribution).}
    \label{fig:poly_normal_form}
\end{figure}

Firstly, as it can be seen from Figures \ref{subfig:poly_big_perturbation_normal} and \ref{subfig:poly_big_perturbation_uniform}, where a large perturbation of the parameter is considered, PC visually captures the regions of non-uniqueness. 
Indeed, the polynomials solving \eqref{eq:nomal_form_system} tend to enlarge their image space when $\mu>0$, for both Gaussian and Uniform perturbations. 
This first result illustrates how each PC solution, without any ad hoc initialization of the numerical solver, inherently encapsulates information about all deterministic solution branches.

Such information, collected by PC solutions for large perturbations of $\mu$, can be further exploited to precisely identify the solution branches values for a given parameter value. 
This can be done by considering a small perturbation centered around the value of interest of the parameter. 
As it is reported in Figures \ref{subfig:poly_little_perturbation_normal} and \ref{subfig:poly_little_perturbation_uniform}, smaller perturbations evaluate the solution branches in a tighter region, where their values do not undergo great changes. 
Thus, the local extrema of the polynomials concentrate around specific local values of the solution branches determined by the perturbation support, and even if many different PC solutions are feasible, they all share the values of their local extrema. 
This way, by varying the mean of a small parameter perturbation and observing the local extrema of a generic PC solution, we obtain a computationally efficient method for branches reconstruction, as it enables their identifications for any parameter value without requiring any expensive continuation method.

\begin{figure}[t!]
    \centering

     \captionsetup[sub][figure]{skip=-30pt,slc=off,margin={0pt,0pt}}
    
    \subcaptionbox{\label{subfig:pdf_normal}}{\includegraphics[trim={0cm 0cm 0cm 0cm}, clip, width=0.47\textwidth]{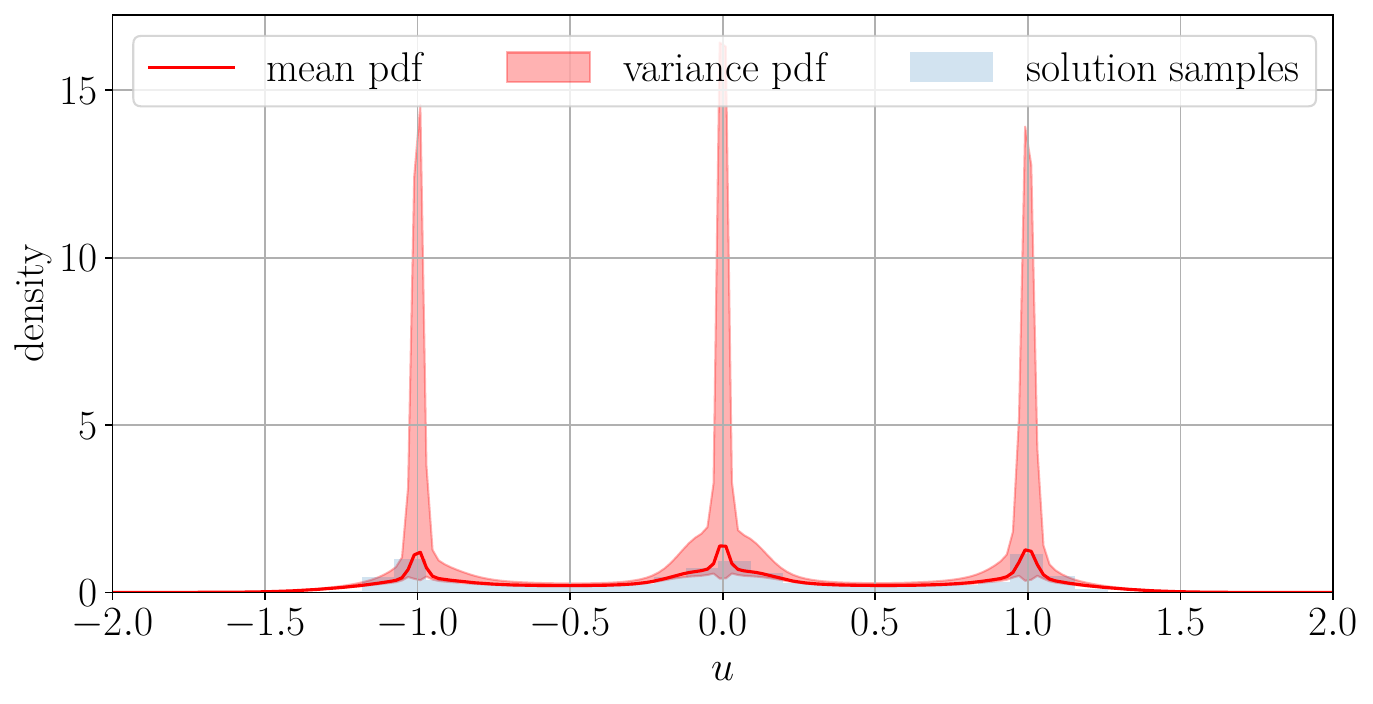}
    }
    \hfill
    \subcaptionbox{\label{subfig:pdf_uniform}}{\includegraphics[trim={0cm 0cm 0cm 0cm}, clip, width=0.47\textwidth]{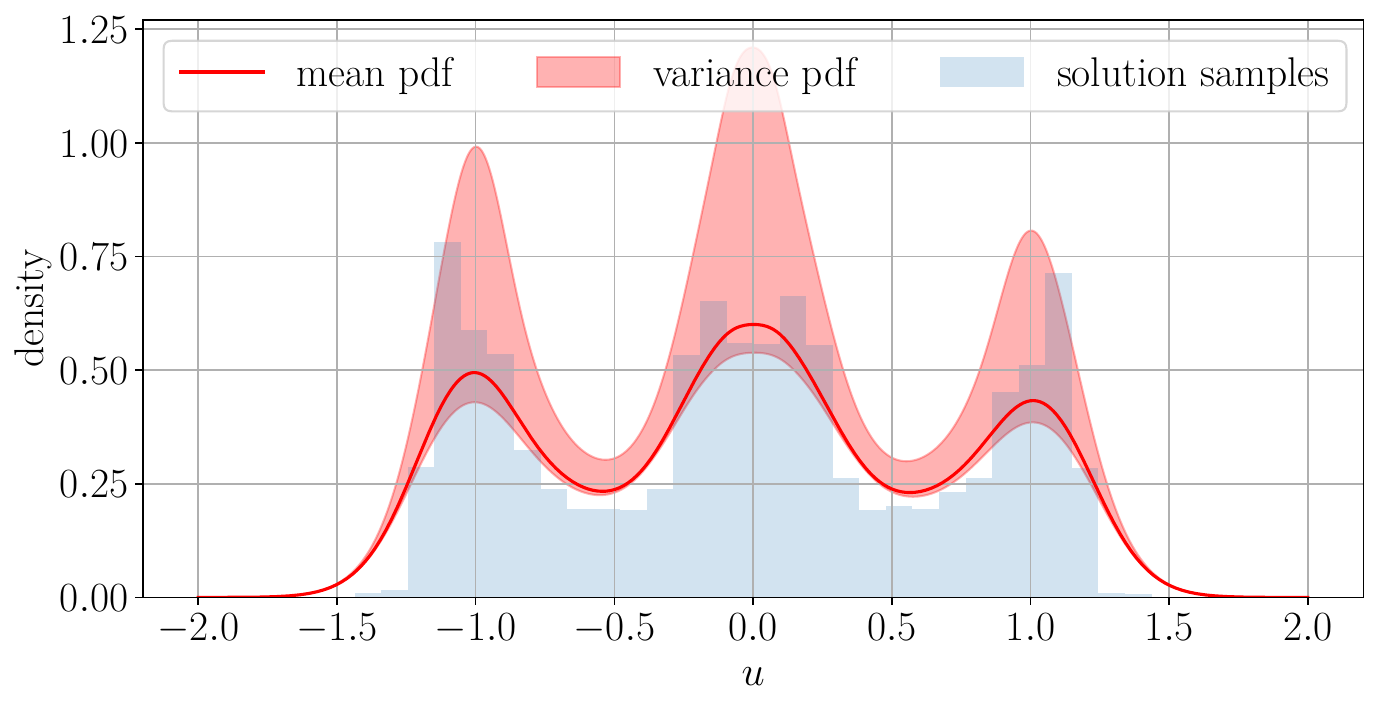}
    }
    \caption{{Mean PDF and its variance for $100$ PC representations corresponding to different random initializations of the solver with $N_{PC}=5$, for $\mu\sim\mathcal{N}(1,0.06)$ and $\mu\sim\mathcal{U}(0.8,1.2)$, left and right respectively.}}
    \label{fig:poly_var}
\end{figure}

These examples also show that when the support of the parameter distribution intersects a bifurcation region, the non-uniqueness of the deterministic solution leads to a corresponding non-uniqueness of the PC solution for the perturbed system.

Nevertheless, as shown in Figure \ref{fig:poly_var}, the probability density functions (PDFs) computed by kernel density estimation (KDE) on different PC solutions exhibit some common characteristics. 
Indeed, the mean PDF shows three peaks precisely in correspondence to the values of the three solution branches for $\mathbb{E}[\mu]=1$, which are $\{0,\pm 1\}$. 
These peaks are directly related to the locations of the local extrema of the polynomials shown in Figure~\ref{fig:poly_normal_form}, where the sampled values from the PC solutions tend to concentrate.

Regarding the variance, there is instead a difference between the two plots: for the Gaussian perturbation the three peaks of the multi-modal PDF exhibit a significantly higher variance. 
This is a consequence of the concentration of samples around the mean value for $\bar{\xi}\sim\mathcal{N}(0,1)$. 
Indeed, the PC solutions corresponding to different solver initialization show different local extrema values around $\mu=1$ (see Figure \ref{subfig:poly_little_perturbation_normal}). 
Since $\mu\sim\mathcal{N}(1,0.06)$, most samples for each PC solution cluster around the value of the correspondent local extrema closest to $\mu=1$. 
Therefore, the PC solutions exhibit different peaks in their PDFs, carrying their high variance as shown in Figure \ref{subfig:pdf_normal}. 
In contrast, the Uniform distribution enables a uniform contribution of all the local extrema of each PC solution, leading to a lower variance of the peaks (see Figure \ref{subfig:pdf_uniform}). 
Nevertheless, in both cases the variance remains low away from the peaks of the mean PDF, indicating that the local extrema of the polynomials do not cluster around any points other than those associated with the deterministic solution branches.

\begin{figure}[htb!]
    \centering
    \includegraphics[width=0.7\linewidth]{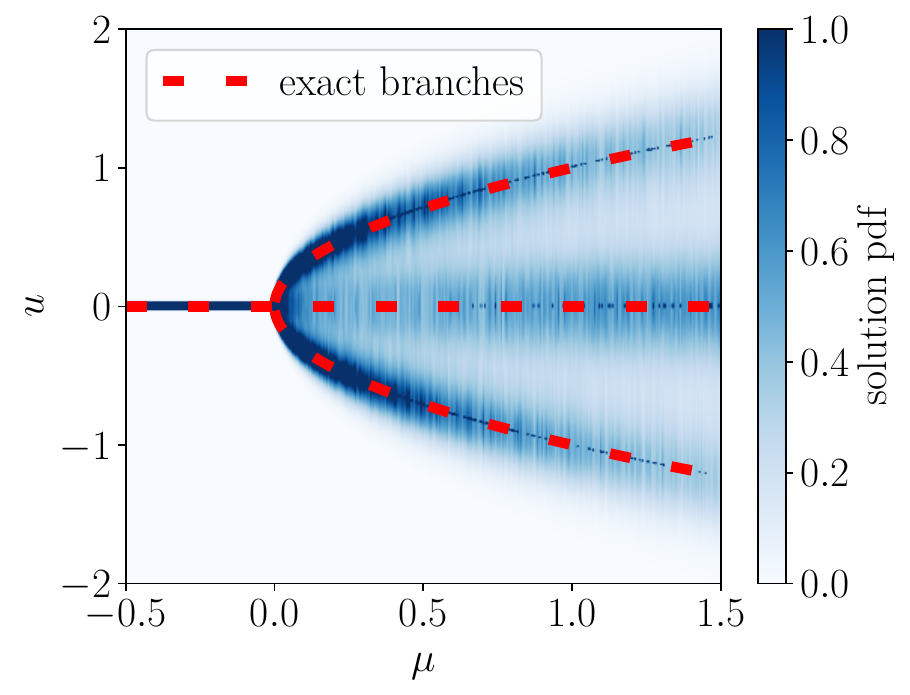}
    \caption{Comparison between deterministic bifurcation diagram and the solution PDFs computed for $500$ different uniform perturbations of $\mu$, with constant low variance and different mean $\mu\sim\mathcal{U}(\overline{\mu}-0.01,\overline{\mu}+0.01)$ with $\overline{\mu} \in [-0.5,1.5]$.}
    \label{fig:stochastic_bif_normal}
\end{figure}

{
Finally, we show the recovery of the whole bifurcation diagram through multiple runs of the perturbative approach. 
The result is shown in Figure \ref{fig:stochastic_bif_normal}, where the exact branches are precisely recovered by the solution PDF peaks, also serving as a quantitative metric of the performance of the approach in approximate these branches in a meaningful way.
This result is obtained with Uniform distributions assigned to $\mu$, solving the system multiple times varying the mean $\bar{\mu}$, with random initialization of the solver. 
Finally, we remark that, in the uniqueness regime, the solution is fully recovered by unimodal distribution of the PC solution. 
}

\section{A revisited benchmark in fluid-dynamics: the Coand\u{a} effect} \label{section:testcase}

We present here a more computationally challenging test to extend our considerations to cases in which the bifurcating behavior is not obvious.

\begin{figure}[b]
    \centering
    \includegraphics[width=0.6\textwidth]{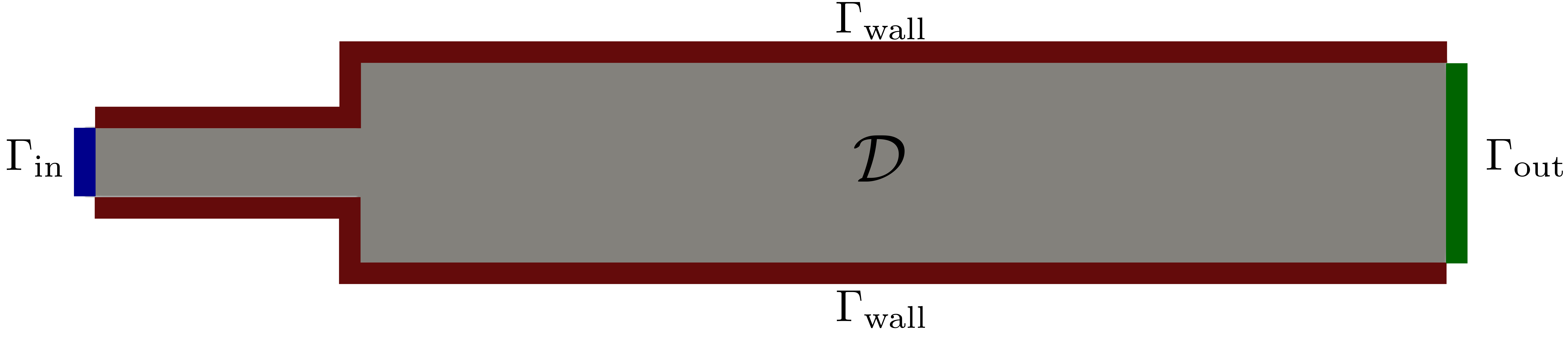}
    \caption{Computational domain representing the 2D sudden-expansion channel.}
    \label{fig:domain}
\end{figure}

Let $\vectspace{D} \in \mathbb{R}^{2}$ be the bounded domain representing the sudden-expansion channel in Figure \ref{fig:domain}.
We consider the steady incompressible Navier-Stokes equations described by:
\begin{equation}\label{eq:coanda}
    \begin{cases}
        -\mu\Delta \vect{v} +\vect{v}\cdot\nabla \vect{v} + \nabla p = 0 &\text{in}\ \ \vectspace{D}, \\
        \vect{\nabla}\cdot \vect{v} = 0 & \text{in}\ \ \vectspace{D}, \\
    \end{cases} \quad \text{with} \quad 
    \begin{cases}
        \vect{v}=\vect{v}_\text{in} & \text{on}\ \ \Gamma_{\text{in}}, \\
        \vect{v}=0 & \text{on}\ \ \Gamma_{\text{wall}}, \\
        -p\vect{n}+(\mu\nabla \vect{v})\vect{n}=0 & \text{on}\ \ \Gamma_{\text{out}}, \\
    \end{cases}\end{equation}
where the unknown state of the fluid  $(\vect{v}, p)$, encompasses its velocity $\vect{v} = (v_x, v_y)$ and pressure $p$ normalized over a constant density, where $v_x$ and $v_y$ stand respectively for the velocity's horizontal and vertical components.
The imposed boundary conditions entail a stress-free condition on the
velocity at the outlet $\Gamma_{\text{out}}$, with the outer normal denoted as
$\vect{n}=(n_{x_1}, n_{x_2})$, a no-slip homogeneous Dirichlet boundary condition
on $\Gamma_{\text{wall}}$, and a non-homogeneous Dirichlet boundary condition
$\vect{v}_\text{in}(x_1, x_2)= [20(5-x_2)(x_2-2.5), 0]$ at the inlet $\Gamma_{\text{in}}$ (see \cite{PichiStrazzullo2022} for more details regarding the benchmark).

The parameter $\mu$ represents the kinematic viscosity and plays a crucial role in the behavior of the system. Indeed, being related to the Reynolds number, it balances viscous and inertial forces, originating a symmetry-breaking phenomenon known as the Coand\u{a} effect.

Specifically, while decreasing the kinematic viscosity, the influence of inertia gradually becomes more significant, leading to a symmetry-breaking phenomenon for a critical viscosity value $\mu^{*}\approx 0.96$. 
The unique and stable symmetric branch bifurcates in two asymmetric coexisting configurations, switching its stability property, and generating the so-called \textit{pitchfork bifurcation}. 
Its study has relevance in various fields of application, ranging from shape-modeling in aerodynamics to cardiovascular blood flow in health care, recently garnering a substantial body of literature, comprising both experimental and deterministic numerical investigations \cite{Sobey,Cherdron,PichiStrazzullo2022, BravoGeometricallyParametrisedReduced2023, Khamlich2022}. 
In particular, it models the ``wall-hugging'' tendency of a viscous fluid, which is related to mitral valve regurgitation in the cardiovascular context \cite{PichiStrazzullo2022}.  

\begin{figure}[tb!]
    \includegraphics[width=\textwidth]{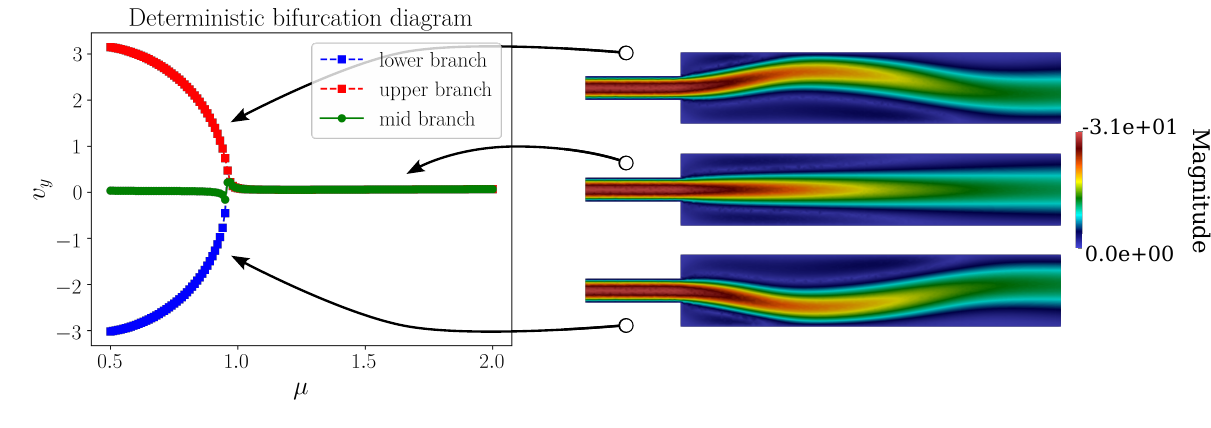}
    \caption{Deterministic FEM bifurcation diagram with the corresponding co-existing solutions belonging to the upper, middle, and lower branches for $\mu = 0.5$.}
\label{fig:bif-diag-deterministic}
\end{figure}

To build a bifurcation diagram for the Coand\u{a} problem, and thus investigate the evolution of the solution branches, one may choose to examine the behavior of the vertical velocity, the component responsible for breaking the symmetry, at a point along the channel's symmetry axis $\widehat{x}$, reported in Figure \ref{fig:bif-diag-deterministic}. 
This diagram is constructed through an iterative method, in which $\mu$ has been varied uniformly in the parametric range $\mathbb{P}=[0.5,2]$ with a step size $\Delta\mu = 0.01$.
This is complex and costly, and involves a continuation strategy to properly select the initial guess for the nonlinear solver at each new parameter value.

{For sample-based approaches to this bifurcation, we refer the reader to Appendix~\ref{sec:MC}. We now proceed with an explanation of how the SSFEM pipeline is applied in this context.}

\subsection{SSFEM approach: randomized viscosity}
We proceed to apply the SSFEM pipeline presented in Section \ref{section:SSFEM} to the Coand\u{a} test case.
We model $\mu$ as a random variable following a given distribution, and we write its K-L expansion as $\mu = \overline{\mu} + \sigma\overline{\xi}$. 
%In the case of $\mu \sim \mathcal{N}(\overline{\mu}, \sigma)$ one has $\overline{\xi_1}=\overline{\xi}\sim \vectspace{N}(0,1)$ and $\lambda_1=\sigma^2$.
Then, we expand the unknowns of the problem in \eqref{eq:coanda} through PC:
\begin{equation} \label{eq:coanda_gpc_variables}
        \vect{v}(\vect{x}, \omega) = \sum_{i=0}^{N_D^{\vect{v}}}\sum_{j=1}^{N_{PC}^{\vect{v}}}\vect{v}_{i,j}\phi_i^{\vect{v}}(\vect{x})\psi_j^{\vect{v}}(\vect{\xi}(\omega)), \quad \text{and} \quad
        p(\vect{x}, \omega) = \sum_{i=0}^{N_{D}^p}\sum_{j=1}^{N_{PC}^p}p_{i,j}\phi_i^{p}(\vect{x})\psi_j^{p}(\vect{\xi}(\omega)),
\end{equation}
where the superscripts denote quantities related to the velocity and pressure fields. 
Finally, applying Galerkin projection both in the FEM space via the test functions $\{\phi_n^{\vect{v}}(\vect{x})\}_{n=1}^{N_D^{\vect{v}}}$, $\{\phi_n^p(\vect{x})\}_{n=1}^{N_D^p}$, and in the probabilistic one with $\{\psi_m^{\vect{v}}(\xi)\}_{m=1}^{N_{PC}^{\vect{v}}}$, $\{\psi_m^p(\vect{\xi})\}_{m=1}^{N_{PC}^p}$, $\forall n\in\{1,\dots,N_D^{\vect{v}}\}$, $\forall m\in\{1,\dots,N_{PC}^{\vect{v}}\}$, we get:
\begin{equation} \label{eq:coanda_SSFEM_first_eq}
\begin{cases}
&\displaystyle\sum_{k=0}^{1}\sum_{i=0}^{N_D^{\vect{v}}}\sum_{j=1}^{N_{PC}^{\vect{v}}}\vect{v}_{i,j}\int_{\set{R}}\sqrt{\lambda_{k}}\xi_{k}\psi_j^{\vect{v}}(\vect{\xi})\psi_m^{(\vect{v})}(\vect{\xi})d\nu_{\vect{\xi}}\int_{\vect{D}}\nabla\phi_i^{\vect{v}}(\vect{x})\cdot\nabla\phi_n^{\vect{v}}(\vect{x})d\vect{x}  \\%[.6em] 
&\displaystyle\quad+ \sum_{i=0}^{N_D^{\vect{v}}}\sum_{j=1}^{N_{PC}^{\vect{v}}}\sum_{l=0}^{N_D^{\vect{v}}}\sum_{h=1}^{N_{PC}^{\vect{v}}} \vect{v}_{i,j}\vect{v}_{l,h}\int_{\set{R}}\psi_j^{\vect{v}}(\vect{\xi})\psi_h^{\vect{v}}(\vect{\xi})\psi_m^{\vect{v}}(\vect{\xi})d\nu_{\vect{\xi}} \int_{\vect{D}}(\phi_i^{\vect{v}}(\vect{x})\cdot \nabla\phi_l^{\vect{v}}(\vect{x}))\phi_n^{\vect{v}}(\vect{x})d\vect{x}  \\%[.6em] 
&\displaystyle\quad-\sum_{i=0}^{N_D^p}\sum_{j=1}^{N_{PC}^p}p_{i,j}\int_{\set{R}}\psi_j^{p}(\vect{\xi})\psi_m^{\vect{v}}(\vect{\xi})d\nu_{\vect{\xi}}\int_{\vect{D}}\phi_i^{p}(\vect{x})\nabla\cdot\phi_n^{\vect{v}}(\vect{x})d\vect{x} = 0, \\%[.6em] 
&\displaystyle\sum_{i=0}^{N_D^{\vect{v}}}\sum_{j=1}^{N_{PC}^{\vect{v}}}\vect{v}_{i,j}\int_{\vect{D}}\nabla\cdot\phi_i^{\vect{v}}(\vect{x})\phi_n^{p}(\vect{x})d\vect{x}\int_{\set{R}}\psi^{\vect{v}}_j(\vect{\xi})\psi^{p}_m(\vect{\xi})d\nu_{\vect{\xi}} = 0.
\end{cases}
\end{equation}
It can be rewritten in the compact form
\begin{equation} \label{eq:coanda_SSFEM_tf}
\begin{cases}\vect{A}^T\cdot\vect{U}\cdot(\vect{E}^{(0)}+\vect{E}^{(1)})+ \vect{I} - \vect{C}^T\cdot\vect{Q}\cdot\vect{G}= \vect{0}, \\
        \vect{D}^T\cdot\vect{U}\cdot\vect{H} = \vect{0}
    \end{cases}
\end{equation}
where we denoted with $\vect{U}\in \set{R}^{N_D^{\vect{v}}\times N_{PC}^{\vect{v}}}$ and $\vect{Q}\in \set{R}^{N_D^p\times N_{PC}^p}$, respectively, the unknown coefficients for the velocity and pressure fields. 
Denoting with $\odot$ the Hadamard product, the FEM and stochastic matrices can be defined in Equations \eqref{eq:coanda_SSFEM_det_mat} and \eqref{eq:coanda_SSFEM_stochastic_mat}, respectively as
\begin{equation}
\begin{gathered} \label{eq:coanda_SSFEM_det_mat}
            A_{j,n} = \int_{\vect{D}}\nabla\phi_j^{\vect{v}}\cdot\nabla\phi_n^{\vect{v}}d\vect{x}, \quad
            B_{n,j,h} = \int_{\vect{D}}(\phi_j^{\vect{v}}\cdot \nabla\phi_h^{\vect{v}})\phi_n^{\vect{v}}d\vect{x}, \quad
            C_{j,n} = \int_{\vect{D}}\phi_j^{p}\nabla\cdot\phi_n^{\vect{v}}d\vect{x}, \\
            D_{j,n} = \int_{\vect{D}}\nabla\cdot\phi_j^{\vect{v}}\phi_n^{p}d\vect{x}, \quad
            I_{n,m} =(\vect{U}^T\cdot\vect{B}_n\cdot\vect{U})\odot\vect{F}_m,
\end{gathered}
\end{equation}

\begin{equation}
\begin{gathered} \label{eq:coanda_SSFEM_stochastic_mat}
            E_{i,m}^{(0)} = \int_{\set{R}}\sqrt{\lambda_0}\xi_0\psi_i^{\vect{v}}\psi_m^{\vect{v}}d\nu_{\vect{\xi}}, \quad
            E_{i,m}^{(1)} = \int_{\set{R}}\sqrt{\lambda_1}\xi_1\psi_i^{\vect{v}}\psi_m^{\vect{v}}d\nu_{\vect{\xi}}, \quad
            F_{m,l,i} = \int_{\set{R}}\psi_i^{\vect{v}}\psi_l^{\vect{v}}\psi_m^{\vect{v}}d\nu_{\vect{\xi}}, \\
            G_{i,m} = \int_{\set{R}}\psi_i^{p}\psi_m^{\vect{v}}d\nu_{\vect{\xi}}, \quad
            H_{i,m} = \int_{\set{R}}\psi^{\vect{v}}_i\psi^{p}_md\nu_{\vect{\xi}}.
\end{gathered}
\end{equation}

{
System \eqref{eq:coanda_SSFEM_tf} is solved using a Newton–Krylov method enhanced with a local line-search strategy.   
Regarding the imposition of the boundary conditions defined in the general formulation \eqref{eq:coanda}, for the inlet and wall nodes $\Gamma_{\text{in}}$ and $\Gamma_{\text{wall}}$ respectively, they are enforced in a mean sense. In particular, the first coefficient of the PC expansion, responsible for the mean value of the solution, is prescribed according to the deterministic boundary conditions. The variance, corresponding to the weighted contribution of all higher-order coefficients, is instead set to zero by imposing all remaining PC coefficients to vanish at those nodes.
}
\section{Results} \label{sec:results}
In this section, we discuss the numerical results of the SSFEM pipeline when applied to the Navier-Stokes equations modeling the Coand\u{a} bifurcation problem.
Consistently with other studies in the deterministic setting \cite{PichiStrazzullo2022,Khamlich2022}, we exploit the vertical component of the velocity field to discuss the non-uniqueness feature of the model. 

To perform an explorative analysis of the problem we start by treating the viscosity parameter as a Gaussian random variable centered near the known bifurcation point with a small variance. This specific choice has been made to study the system's response to little perturbations, in the case in which the multiple coexisting branches (see Figure \ref{fig:bif-diag-deterministic}) are still reasonably close to each other. Note that, as underlined in Section \ref{sec:pitchfork}, in the case of no prior knowledge of the bifurcation point one could largely perturb the parameter and observe the resulting PC solutions to have hints on its position in the parameter space.

Successively,  we show the effect of the parameter distribution on the solution, presenting also results for uniform perturbations of the viscosity. 
Moreover, we discuss the effects of coarseness and symmetry of the grid, recovering some known results in the literature about the computational grid's requirements for bifurcation detection.

Finally, we extend the analysis to multiple viscosity mean values far from the bifurcation point, and discuss the discovered strict and surprising relationship between the PC representation and the deterministic bifurcation diagram.

The numerical simulations are conducted thanks to SSFEM implementation based on the deal.II library. The code\footnote{https://github.com/ICGonnella/SSFEM-Coanda-Effect.git} has a Python interface and it is parallelized to enable faster computations on larger meshes.
The space discretization has been carried out using the Taylor-Hood ($\mathbb{P}_2$-$\mathbb{P}_1$) elements, and the nonlinear solver is the Newton-Raphson one implemented in deal.II as the default option of the \texttt{NonlinearSolverSelector} class, enriched with periodic linear search steps.

\subsection{Statistical moments around the bifurcation point} \label{sec:moments_around_bifurcation}

As a first step, we are interested in understanding whether the statistical moments inferred from the PC representation of the velocity solution are informative on its bifurcating behavior. 

To this aim, we perform SSFEM computation assigning a Gaussian distribution to the viscosity $\mu\sim\vectspace{N}(\overline{\mu},\sigma)$, where the mean is located in the non-uniqueness regime of the parametric space. Therefore, following the intuitive interpretation of a bifurcating behavior where drastic changes follow a small variation of the parameter, we apply a small stochastic perturbation to the kinematic viscosity $\mu$. In particular, we expect that the stochastic solution exhibits greater variance at the spatial coordinates for which the discrepancy between the coexisting deterministic solutions is maximized.

\begin{figure}[b]
\centering
        \captionsetup[sub][figure]{oneside,margin={-7cm,0cm},skip=-1.6cm}
        \subcaptionbox{\label{subfig:vel_mean_0.9_1_Normal}}{        \includegraphics[trim={12cm 10.3cm 6.0cm 11.cm}, clip, width=.65\textwidth]{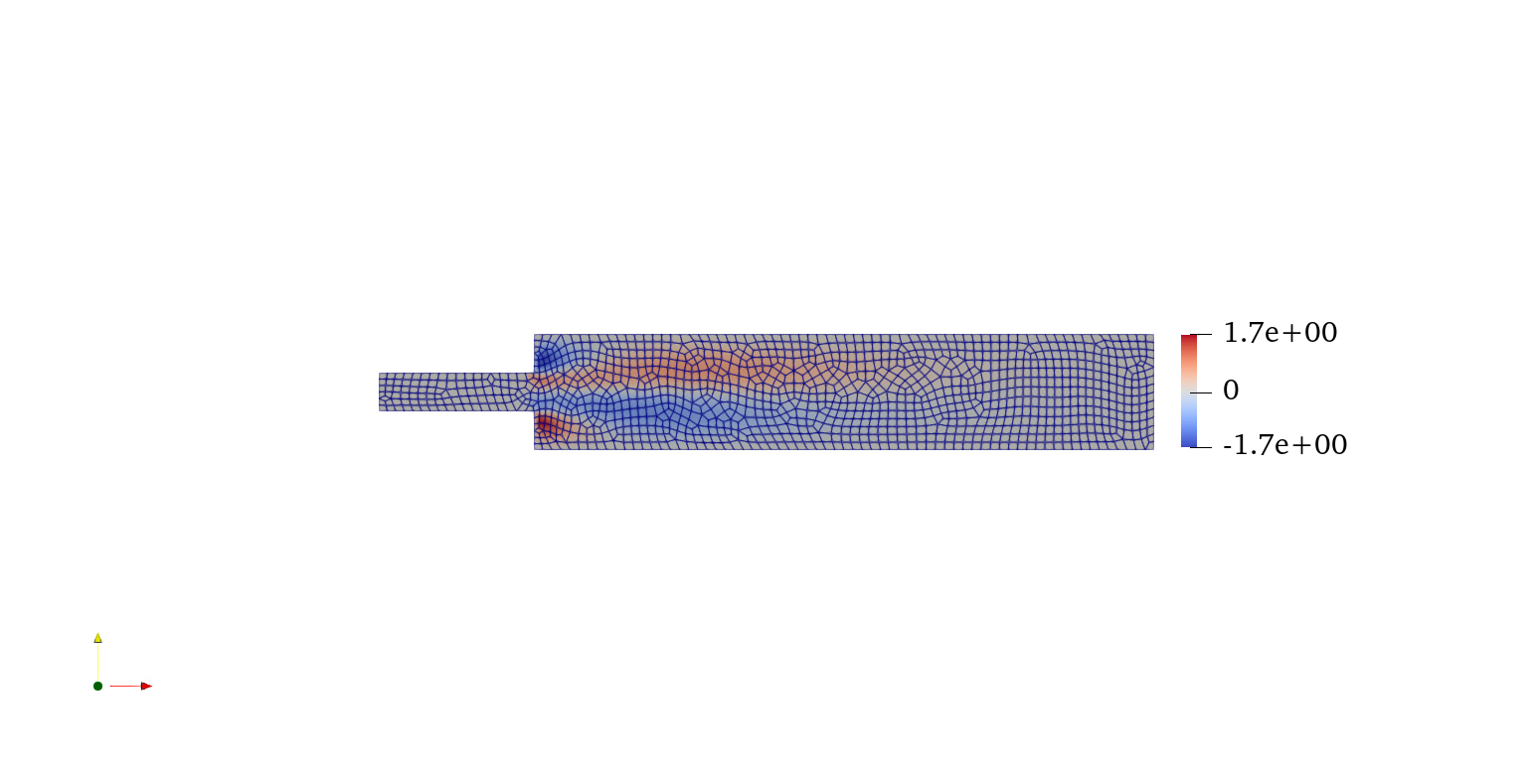}
        }
    \hfill
        \subcaptionbox{\label{subfig:vel_variance_0.9_1_Normal}}{\includegraphics[trim={12cm 10.3cm 6.0cm 11.cm}, clip, width=.65\textwidth]{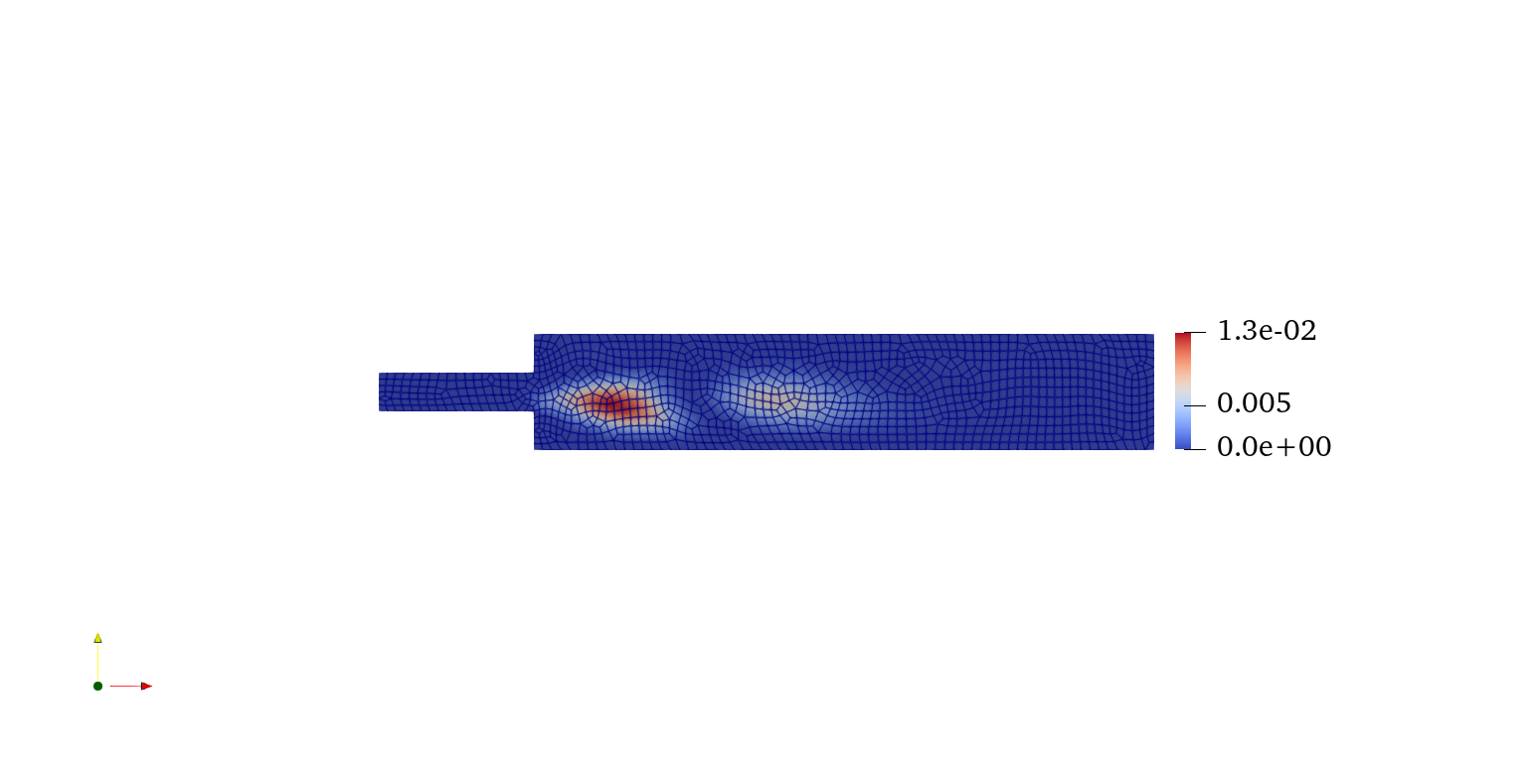}
        \captionsetup[sub][figure]{oneside,margin={-7cm,0cm},skip=-1.6cm}}
        \captionsetup{subrefformat=parens}
    \caption{SSFEM statistics \cref{subfig:vel_mean_0.9_1_Normal} mean and  \cref{subfig:vel_variance_0.9_1_Normal} variance on a quadrilateral grid of $1275$ nodes, obtained using a linear PC expansion ($N_{PC}=1$) and $\mu\sim\vectspace{N}(0.9,0.001)$.}
    \label{fig:vel_stats_0.9_Normal}
\end{figure}

In Figure \ref{fig:vel_stats_0.9_Normal} we show the mean and the variance of $v_y$ computed relying only on a linear polynomial expansion of the velocity. Interestingly, it can be noticed that the grid point corresponding to the highest variance value is located in the post-inlet region, where the vertical components of the admissible velocity solutions break the symmetry. Therefore, having this statistical information, one can easily locate the point of maximum discrepancy of the coexisting solutions directly in the computational domain. 

Nevertheless, from the above variance plots, it is not possible to distinguish between a bifurcating behavior of the solution and a simple high variability induced by other factors. 
To investigate this difference, possibly encoded in the PC solution, we search for some insights on its polynomials. 
More precisely, we analyze the shape of the PC solution evaluated in the regions of the domain showing higher variance.

We define the \textit{sampling zone} as the subset of supp$(\vect{\xi})$ at which the $99\%$ of the measure $\nu_{\vect{\xi}}$ is concentrated (e.g.\ it corresponds to the interval $(-3,3)$ for standard normal basic random variables $\vect{\xi}$).
In particular, we aim to recover the correlation between the minima and maxima of the polynomials found in the sampling zone, and the presence of a bifurcating behavior, as it happened with the pitchfork normal form in section \ref{sec:pitchfork}.

\begin{figure}[t]
    \centering
    \captionsetup[sub][figure]{skip=-68pt,slc=off,margin={0pt,0pt}}
    \subcaptionbox{\label{subfig:vel_variance_0.9_2_Normal}}{\includegraphics[trim={12cm 9.3cm 9.6cm 11.5cm}, clip, width=0.65\textwidth]{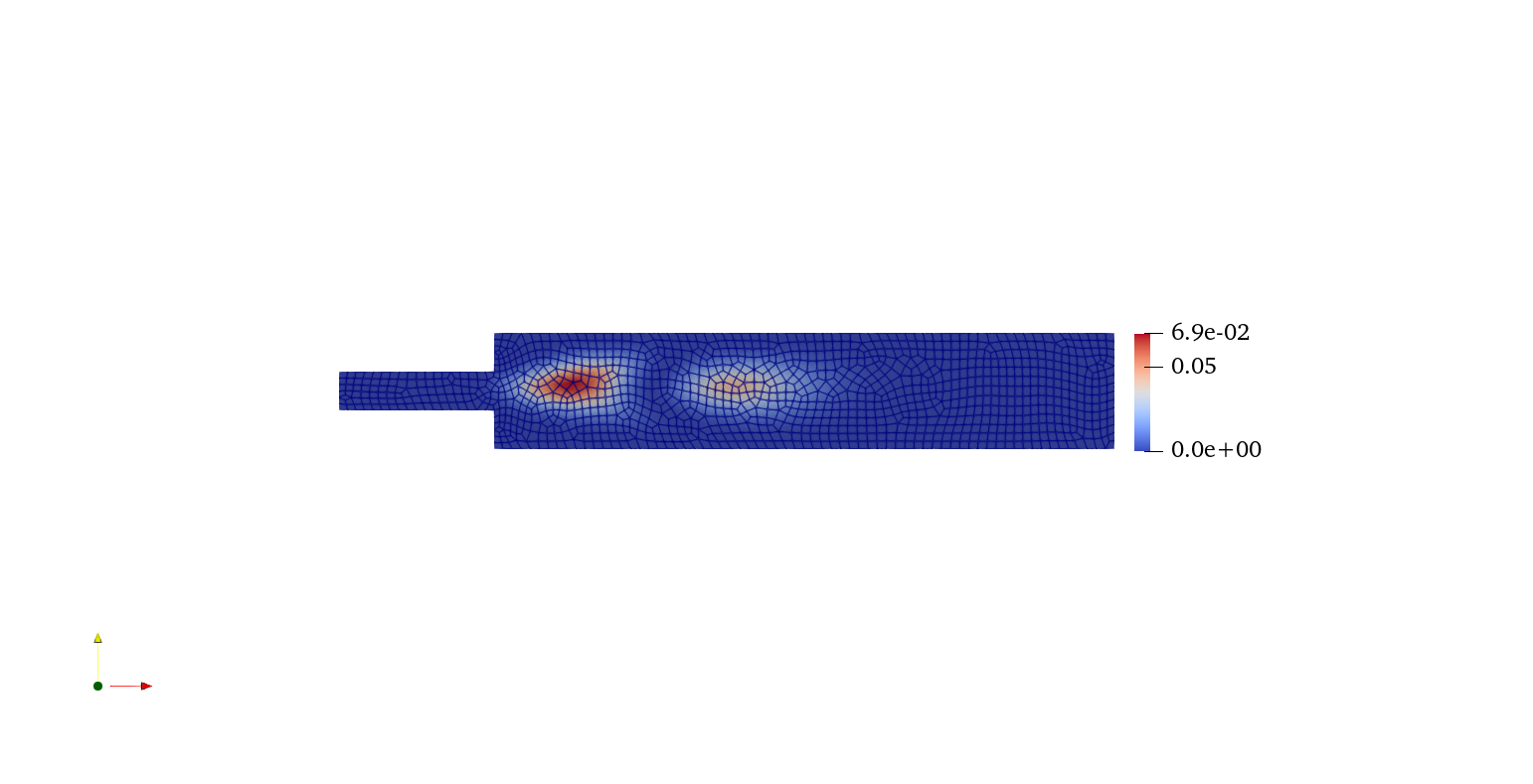}
    %\captionsetup[sub][figure]{oneside,margin={-2cm,-2cm},skip=-2cm}
    }
    \hfill
    \subcaptionbox{\label{subfig:poly_vel_variance_0.9_2_Normal}}{\includegraphics[trim={0cm 0cm 0cm 0cm}, clip, width=0.3\textwidth]{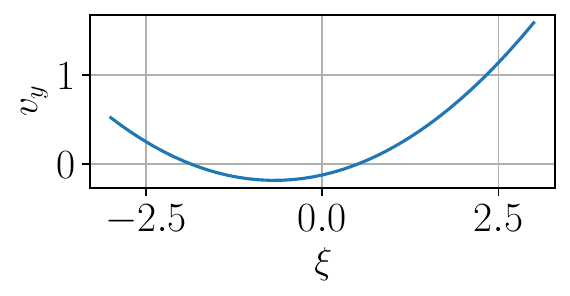}
    %\captionsetup[sub][figure]{oneside,margin={-6cm,0cm},skip=-2cm}
    }
    
    \subcaptionbox{\label{subfig:vel_variance_0.9_3_Normal}}{\includegraphics[trim={12cm 9.3cm 9.6cm 11.5cm}, clip, width=0.65\textwidth]{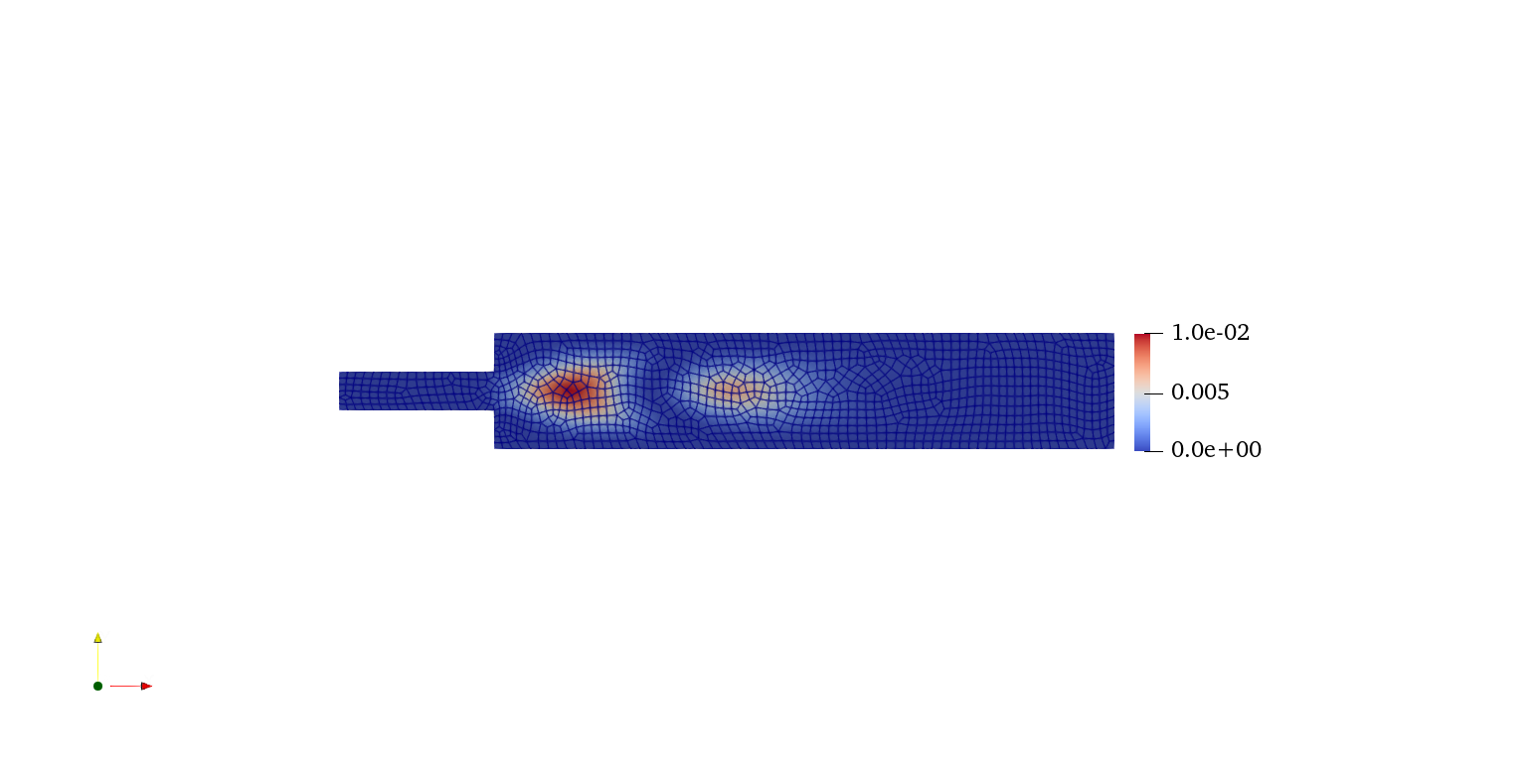}
    %\captionsetup[sub][figure]{oneside,margin={-6cm,0cm},skip=-2cm}
    }
    \hfill
    \subcaptionbox{\label{subfig:poly_vel_variance_0.9_3_Normal}}{\includegraphics[trim={0cm 0cm 0cm 0cm}, clip, width=0.3\textwidth]{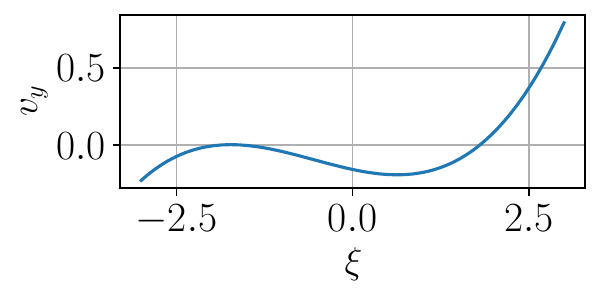}
    %\captionsetup[sub][figure]{oneside,margin={-6cm,0cm},skip=-2cm}
    }
    
    \subcaptionbox{\label{subfig:vel_variance_0.9_4_Normal}}{\includegraphics[trim={12cm 9.3cm 9.6cm 11.5cm}, clip, width=0.65\textwidth]{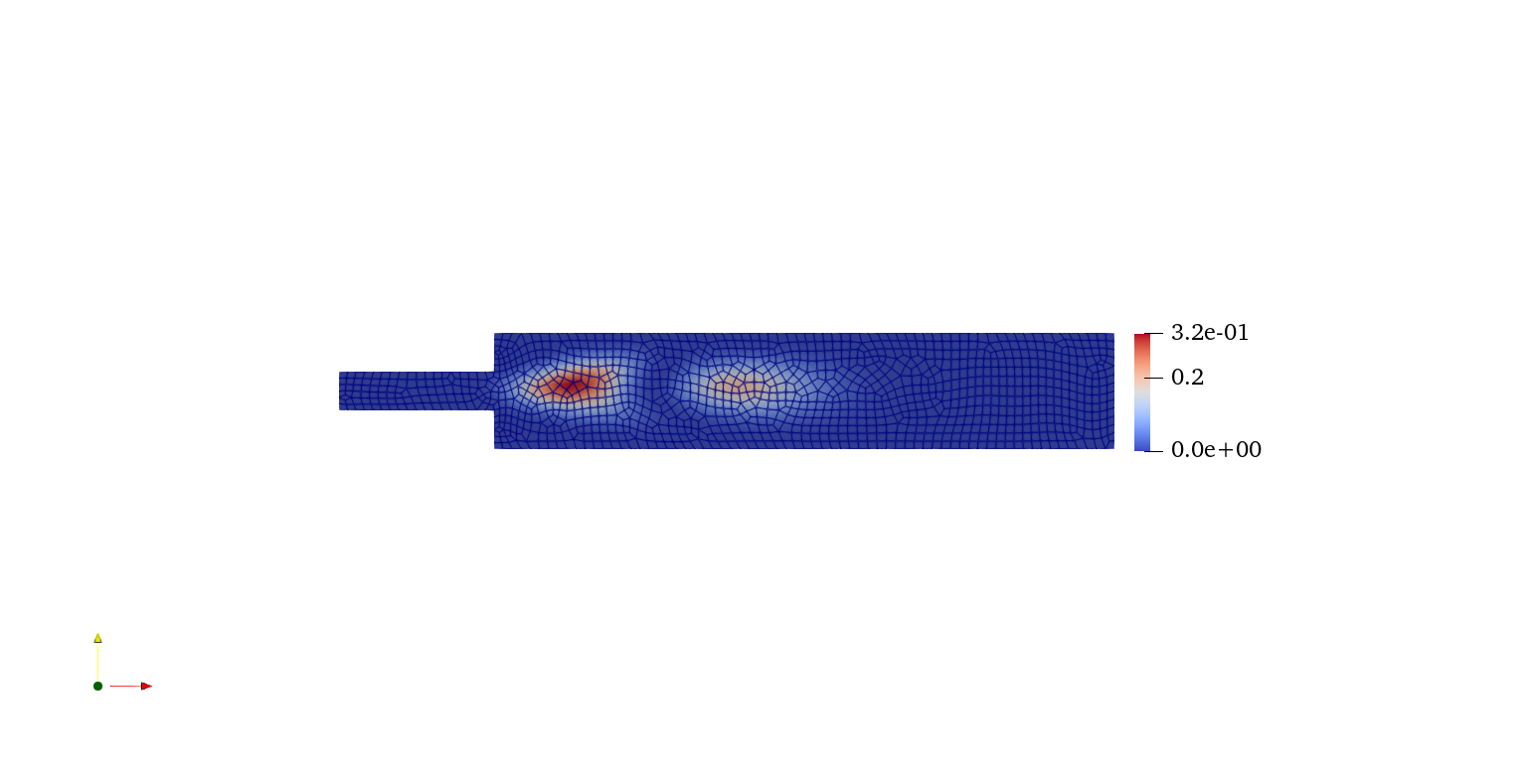}
    %\captionsetup[sub][figure]{oneside,margin={-6cm,0cm},skip=-2cm}
    }
    \hfill
    \subcaptionbox{\label{subfig:poly_vel_variance_0.9_4_Normal}}{\includegraphics[trim={0cm 0cm 0cm 0cm}, clip, width=0.3\textwidth]{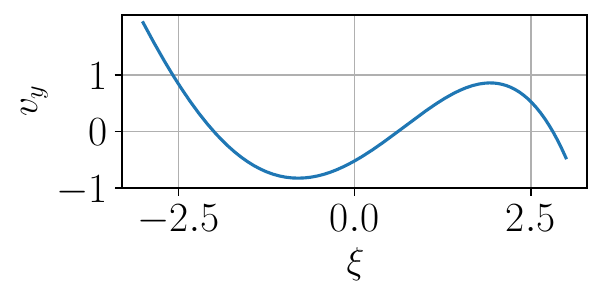}
    %\captionsetup[sub][figure]{oneside,margin={-6cm,-2cm},skip=0cm}
    }
    \captionsetup{subrefformat=parens}

    \caption{SSFEM results on a quadrilateral grid of $1275$ nodes with $\mu\sim\vectspace{N}(0.9,0.001)$. Variance magnitude and solution at $\bar{\boldsymbol{x}} = (15, 3.75)$ (point of maximum variance), for $N_{PC}$=2 \cref{subfig:vel_variance_0.9_2_Normal}\cref{subfig:poly_vel_variance_0.9_2_Normal}, $N_{PC}$=3 \cref{subfig:vel_variance_0.9_3_Normal}\cref{subfig:poly_vel_variance_0.9_3_Normal}, and $N_{PC}$=4 \cref{subfig:vel_variance_0.9_4_Normal}\cref{subfig:poly_vel_variance_0.9_4_Normal}.} 
    % \cref{subfig:vel_variance_0.9_2_Normal}\cref{subfig:vel_variance_0.9_3_Normal}\cref{subfig:vel_variance_0.9_4_Normal}, and solution at $\bar{\boldsymbol{x}}$ on the right in Figures \cref{subfig:poly_vel_variance_0.9_2_Normal}\cref{subfig:poly_vel_variance_0.9_3_Normal}\cref{subfig:poly_vel_variance_0.9_4_Normal} respectively for $N_{PC}$ taking the values $3,4,5$.}
    \label{fig:vel_variance_0.9_Normal}
\end{figure}

As it can be seen in Figures \ref{fig:vel_stats_0.9_Normal}, \ref{fig:vel_variance_0.9_Normal}, the variance magnitude tends to increase in absolute value as the polynomial degree increases but remains localized at the same grid points. 
Moreover, along with the variance plots, we show the shapes of the polynomials for the high variance domain point $\bar{\boldsymbol{x}} = (15, 3.75)$, belonging to the horizontal symmetry axis. 
It can be noticed that up to degree two there is no evidence of bifurcating behavior, as only one local extrema appears. 
Instead, when reaching degree three, multiple local extrema are visible in the sampling area, whose values for degree four (see Figure \ref{subfig:poly_vel_variance_0.9_4_Normal}) seem to tend to the values of the two stable solution branches in that specific domain point (see Figure \ref{fig:bif-diag-deterministic}). 

Therefore, a correlation between the value of the local extrema in the sampling region and the branching solutions is confirmed, suggesting that a bifurcating behavior is characterized by the joint occurrence of both a high variance and multiple local extrema for the solution's PC representation. Moreover, one can note that the high-variance region keeps improving its symmetry w.r.t.\ the horizontal axis starting from degree two, while for degree one (see Figure \ref{subfig:vel_variance_0.9_1_Normal}) it is not perfectly centered yet. This, together with the progressive adjustment of the peaks values versus the stable branches values (from Figure \ref{subfig:poly_vel_variance_0.9_3_Normal} to Figure \ref{subfig:poly_vel_variance_0.9_4_Normal}), remarks even more the dependence of the PC representation's accuracy of the bifurcating phenomenon to the polynomial degree considered.

\begin{figure}[hbt]
    \centering
    \includegraphics[width=0.44\textwidth]{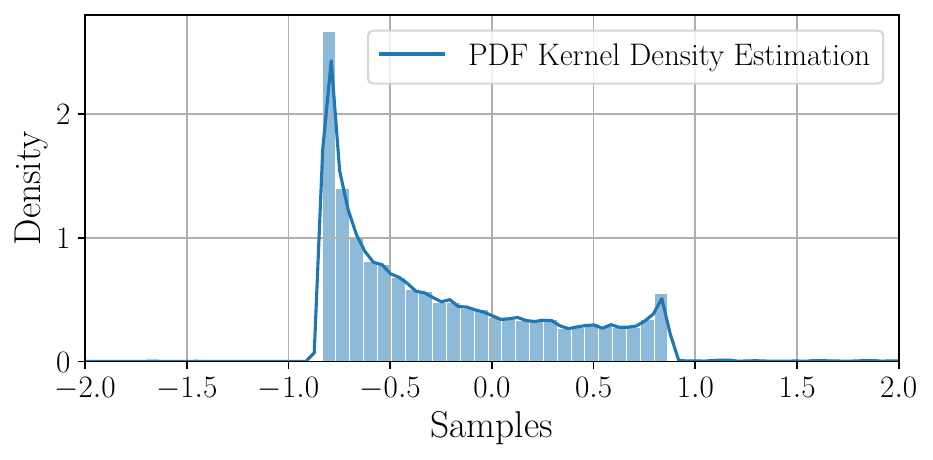}
    \caption{KDE of the velocity PDF for $\mu\sim\vectspace{N}(0.9,0.001)$ at point $\bar{\boldsymbol{x}}$ of maximum variance.}
    \label{fig:kde_Normal}
\end{figure}

The PDF of the vertical component of the velocity solution is reconstructed in Figure \ref{fig:kde_Normal} thanks to its PC representation, as a large number of samples can be collected and used for the KDE. 
Two peaks can be observed at the two local extrema of the correspondent polynomial of Figure \ref{subfig:poly_vel_variance_0.9_4_Normal}, around the values $\{-1,1\}$. 
Indeed, as anticipated before, samples of the solution are obtained by repeatedly applying the polynomial function to samples collected from a normal random variable $\xi\sim\mathcal{N}(0,1)$. 
Thus, it is reasonable that we observe a higher peak around $-1$ in the density plot, as the abscissa of the corresponding local minimum is close to zero, and thus significantly more frequent in the sampling process of $\xi$. 
This has already been noted in Section \ref{sec:pitchfork}, and justifies the higher variance of PC solutions PDFs in case of Gaussian perturbation of $\mu$.

% \fp{nella caption il punto di massima varianza non per $\mu$ ma per la soluzione}
%Given the analysis of the influence of a small Gaussian noise applied to the viscosity parameter also in a computationally complex case, we aim to investigate the influence of assuming a different a-priori distribution: the Uniform one. Indeed, a big advantage of SSFEM representation is that it is not restrictive to Gaussian variables, but can be generalized using gPC expansion (see Section \ref{section:gPC_expansion}). We conduct such analysis in the following section, where we discuss uniform distributions for the viscosity and the impact of different computational meshes. 
% In this way we propose to deepen the roles of PC expansion and mesh parameters in the probabilistic solution accuracy.

\subsection{Influence of parameter distributions and mesh discretizations}
Motivated by the asymmetry of the peaks in Figure \ref{fig:kde_Normal} and by the analysis conducted on the pitchfork normal form in Section \ref{sec:pitchfork}, we report the solution variance and the correspondent polynomial analysis in case of $\mu\sim\mathcal{U}(a,b)$, where $a, b$ are chosen to have the distribution centered at value $0.9$ with same variance as in the previous section. 

\begin{figure}[ht]
    \centering
    \captionsetup[sub][figure]{skip=-68pt,slc=off,margin={0pt,0pt}}
    \subcaptionbox{\label{subfig:vel_variance_0.9_2_uniform}}{\includegraphics[trim={12cm 9.3cm 9.6cm 11.5cm}, clip, width=0.65\textwidth]{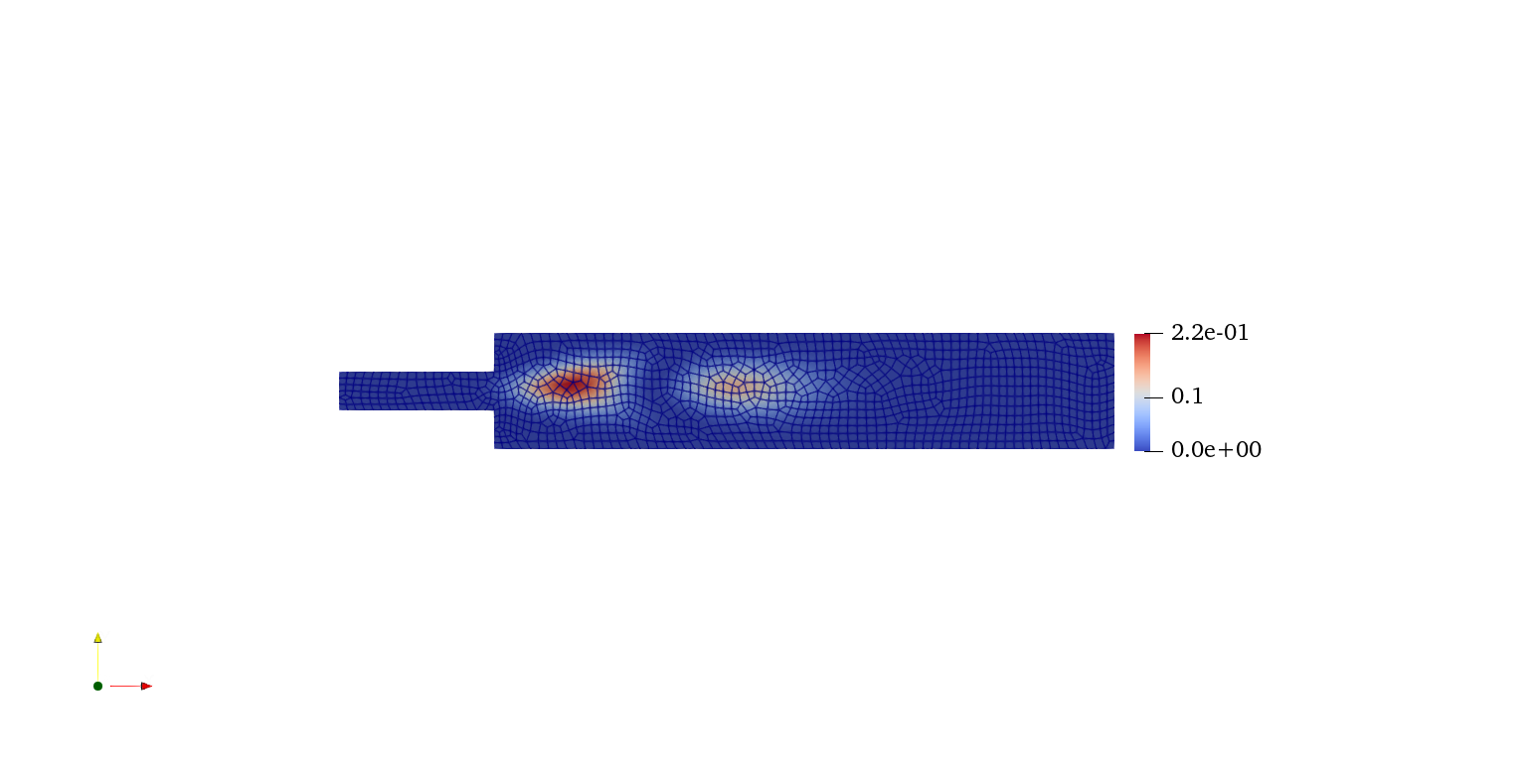}
    }
    \hfill
    \subcaptionbox{\label{subfig:poly_vel_variance_0.9_2_uniform}}{\includegraphics[trim={0cm 0cm 0cm 0cm}, clip, width=0.3\textwidth]{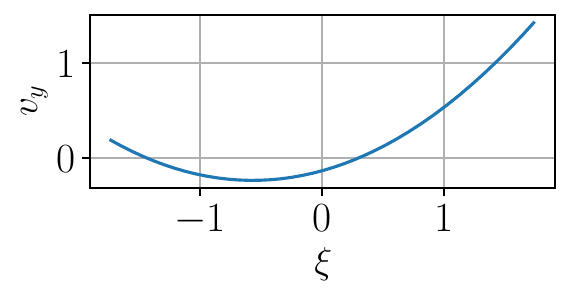}
    }

    \subcaptionbox{\label{subfig:vel_variance_0.9_3_uniform}}{\includegraphics[trim={12cm 9.3cm 9.6cm 11.5cm}, clip, width=0.65\textwidth]{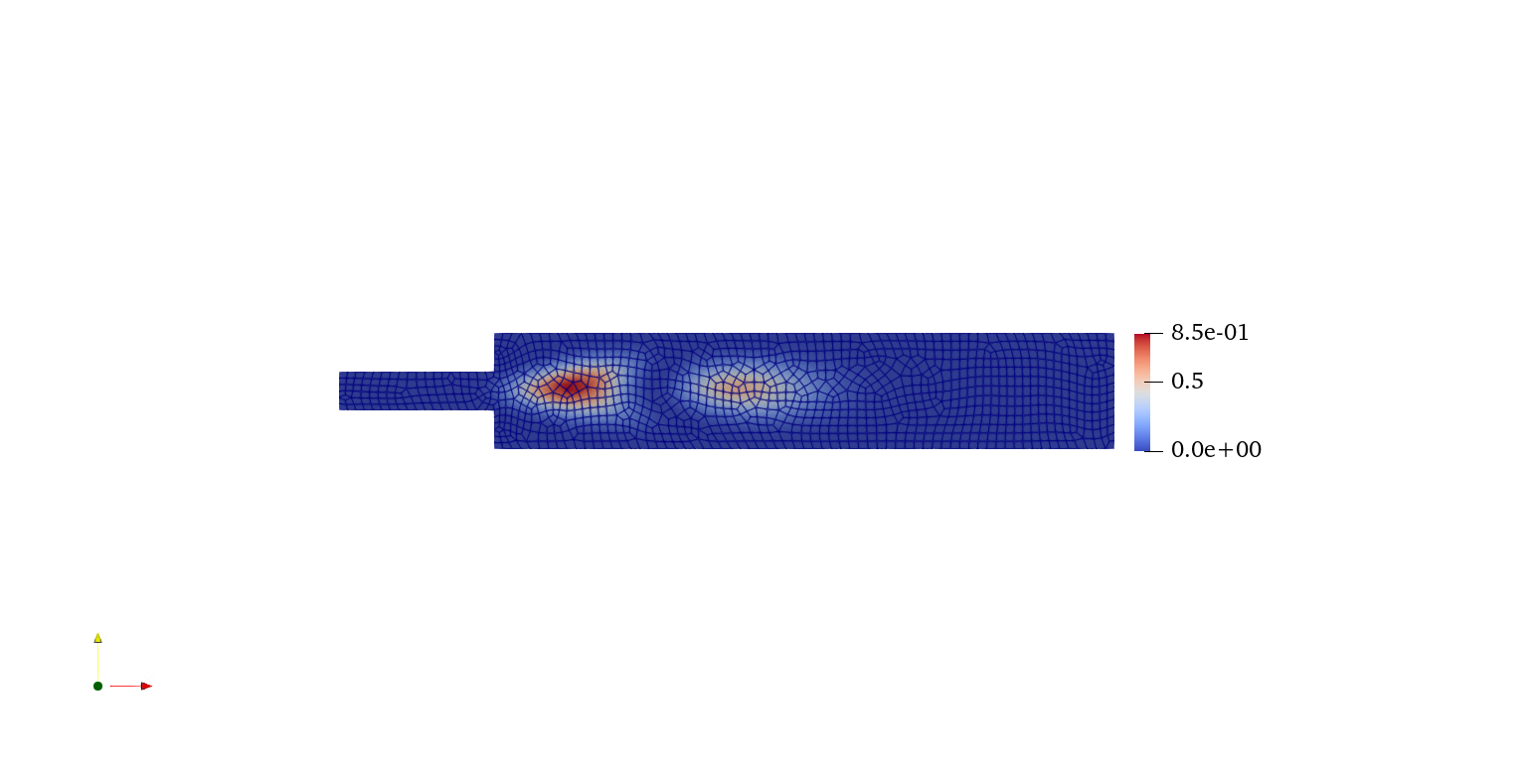}
    }
    \hfill
    \subcaptionbox{\label{subfig:poly_vel_variance_0.9_3_uniform}}{\includegraphics[trim={0cm 0cm 0cm 0cm}, clip, width=0.3\textwidth]{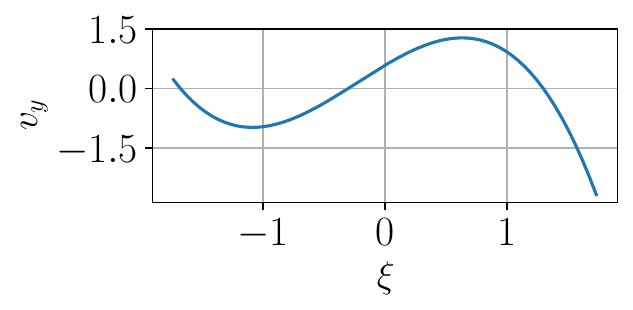}
    }
    
    \subcaptionbox{\label{subfig:vel_variance_0.9_4_uniform}}{\includegraphics[trim={12cm 9.3cm 9.6cm 11.5cm}, clip, width=0.65\textwidth]{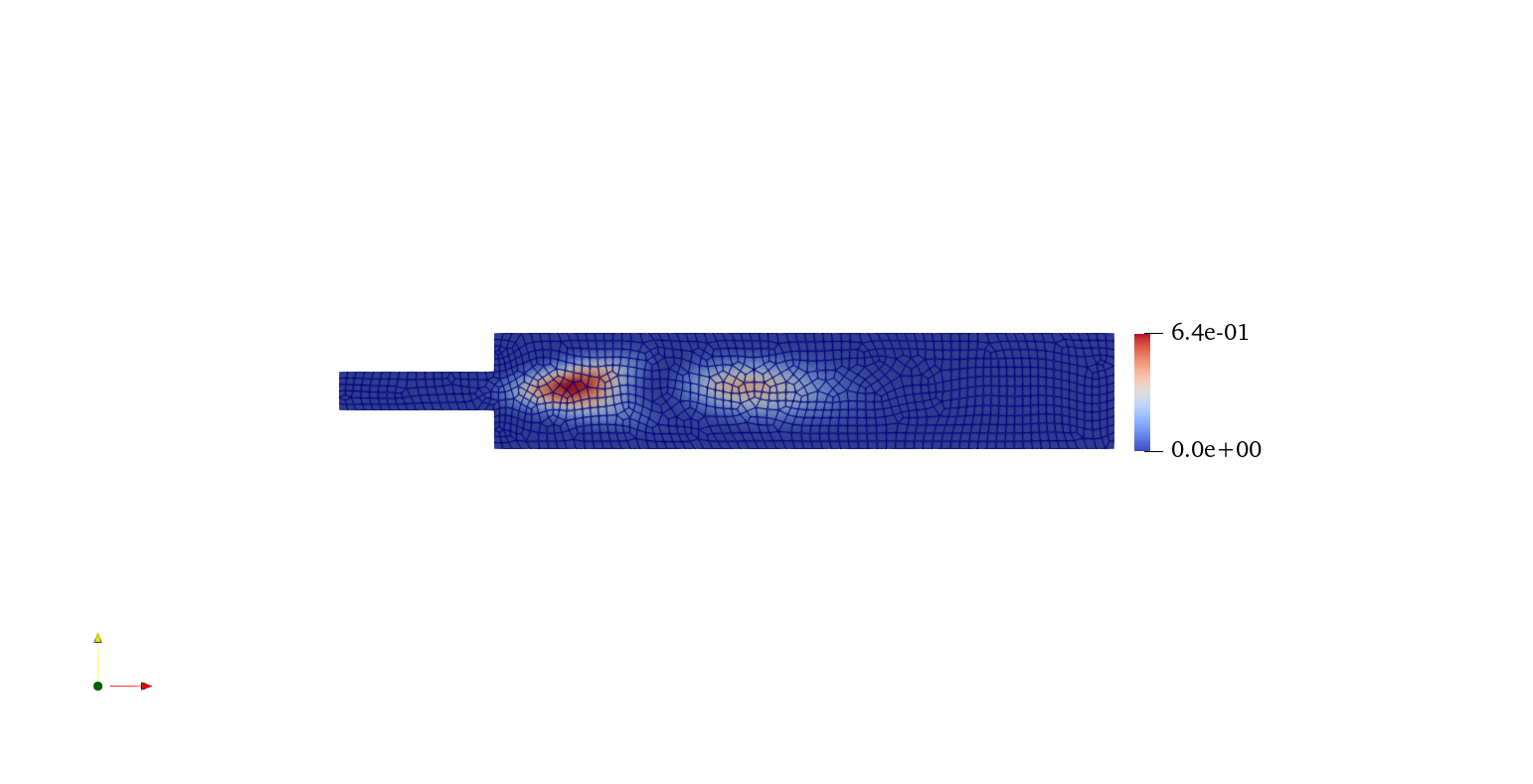}
    }
    \hfill
    \subcaptionbox{\label{subfig:poly_vel_variance_0.9_4_uniform}}{\includegraphics[trim={0cm 0cm 0cm 0cm}, clip, width=0.3\textwidth]{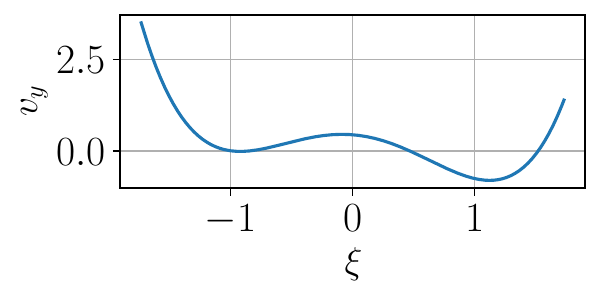}
    }

    \subcaptionbox{\label{subfig:vel_variance_0.9_5_uniform}}{\includegraphics[trim={12cm 9.3cm 9.6cm 11.5cm}, clip, width=0.65\textwidth]{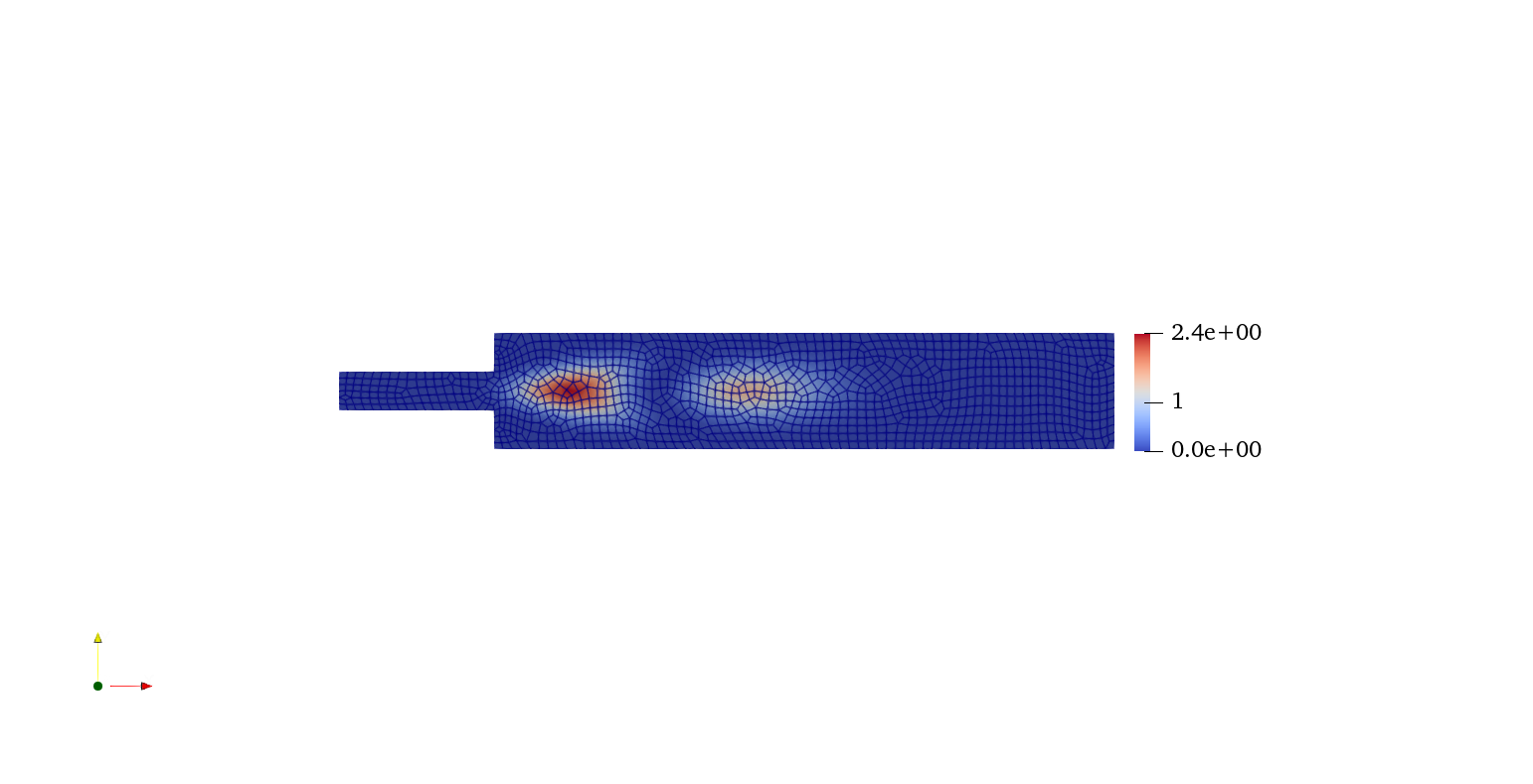}
    }
    \hfill
    \subcaptionbox{\label{subfig:poly_vel_variance_0.9_5_uniform}}{\includegraphics[trim={0cm 0cm 0cm 0cm}, clip, width=0.3\textwidth]{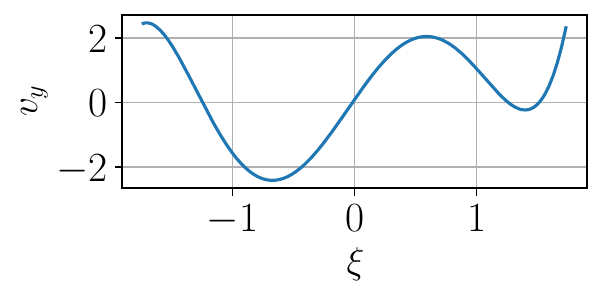}
    }
    \captionsetup{subrefformat=parens}

    \caption{SSFEM results with $\mu\sim\vectspace{U}(0.845, 0.955)$. Variance magnitude
    and solution at $\bar{\boldsymbol{x}}$, for $N_{PC}=2$ \cref{subfig:vel_variance_0.9_2_uniform}\cref{subfig:poly_vel_variance_0.9_2_uniform}, $N_{PC}=3$ \cref{subfig:vel_variance_0.9_3_uniform}\cref{subfig:poly_vel_variance_0.9_3_uniform}, $N_{PC}=4$ \cref{subfig:vel_variance_0.9_4_uniform}\cref{subfig:poly_vel_variance_0.9_4_uniform}, $N_{PC}=5$ \cref{subfig:vel_variance_0.9_5_uniform}\cref{subfig:poly_vel_variance_0.9_5_uniform}.}

    % \caption{SSFEM results on a quadrilateral grid of $1275$ nodes with $\mu\sim\vectspace{U}(0.845, 0.955)$. Variance magnitude on the left in Figures \cref{subfig:vel_variance_0.9_2_uniform}\cref{subfig:vel_variance_0.9_3_uniform}\cref{subfig:vel_variance_0.9_4_uniform}\cref{subfig:vel_variance_0.9_5_uniform}, and solution at $\bar{\boldsymbol{x}}$ on the right in Figures \cref{subfig:poly_vel_variance_0.9_2_uniform}\cref{subfig:poly_vel_variance_0.9_3_uniform}\cref{subfig:poly_vel_variance_0.9_4_uniform}\cref{subfig:poly_vel_variance_0.9_5_uniform} respectively for $N_{PC}$ taking the values $3,4,5,6$.}
    \label{fig:vel_variance_0.9_uniform}
\end{figure}

As we can see from Figure \ref{fig:vel_variance_0.9_uniform}, also for such viscosity distribution the variance of the vertical velocity is concentrated in the post-inlet region.  
Moreover, for polynomials of degree greater than one multiple local extrema appear in the sampling region, and for increasing degrees the peaks stabilize their values around the three deterministic branches for the mean value of the $\mu$ parameter (see Figure \ref{subfig:poly_vel_variance_0.9_5_uniform}), respectively $\{-2,0,2\}$.

\begin{figure}[th]
\centering
\captionsetup[sub][figure]{skip=-98pt,slc=off,margin={0pt,0pt}}
    \subcaptionbox{\label{fig:kde_uniform}}{\includegraphics[width=0.44\textwidth]{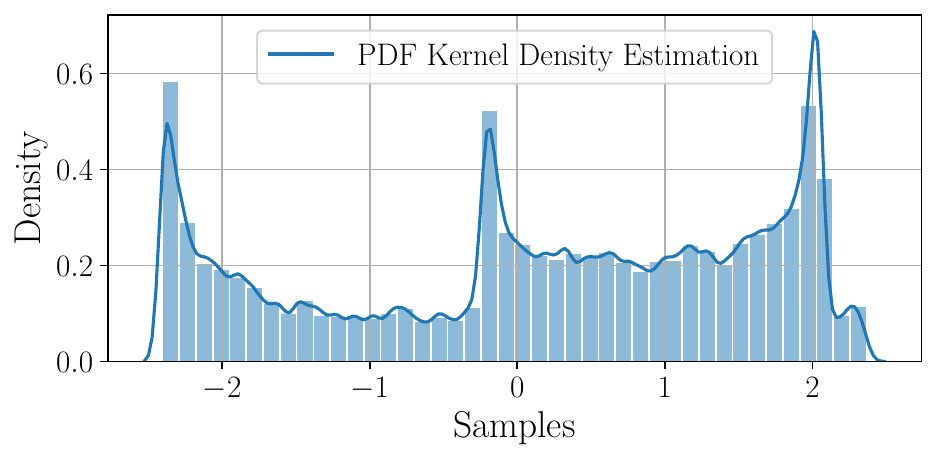}
    }
    \hfill
    \subcaptionbox{\label{fig:kde_nonsymm}}{\includegraphics[width=0.44\textwidth]{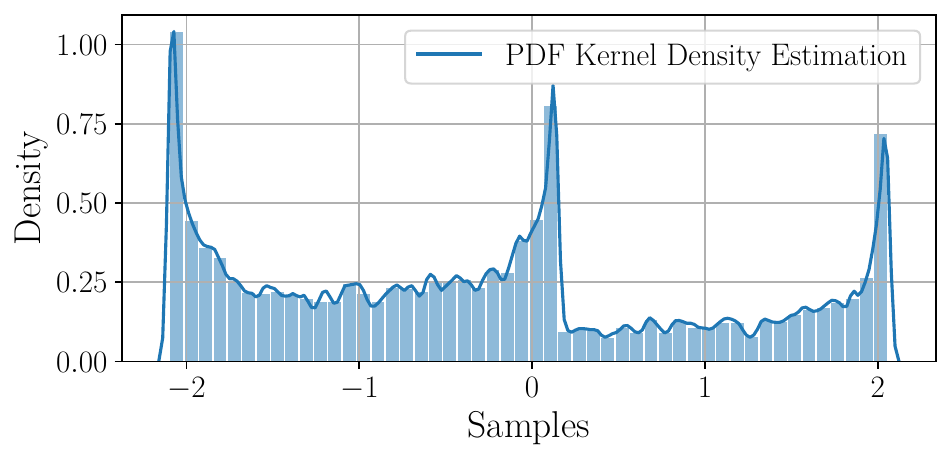}
    }
    \captionsetup{subrefformat=parens}
 
\caption{Kernel Density Estimation of the velocity PDF at the domain point of maximum variance for $\mu\sim\vectspace{U}(0.845,0.955)$, for non-symmetric meshes of \cref{fig:kde_uniform} $1275$  and \cref{fig:kde_nonsymm} $1541$ nodes.}
\end{figure}

Sampling the polynomials according to $\overline{\xi}$ leads to Figure \ref{fig:kde_uniform}, where three evident peaks are visible, corresponding precisely to the three extrema of the polynomial. %We remark that a minor peak seems also to be present for a viscosity value larger than $2$. Nevertheless, this is clearly negligible with respect to the others, as it derives from the maximum of Figure \ref{subfig:poly_vel_variance_0.9_5_uniform} barely appearing at the left margin of the sampling region. 
Thus, by assigning a uniform distribution to $\mu$, which is not strictly peaked on the mean value as in the Gaussian case, Legendre polynomials allow inferring the unstable branch too, capturing a larger variance of the solution.

\begin{figure}[tb!]
    \captionsetup[sub][figure]{skip=-68pt,slc=off,margin={0pt,0pt}}
    \subcaptionbox{\label{subfig:variance_vel_0.9_nonsymm}}{\includegraphics[trim={12cm 9.3cm 9.6cm 11.5cm}, clip, width=0.65\textwidth]{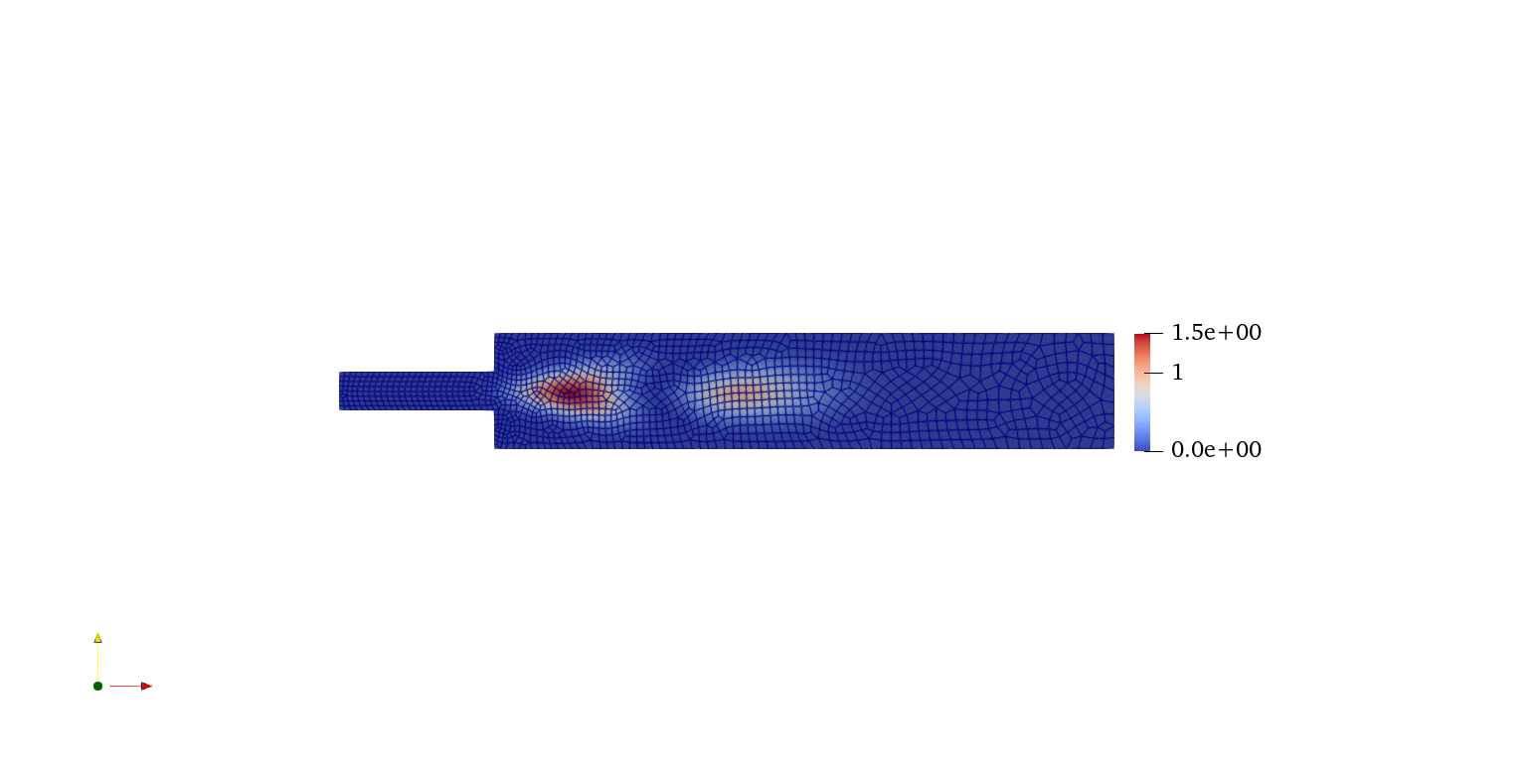}
    }
    \hfill
    \subcaptionbox{\label{subfig:poly_vel_0.9_nonsymm}}{\includegraphics[trim={0cm 0cm 0cm 0cm}, clip, width=0.3\textwidth]{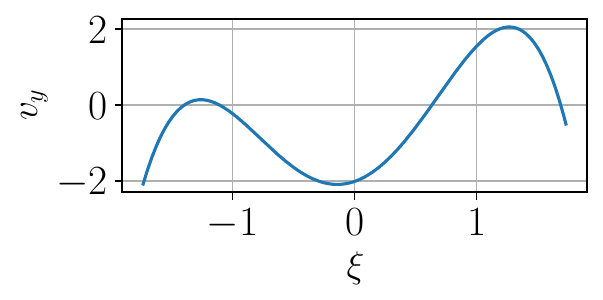}
    }
    \captionsetup{subrefformat=parens}
    
    \caption{SSFEM results on a quadrilateral grid of $1541$ nodes with $\mu\sim\vectspace{U}(0.845, 0.955)$: \cref{subfig:variance_vel_0.9_nonsymm} variance magnitude, and \cref{subfig:poly_vel_0.9_nonsymm} solution at $\bar{\boldsymbol{x}}$ for $N_{PC}=4$.}
    \label{fig:vel_09_nonsymm}
    \end{figure} 

Furthermore, we find that dealing with more dense grids as the one displayed in Figure \ref{fig:vel_09_nonsymm} can bring advantages in branch detection, requiring lower degree polynomials. 
Indeed, as it can be seen in Figure \ref{fig:kde_nonsymm}, the reconstructed PDF corresponding to the dense unstructured grid yields the same information than \ref{fig:kde_uniform}, but using polynomials of degree $4$ instead of degree $5$.

\begin{figure}[bht]
    \centering

    %------------------ FIGURA ------------------
    \begin{minipage}[c]{0.47\textwidth}
        \centering
        \includegraphics[width=.9\linewidth]{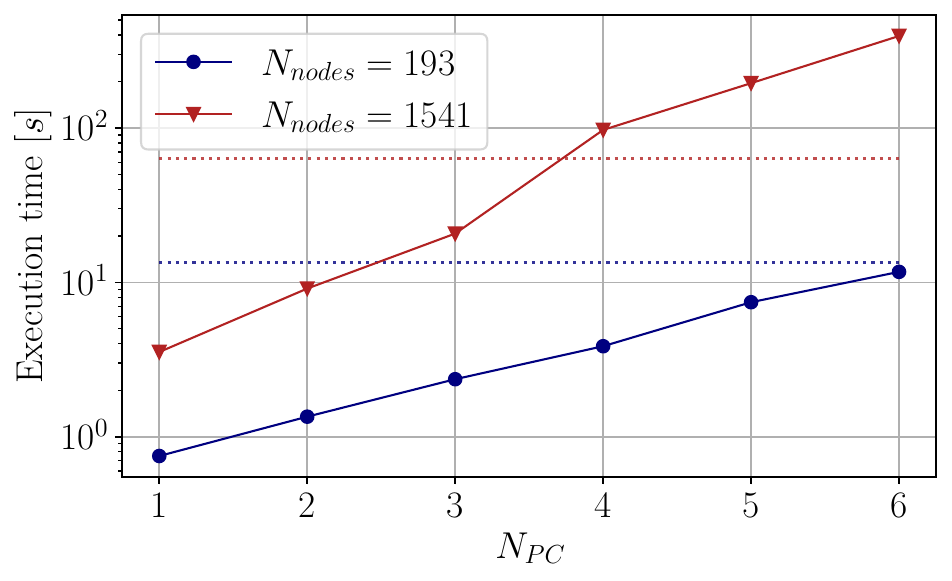}
        \captionsetup{width=\linewidth}
        \captionof{figure}{{Plot of the scaling of the execution time for coarse and fine grids, with \(\mu \sim \mathcal{U}(0.845,\,0.955)\).}}
        \label{fig:scaling}
    \end{minipage}
    \hfill
    %------------------ TABELLA ------------------
    \begin{minipage}[c]{0.47\textwidth}
        \centering
        \vspace{-0.2cm}
        \vspace{0.5em}
        \begin{tabular}{c c c} 
            \hline
            % \rowcolor{mplblue!20}
            \(N_{\mathrm{PC}}\) & $N_{\text{nodes}}=193$ & $N_{\text{nodes}}=1541$ \\
            \hline
            1 & $0.751s$ & $3.54s$ \\
            2 & $1.35s$ & $9.12s$ \\
            3 & $2.36s$ & $20.7s$ \\
            4 & $3.86s$ & $97.1s$ \\
            5 & $7.44s$ & $195s$ \\
            6 & $11.7s$ & $394s$ \\
            \hline
        \end{tabular}
        \vspace{.6cm}
        \captionsetup{width=\linewidth}
        \captionof{table}{{Execution times for increasing \(N_{\mathrm{PC}}\) for coarse and fine grids, with \(\mu \sim \mathcal{U}(0.845,\,0.955)\).}}
        \label{tab:scaling}
    \end{minipage}

\end{figure}

\begin{figure}[b]
\centering
\captionsetup[sub][figure]{skip=-78pt,slc=off,margin={0pt,0pt}}
        \subcaptionbox{\label{subfig:vel_symm_variance}}{\includegraphics[trim={12cm 7.3cm 7.6cm 11.5cm}, clip, width=0.66\textwidth]{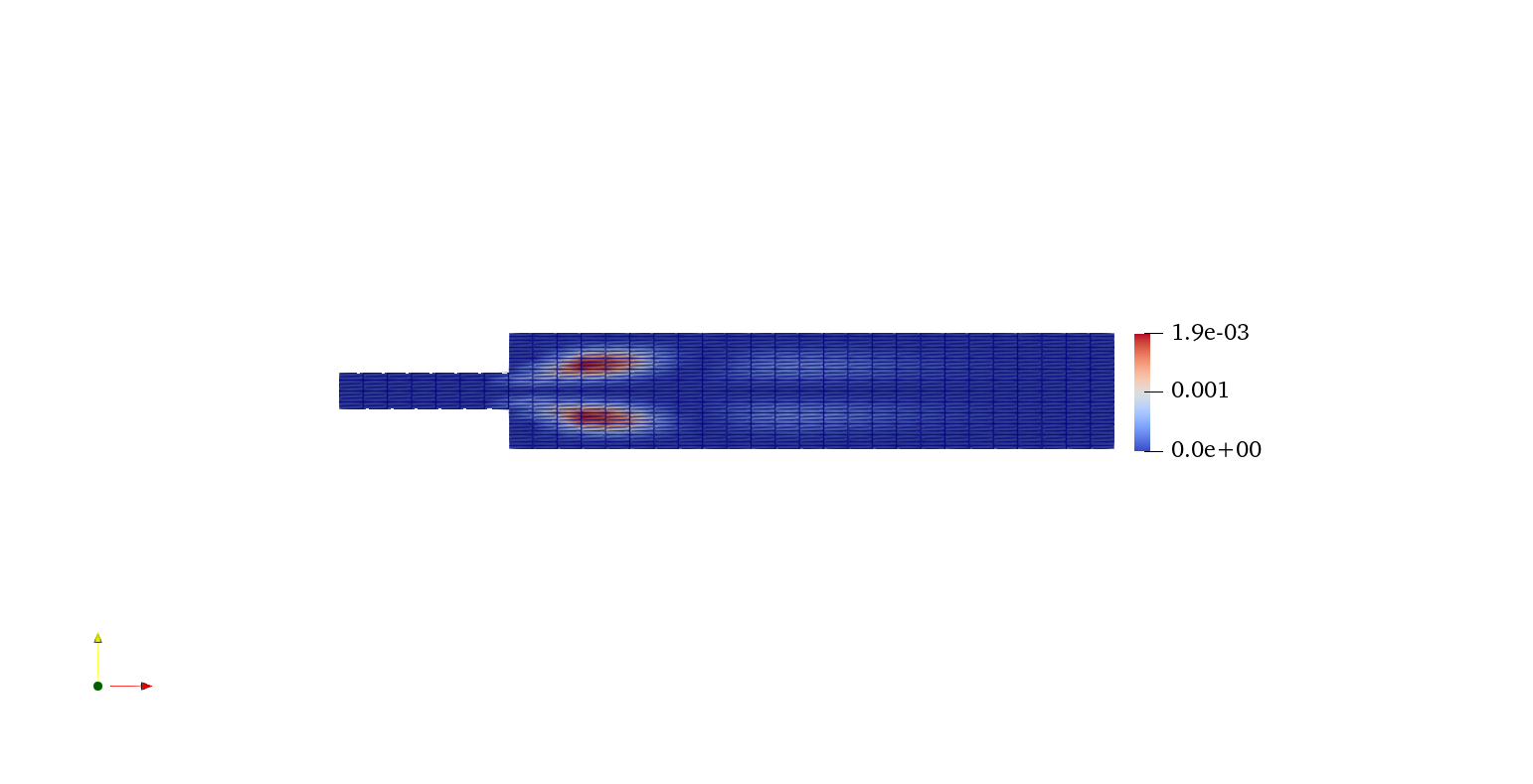}
    }
    \hfill
    \captionsetup[sub][figure]{skip=-68pt,slc=off,margin={0pt,0pt}}
    \subcaptionbox{\label{subfig:poly_vel_symm}}{\includegraphics[trim={0cm 0.cm 0cm 0cm}, clip, width=0.3\textwidth]{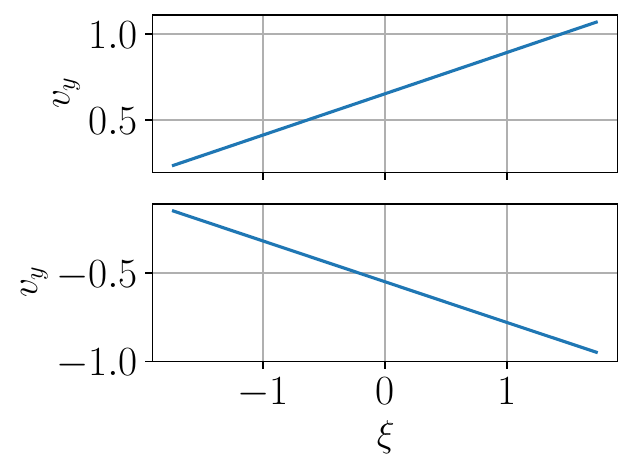}
    }
    \captionsetup{subrefformat=parens}

    \caption{SSFEM results on a quadrilateral symmetric grid of $935$ nodes with $\mu\sim\vectspace{U}(0.845, 0.955)$: \cref{subfig:vel_symm_variance} variance magnitude, and \cref{subfig:poly_vel_symm}
    solution at $\bar{\boldsymbol{x}}_1$ and $\bar{\boldsymbol{x}}_2$ for $N_{PC}=4$.}
    \label{fig:vel_symm}
\end{figure}

{Nevertheless, denser grids entail a higher computational cost, which increases with the number of PC modes $N_{PC}$, as shown in Fig.~\ref{fig:scaling} and Table \ref{tab:scaling}. 
Results are reported for coarse and fine grids, and all timings refer to a serial implementation. 
The growth of the computational time with increasing $N_{PC}$ is due both to the larger number of degrees of freedom, leading to more expensive evaluations of the residual and the Jacobian, and to the increasing difficulty of Newton convergence for larger coupled systems.

Although the execution time increases with both $N_{PC}$ and the number of grid nodes, it is essential to compare these costs with those required by standard continuation-based approaches. To this end, in Figure \ref{fig:scaling} we also report horizontal dotted lines representing an estimate of the computational cost of classical continuation methods. These lines are obtained by multiplying the execution time of a single deterministic FEM solve ($N_{PC}=1$) by the number of bifurcation branches (three) and by the number of continuation steps, assuming to start at the bifurcation point $\mu^{*}\simeq 0.96$ and taking $\Delta\mu = 0.01$, as done to construct Figure~\ref{fig:bif-diag-deterministic}.  

Under this comparison, the SSFEM approach exhibits a clear computational advantage: it outperforms classical continuation methods for all values of $N_{PC}$ on the coarse grid and up to $N_{PC}=3$ on the finer grid. Moreover, this comparison is in fact unfavorable to our method, since it neglects additional computational costs inherent to standard continuation techniques, such as those associated with deflation strategies \cite{FarrellDeflationTechniquesFinding2015,PichiDeflationbasedCertifiedGreedy2025a}.  

Despite this, the proposed approach already outperforms a lower-bound estimate of classical methods. Its relative advantage is expected to further increase away from the bifurcation point, as illustrated in Figure~\ref{fig:scaling_far}.
}
As a final consideration, in Figure \ref{fig:vel_symm} we show the case of a highly structured symmetric grid. 
One can note that the variance plot is significantly different, since it reports a higher variance not on the horizontal symmetric axis of the inlet region, but instead on two symmetric regions at the top and bottom of it. 
In particular, the polynomials in the points of maximum variance are linear, and no local extremum is visible, suggesting that no bifurcating behavior is detected. 
This is due to the symmetry of the mesh, and it is in line with what has already been observed in deterministic bifurcation detection methods by previous works \cite{HessLocalizedReducedorderModeling2019,PichiPhd}.

\subsection{Deterministic bifurcation diagram inference}

In this final section, we extend the previous analysis to different viscosity mean values.

\begin{figure}[t]
\centering

\captionsetup[sub][figure]{skip=-78pt,slc=off,margin={0pt,0pt}}
        \subcaptionbox{\label{subfig:vel_13_varianceh}}{\includegraphics[trim={12cm 7.3cm 7.6cm 11.5cm}, clip, width=0.66\textwidth]{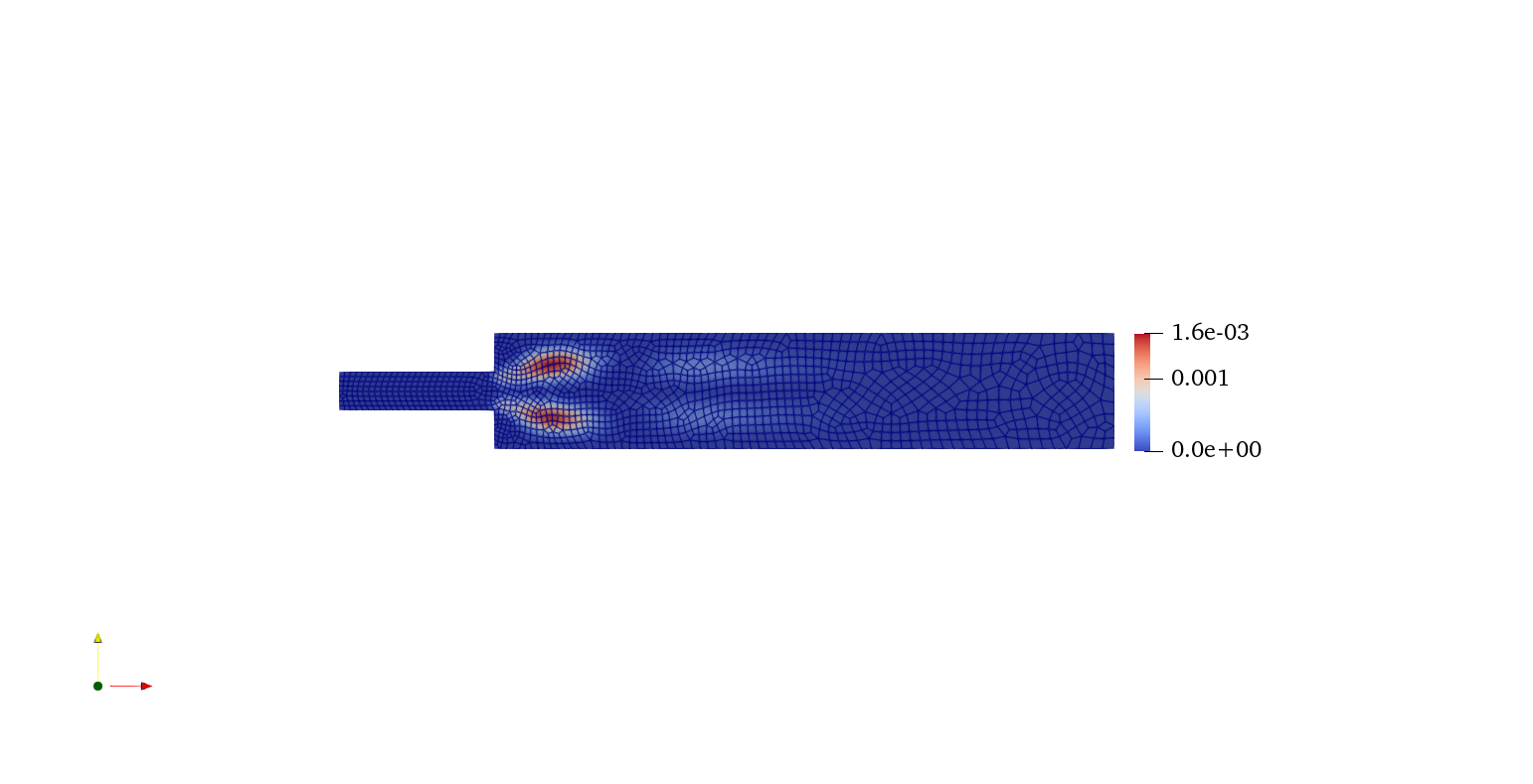}
    }
    \hfill
    \captionsetup[sub][figure]{skip=-68pt,slc=off,margin={0pt,0pt}}
    \subcaptionbox{\label{subfig:poly_vel_13h}}{\includegraphics[trim={0cm 0.cm 0cm 0cm}, clip, width=0.3\textwidth]{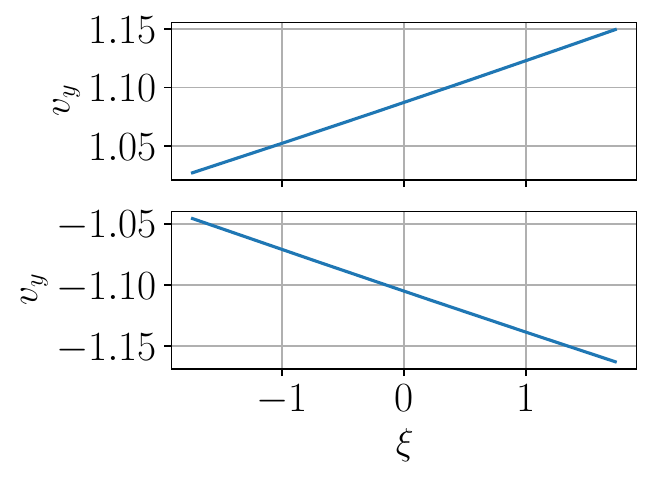}
    }
    \captionsetup{subrefformat=parens}

    \caption{SSFEM results on a quadrilateral grid of $1541$ nodes with $\mu\sim\vectspace{U}(1.245, 1.355)$: \cref{subfig:vel_13_varianceh} variance magnitude, and \cref{subfig:poly_vel_13h} solution at $\bar{\boldsymbol{x}}_1$ and $\bar{\boldsymbol{x}}_2$.}
    \label{fig:vel_13h}
\end{figure}

First, we are interested in checking the behavior of the velocity variance in the uniqueness regime. 
In Figure \ref{fig:vel_13h} the SSFEM results are shown for $\mu\sim\vectspace{U}(1.245, 1.355)$, thus quite far from the bifurcation point. 
When the uniqueness of the deterministic solution is ensured, the variance is significantly lower (about three orders of magnitude) than the one observed in the non-uniqueness regime. 
Moreover, the polynomial solution at the domain points of highest variance is linear, confirming the absence of the bifurcating behavior.

\begin{figure}[t]

    \captionsetup[sub][figure]{skip=-68pt,slc=off,margin={0pt,0pt}}
    \subcaptionbox{\label{subfig:vel_1_variance}}{\includegraphics[trim={12cm 9.3cm 9.6cm 11.5cm}, clip, width=0.65\textwidth]{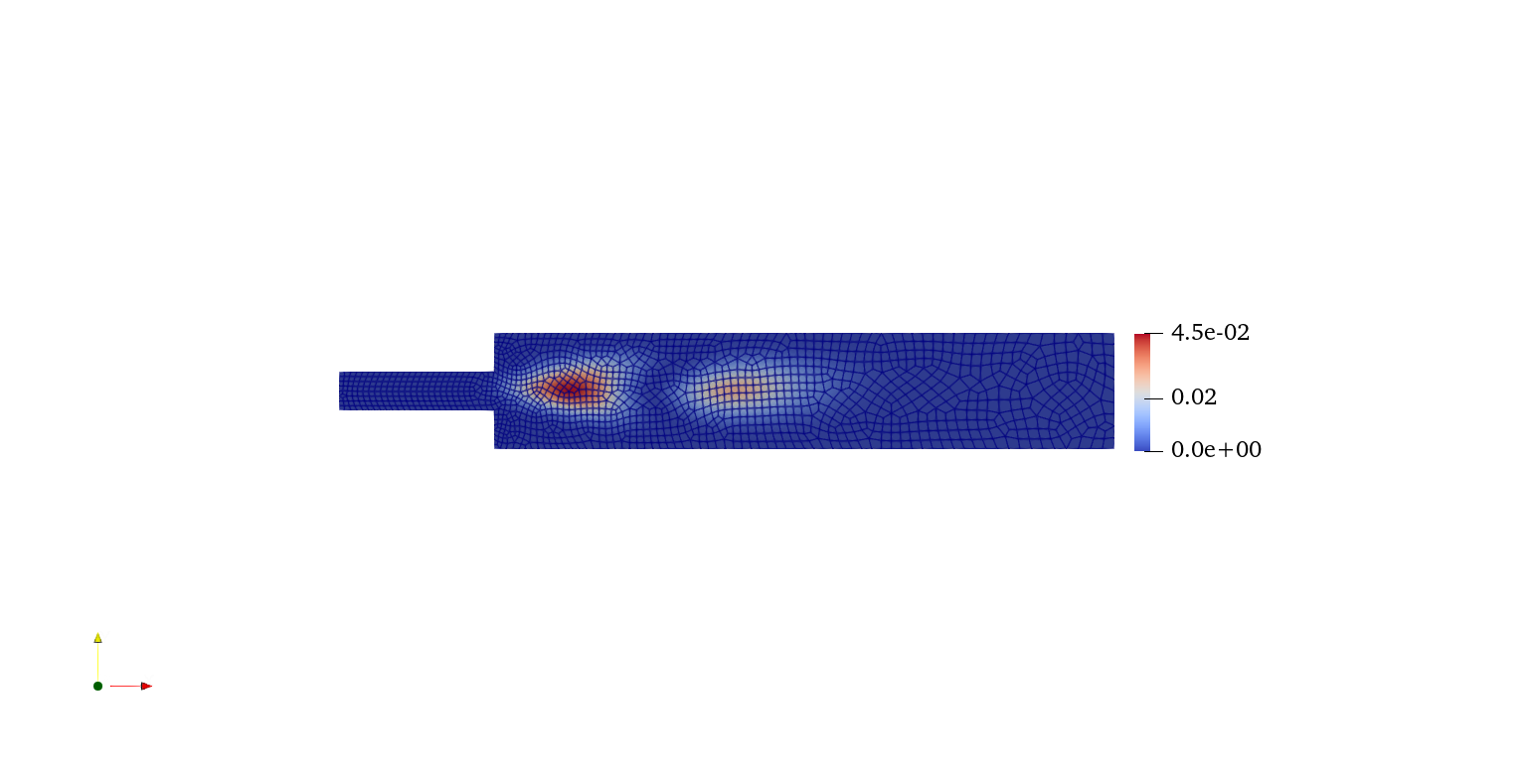}
    }
    \hfill
    \subcaptionbox{\label{subfig:vel_1_poly}}{\includegraphics[trim={0cm 0cm 0cm 0cm}, clip, width=0.3\textwidth]{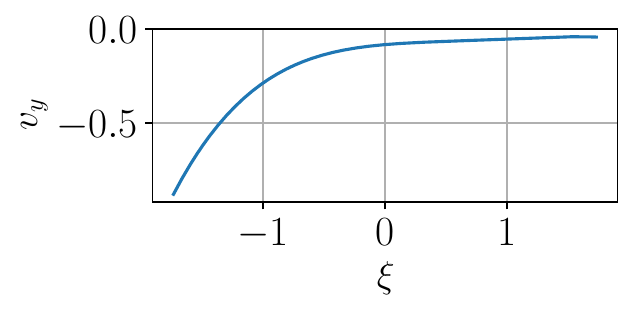}
    }
    \captionsetup{subrefformat=parens}
    
    \caption{SSFEM results on a quadrilateral grid of $1541$ nodes with $\mu\sim\vectspace{U}(0.945, 1.055)$. Variance magnitude \cref{subfig:vel_1_variance}, and solution at $\bar{\boldsymbol{x}}$ \cref{subfig:vel_1_poly}.}
    \label{fig:vel_1}
\end{figure}

When moving the viscosity value towards the bifurcation point, as expected, the variance increases as shown in Figure \ref{fig:vel_1} for $\mu\sim\vectspace{U}(0.945, 1.055)$. 
Nevertheless, in both cases, the corresponding polynomials, evaluated in the domain points of maximum variance, do not present multiple local extrema.  

\begin{figure}[hbt!]
\centering

    \captionsetup[sub][figure]{skip=-68pt,slc=off,margin={0pt,0pt}}
    \subcaptionbox{\label{subfig:variance_vel_0.8}}{\includegraphics[trim={12cm 9.3cm 9.6cm 11.5cm}, clip, width=0.65\textwidth]{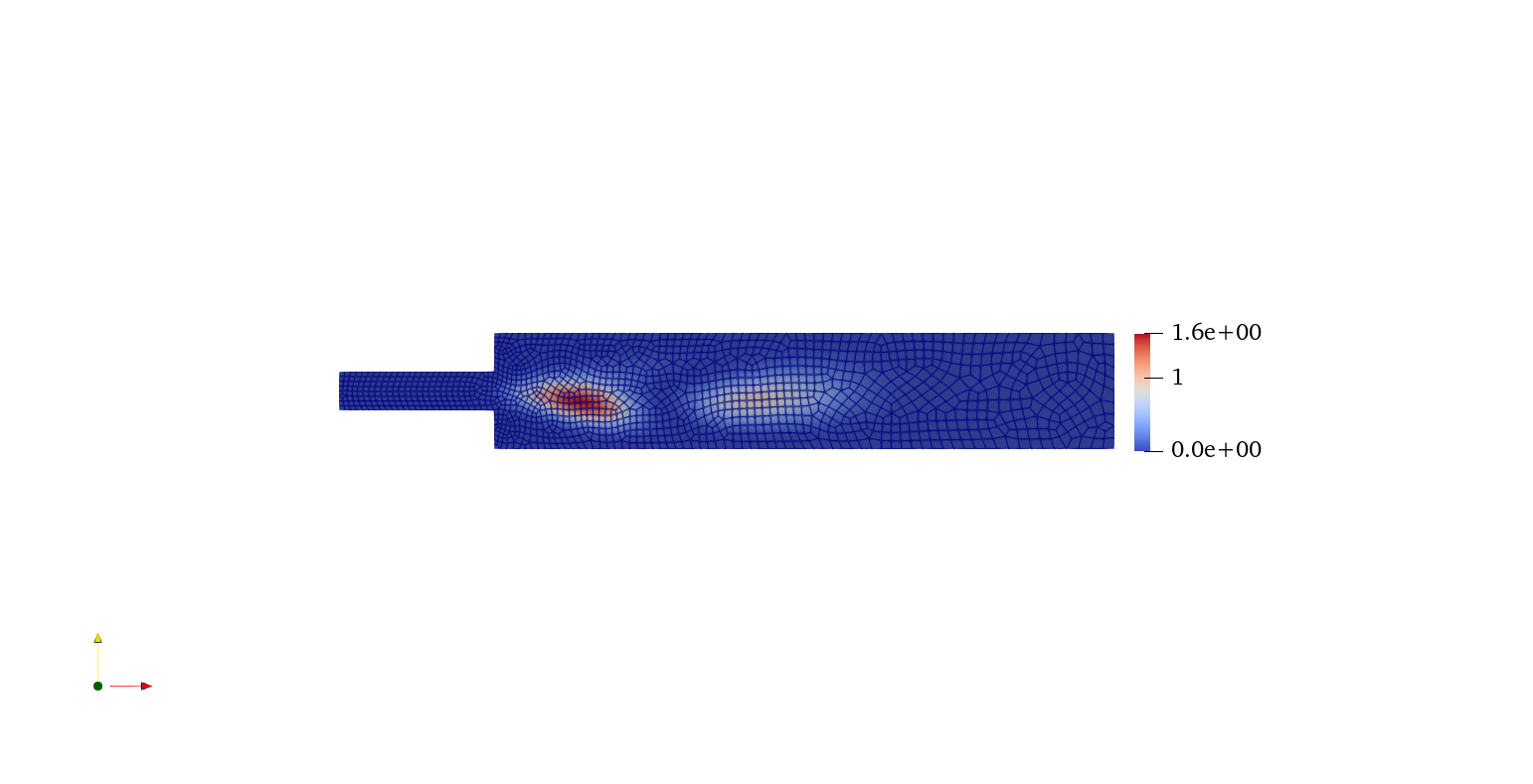}
    }
    \hfill
    \subcaptionbox{\label{subfig:poly_vel_0.8}}{\includegraphics[trim={0cm 0cm 0cm 0cm}, clip, width=0.3\textwidth]{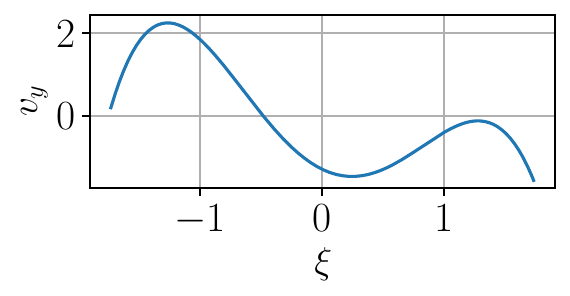}
    }

    \subcaptionbox{\label{subfig:variance_vel_0.7}}{\includegraphics[trim={12cm 9.3cm 9.6cm 11.5cm}, clip, width=0.65\textwidth]{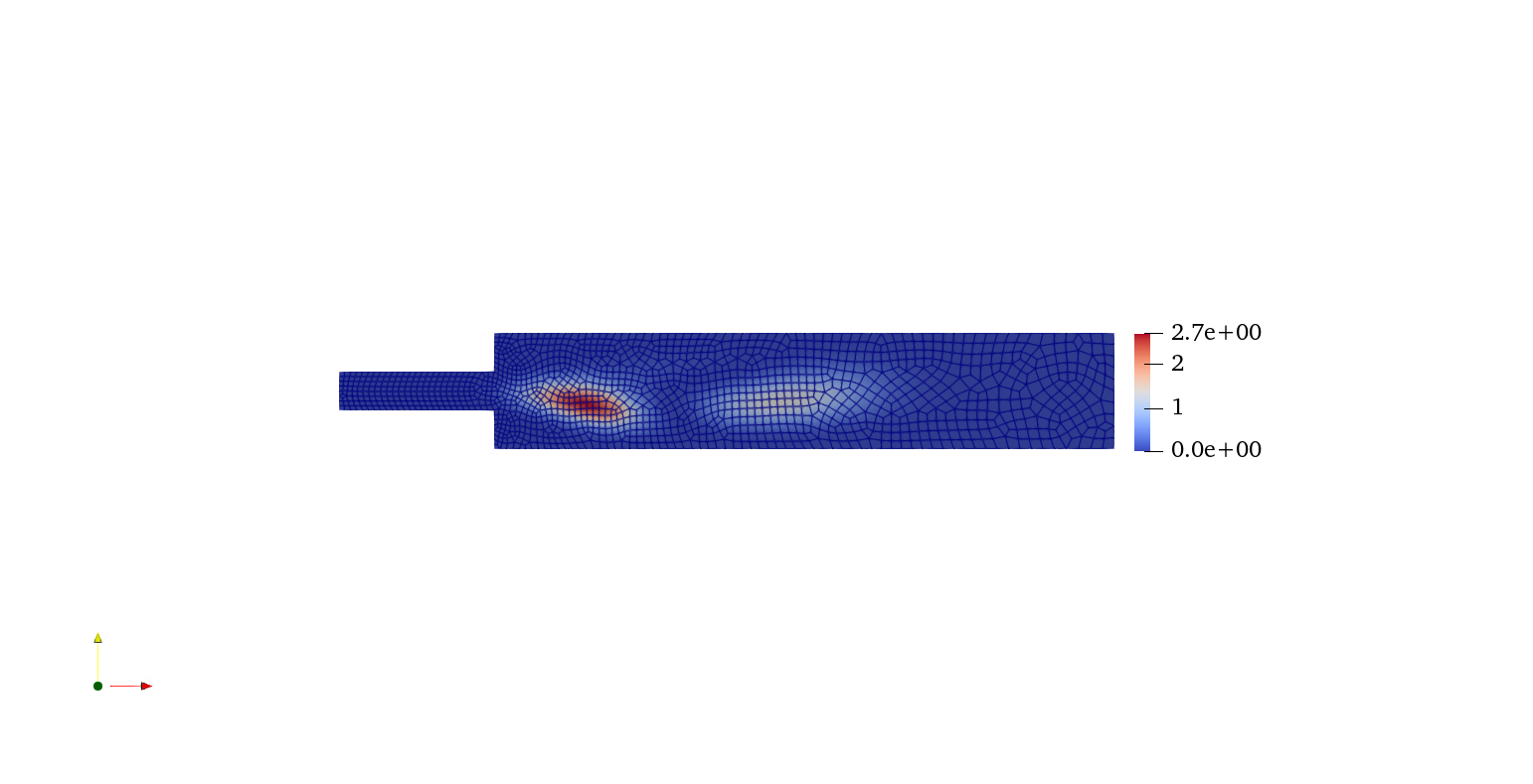}
    }
    \hfill
    \subcaptionbox{\label{subfig:poly_vel_0.7}}{\includegraphics[trim={0cm 0cm 0cm 0cm}, clip, width=0.3\textwidth]{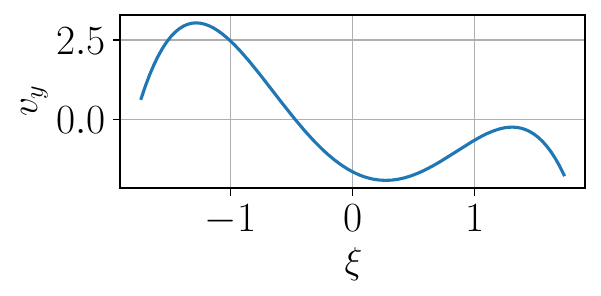}
    }
    
    \captionsetup{subrefformat=parens}
    
    \caption{SSFEM results on a grid of $1541$ nodes: variance magnitude and solution at $\bar{\boldsymbol{x}}$, for 
    $\mu\sim\vectspace{U}(0.745, 0.855)$ in \cref{subfig:variance_vel_0.8}\cref{subfig:poly_vel_0.8} and $\mu\sim\vectspace{U}(0.645, 0.755)$ in \cref{subfig:variance_vel_0.7}\cref{subfig:poly_vel_0.7}.}
    \label{fig:vel_variance_nonunique_solution}
\end{figure}

Towards a complete investigation of the bifurcation diagram, we are interested in reconstructing the branches for any value of viscosity, without relying on any prior information. 
To this end, we report the results of the SSFEM computations for $\mu\sim\vectspace{U}(0.745,0.855)$ and $\mu\sim\vectspace{N}(0.645,0.755)$ in Figure \ref{fig:vel_variance_nonunique_solution}.
The magnitude of the variance increases with respect to Figure \ref{fig:vel_09_nonsymm}, and progressively becomes higher as the viscosity mean decreases. 
This is consistent with the progressive separation of the branches already in the deterministic bifurcation diagram.
The correspondent PDFs and are shown in Figures \ref{subfig:poly_vel_0.8} and \ref{subfig:poly_vel_0.7}.

\begin{figure}[t]
    \centering

    %------------------ FIGURA ------------------
    \begin{minipage}[c]{0.47\textwidth}
        \centering
        \includegraphics[width=.9\linewidth]{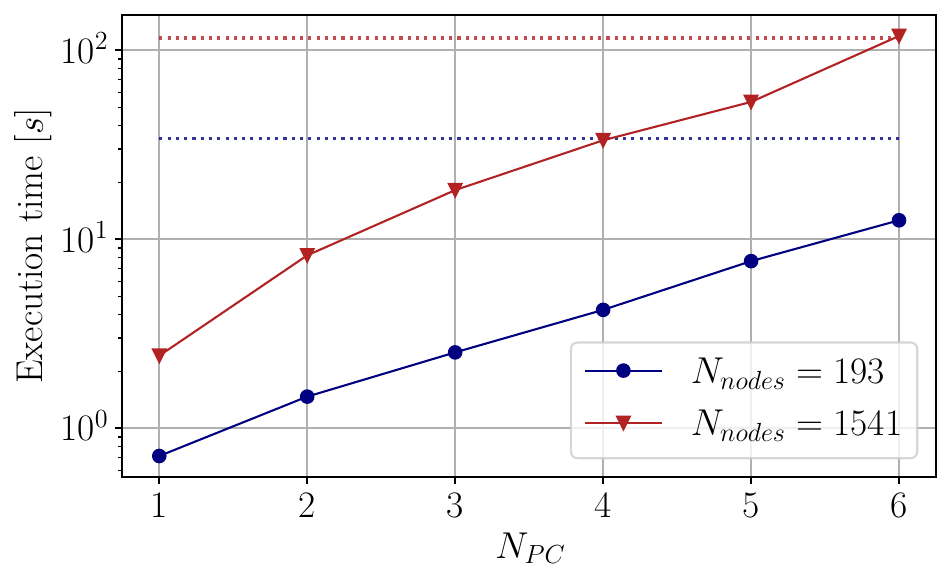}
        \captionsetup{width=\linewidth}
        \captionof{figure}{{Plot of the scaling of the execution time for coarse and fine grids, with $\mu \sim \mathcal{U}(0.745,\,0.855)$.}}
        \label{fig:scaling_far}
    \end{minipage}
    \hfill
    %------------------ TABELLA ------------------
    \begin{minipage}[c]{0.47\textwidth}
        \centering
        \vspace{-0.2cm}
        \vspace{0.5em}
        \begin{tabular}{c c c} 
            \hline
            % \rowcolor{mplblue!20}
            \(N_{\mathrm{PC}}\) & $N_{\text{nodes}}=193$ & $N_{\text{nodes}}=1541$ \\
            \hline
            1 & $0.713s$ & $2.42s$ \\
            2 & $1.47s$ & $8.22s$ \\
            3 & $2.52s$ & $18.2s$ \\
            4 & $4.23s$ & $33.4s$ \\
            5 & $7.66s$ & $53.3s$ \\
            6 & $12.6s$ & $119s$ \\
            \hline
        \end{tabular}
        \vspace{.6cm}
        \captionsetup{width=\linewidth}
        \captionof{table}{{Execution times for increasing $N_{\mathrm{PC}}$ for coarse and fine grids, with $\mu \sim \mathcal{U}(0.745,\,0.855)$.}}
        \label{tab:scaling_far}
    \end{minipage}

\end{figure}

{
Moreover, Figure~\ref{fig:scaling_far} and Table~\ref{tab:scaling_far} show that the SSFEM approach significantly outperforms traditional continuation methods for all tested values of $N_{PC}$ and for both computational grids when considering $\mu \sim \mathcal{U}(0.745,\,0.855)$, i.e., far from the bifurcation point but still in the non-uniqueness regime. 

Indeed, while the computational cost of continuation methods increases as the distance from the bifurcation point grows, due to the need to track solution branches over a wider parameter range, the cost of the SSFEM approach is not affected. On the contrary, a reduction in solver time is observed with respect to Table \ref{tab:scaling}, likely due to the fact that the problem becomes numerically easier once the system operates far from the bifurcation point.}

By looking at the polynomials for lower viscosity values, we observe that SSFEM is able to detect three local extrema in the sampling region. 
Moreover, their values follow the peaks identified by the trend of the branches in the deterministic diagram. 
Nevertheless, given the difficulties highlighted before, the symmetry between positive and negative values for the two stable branches is not exactly recovered with this mesh, as the wall-hugging bottom flow is identified but not accurately reached in magnitude. 
Nevertheless, the unstable symmetric branch seems to be always recovered.

\begin{figure}[ht]
    \centering
    \includegraphics[width=0.5\textwidth]{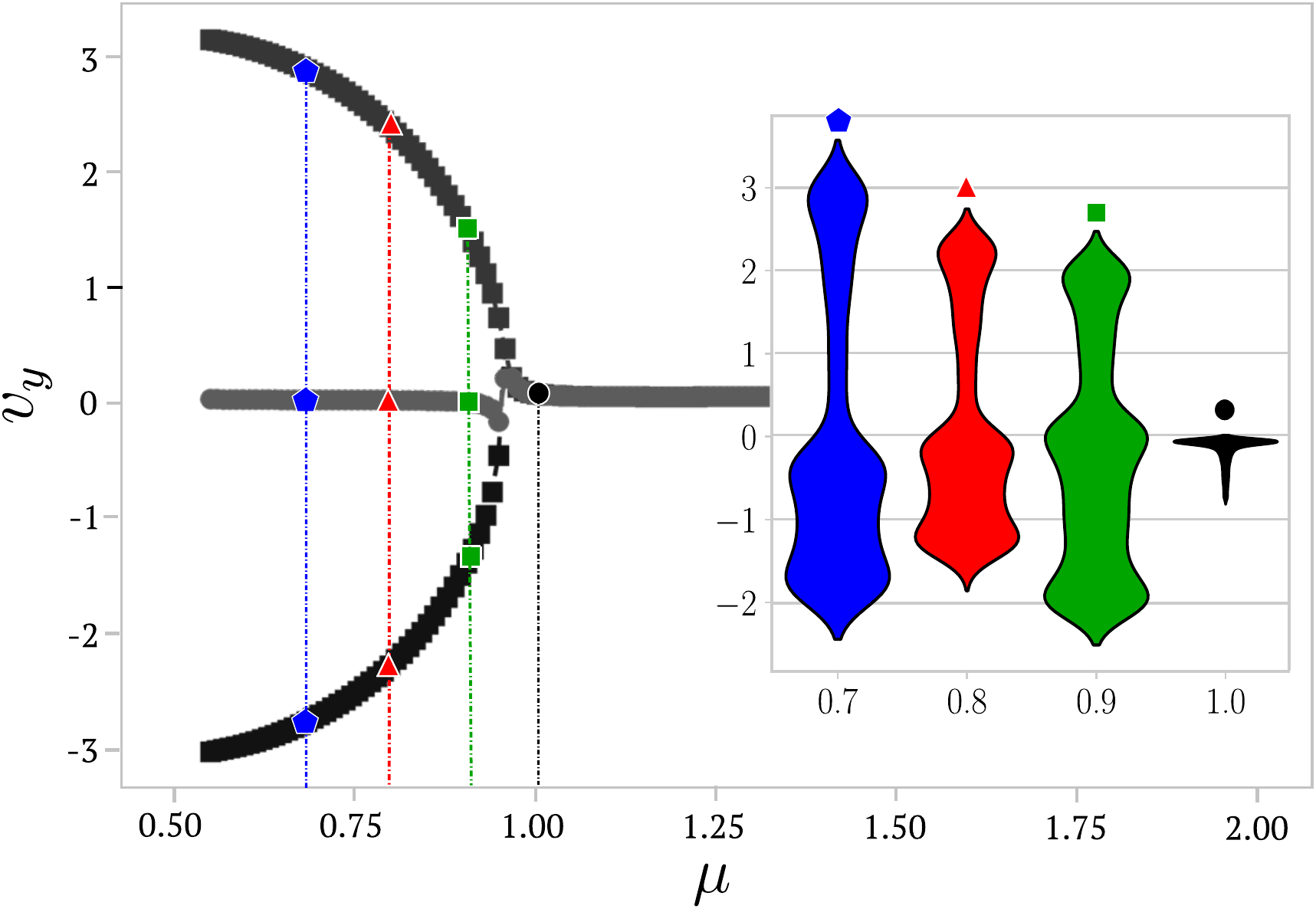}
    \caption{SSFEM probabilistic bifurcation diagram for different viscosity values.}
    \label{fig:SSFEM_bifplot}
\end{figure}

We can exploit the results obtained to visualize a probabilistic version of the bifurcation diagram. 
In Figure \ref{fig:SSFEM_bifplot} we report this comparison, highlighting with different markers the admissible solutions for different viscosity parameters. 
We observe that the branches of the solution are exactly recovered in the case of $\mu=0.9$, and that for $\mu=1.0$ the density of the solution is concentrated on a single value, corresponding to the unique solution. 
Even when considering viscosity values in the ``strong" bifurcation regime, namely $\mu=0.7$ and $\mu=0.8$, the evolution of the branches is followed, although the lower branch is not perfectly recovered in its symmetry. 
{We remark that, given the observed influence of grid density and PC maximum degree in the approximation accuracy, in order to ensure the quality of recovered branches an important task is to carefully design those quantities.}
Still, we conclude that for reasonable mesh density and PC degrees, also in this complex benchmark, one can retrieve significant information on the bifurcation diagram.

% This analysis therefore confirms the potentialities of investigating complex bifurcating phenomena under the stochastic point of view via probabilistic-aware surrogate approaches, as a way of inferring features of the deterministic model through stochastic perturbations. 
% Indeed, this approach is completely agnostic of the bifurcating behavior of the system and can be used for both bifurcation point and branch detection instead of continuation-based methods. 
% Moreover, it also opens a research path to address the performance of the PC expansion accuracy w.r.t.\ the properties of the mesh and the type/order of the polynomials. 

\section{Conclusions and Perspectives}

In this work, we propose a stochastic-perturbation approach to deterministic bifurcation problems. 
We implement an intrusive pipeline based on SSFEM to retrieve the PC coefficients of the probabilistic solution, and we provide computational evidence that it successfully encodes accurate information about the deterministic bifurcation diagram. 
In particular, we find that big perturbations are useful for parameter space exploration, while smaller ones can serve to identify the branches for a precise parameter value.

Concerning the PC surrogate model of a stochastic bifurcating PDE, we show the strict relation between the branches of coexisting solutions and the peaks of the multimodal probability density function obtained by sampling the PC solution. 
Indeed, despite the non-uniqueness of the PC solution in a bifurcating regime, the solution's PDF appears to be uniquely determined and captured, and it demonstrates a correspondence with the solution branches for the specific parameter value at which its distribution is centered. 
That is actually because the PC polynomials' local extrema tend to concentrate around the values of the branches, leading to the formation of clusters of points around them when the PC solution is sampled.

Therefore, we show that such an intrusive stochastic approach overcomes an important issue arising in the computational bifurcating context, as it does not require prior knowledge about the solution branches to give back valid information about the solution manifold. 
Indeed, with classical continuation methods, retrieving samples of the non-unique parametric solution requires some prior information about the bifurcating behavior.
This precaution is essential to allow the numerical solver to converge to different branches, and thus to provide a complete overview of the solution manifold via a meaningful and representative ensemble of admissible configurations. 
As a result, SSFEM offers a significant advantage by being independent of problem-specific initialization, while still providing a precise representation of the complete solution manifold. 
In particular, in this work, we address the computational study of its performance in the case of the pitchfork bifurcation. 
We refer to future works for additional investigations on the method's behavior and reliability with different types of bifurcations.

On the other hand, from a strictly computational point of view, the SSFEM intrusive approach is more expensive then a single FEM simulation, due to the dimension of the system to be solved. 
{This poses new computational challenges, which can be mitigated through ad-hoc convergence strategies and preconditioning methods, as shown in previous works also addressing the Navier-Stokes equations \cite{sousedik2016, Lee2019}. }
Therefore, significant improvements could be achieved, enabling faster computation of higher order coefficients and allowing for the analysis of more complex behaviors, such as the case of higher Reynolds numbers with successive bifurcations in the Navier-Stokes setting. 

{A topic worthy of further investigation concerns the role of mesh symmetry and density in the detection of bifurcations when using the SSFEM approach.}
Moreover, additional investigations could focus on the source of stochasticity in the system. 
For instance, random initial and/or boundary conditions may play an important role in the evolution of the system for time-dependent problems, and it could be used to infer a probabilistic overview of the model's behavior in time. 
Further tests can also be conducted to apply the perturbation method to other types of bifurcations or to couple it with traditional continuation methods. 
In this way, one could extract ad-hoc initializations for the numerical solver to enable a precise evaluation of branches in a larger portion of the parameter space, even far from the bifurcation point.

Finally, given the reliability of the PC surrogate model for the probability distribution of the solution, it is worth exploring a Bayesian approach to the bifurcation problem starting from this prior assumption, possibly following the framework proposed in \cite{Girolami2021} and adapting it to the polynomial basis. 
Indeed, incorporating real data into the analysis could be particularly important to model the complex behavior of a real-world parametric system, as it could bring the model's improvements in its reality reproduction capability.

\section*{Declarations}

\textbf{Funding} This work has been conducted within the research activities of the consortium iNEST (Interconnected North-East Innovation Ecosystem) and its CC5 Young Researchers initiative, Piano Nazionale di Ripresa e Resilienza (PNRR) – Missione 4 Componente 2, Investimento 1.5 – D.D. 1058 23/06/2022, ECS00000043 – CUP G93C22000610007, supported by the European Union's NextGenerationEU program.\\
%\section*{Acknowledgment}
\textbf{Acknowledgment}
The authors thank Pasquale Claudio Africa for the fruitful discussions and fundamental support with the deal.II implementation of the method, and acknowledge the support by Gruppo Nazionale di Calcolo Scientifico (INdAM-GNCS).

%\noindent\textbf{Conflict of interest} The authors declare that they have no conflict of interest.\\

%\noindent\textbf{Data availability} The code to generate the data and results are available at \url{https://github.com/ICGonnella/SSFEM-Coanda-Effect.git}.\\

\bibliographystyle{abbrv}
\bibliography{main.bib}

@book{Seydel,
  author =        {Seydel, R.},
  publisher =     {Springer New York},
  series =        {Interdisciplinary Applied Mathematics},
  title =         {Practical Bifurcation and Stability Analysis},
  year =          {2009},
  doi =           {10.1007/978-1-4419-1740-9},
  isbn =          {9781441917393},
}

@book{allgower2003,
  author =        {Allgower, Eugene L and Georg, Kurt},
  publisher =     {SIAM},
  title =         {Introduction to numerical continuation methods},
  year =          {2003},
  doi =           {10.1137/1.9780898719154},
}

@phdthesis{PichiPhd,
  address =       {Trieste},
  author =        {Pichi, F.},
  school =        {Scuola Internazionale Superiore di Studi Avanzati},
  title =         {Reduced Order Models for parametric Bifurcation
                   Problems in Nonlinear {PDE}s},
  year =          {2020},
}

@article{CrisovanModelOrderReduction2019,
  author =        {Crisovan, R. and Torlo, D. and Abgrall, R. and
                   Tokareva, S.},
  journal =       {Journal of Computational and Applied Mathematics},
  pages =         {466--489},
  title =         {Model Order Reduction for Parametrized Nonlinear
                   Hyperbolic Problems as an Application to Uncertainty
                   Quantification},
  volume =        {348},
  year =          {2019},
  doi =           {10.1016/j.cam.2018.09.018},
}

@incollection{VenturiWeightedReducedOrder2019a,
  address =       {Cham},
  author =        {Venturi, Luca and Torlo, Davide and
                   Ballarin, Francesco and Rozza, Gianluigi},
  booktitle =     {Uncertainty Modeling for Engineering Applications},
  editor =        {Canavero, Flavio},
  pages =         {27--40},
  publisher =     {Springer International Publishing},
  title =         {Weighted {{Reduced Order Methods}} for {{Parametrized
                   Partial Differential Equations}} with {{Random
                   Inputs}}},
  year =          {2019},
  doi =           {10.1007/978-3-030-04870-9_2},
  isbn =          {978-3-030-04870-9},
}

@article{CarereWeightedPODreductionApproach2021,
  author =        {Carere, Giuseppe and Strazzullo, Maria and
                   Ballarin, Francesco and Rozza, Gianluigi and
                   Stevenson, Rob},
  journal =       {Computers \& Mathematics with Applications},
  publisher =     {Elsevier},
  title =         {A weighted {POD}-reduction approach for parametrized
                   PDE-constrained optimal control problems with random
                   inputs and applications to environmental sciences},
  volume =        {102},
  year =          {2021},
  doi =           {10.1016/j.camwa.2021.10.020},
}

@inbook{Arnold1992,
  address =       {Boston, MA},
  author =        {Arnold, Ludwig and Boxler, Petra},
  booktitle =     {Diffusion Processes and Related Problems in Analysis,
                   Volume II: Stochastic Flows},
  pages =         {241--255},
  publisher =     {Birkh{\"a}user Boston},
  title =         {Stochastic bifurcation: instructive examples in
                   dimension one},
  year =          {1992},
  doi =           {10.1007/978-1-4612-0389-6_10},
  isbn =          {978-1-4612-0389-6},
}

@article{KuehnUncertaintyQuantificationAnalysis2024,
  author =        {Kuehn, Christian and Piazzola, Chiara and
                   Ullmann, Elisabeth},
  journal =       {arXiv preprint arXiv:2404.04639},
  title =         {Uncertainty Quantification Analysis of Bifurcations
                   of the {{Allen--Cahn}} Equation with Random
                   Coefficients},
  year =          {2024},
  doi =           {10.1016/j.physd.2024.134390},
}

@article{Sudret2008,
  author =        {Sudret, Bruno},
  journal =       {Reliability Engineering \&amp; System Safety},
  month =         {July},
  number =        {7},
  pages =         {964–979},
  publisher =     {Elsevier BV},
  title =         {Global Sensitivity Analysis using {Polynomial Chaos}
                   Expansions},
  volume =        {93},
  year =          {2008},
  doi =           {10.1016/j.ress.2007.04.002},
  issn =          {0951-8320},
}

@article{Xiu2003,
  author =        {Xiu, Dongbin and Karniadakis, George Em},
  journal =       {Journal of Computational Physics},
  month =         {May},
  number =        {1},
  pages =         {137–167},
  publisher =     {Elsevier BV},
  title =         {Modeling Uncertainty in Flow Simulations via
                   generalized {Polynomial Chaos}},
  volume =        {187},
  year =          {2003},
  doi =           {10.1016/S0021-9991(03)00092-5},
  issn =          {0021-9991},
}

@article{Stefanou2009,
  author =        {Stefanou, George},
  journal =       {Computer Methods in Applied Mechanics and
                   Engineering},
  month =         {February},
  number =        {9–12},
  pages =         {1031–1051},
  publisher =     {Elsevier BV},
  title =         {The Stochastic Finite Element Method: Past, present
                   and future},
  volume =        {198},
  year =          {2009},
  doi =           {10.1016/j.cma.2008.11.007},
  issn =          {0045-7825},
}

@article{QuainiSymmetryBreakingPreliminary2016,
  author =        {Quaini, A. and Glowinski, R. and {\v C}ani{\'c}, S.},
  journal =       {International Journal of Computational Fluid
                   Dynamics},
  number =        {1},
  pages =         {7--19},
  publisher =     {{Taylor \& Francis}},
  title =         {Symmetry Breaking and Preliminary Results about a
                   {Hopf} Bifurcation for Incompressible Viscous Flow in
                   an Expansion Channel},
  volume =        {30},
  year =          {2016},
  doi =           {10.1080/10618562.2016.1144877},
}

@article{Breden2023,
  author =        {Breden, Maxime},
  journal =       {SIAM Journal on Applied Dynamical Systems},
  month =         {June},
  number =        {2},
  pages =         {765–801},
  publisher =     {Society for Industrial & Applied Mathematics (SIAM)},
  title =         {A Posteriori Validation of Generalized Polynomial
                   Chaos Expansions},
  volume =        {22},
  year =          {2023},
  doi =           {10.1137/22m1493197},
  issn =          {1536-0040},
  url =           {http://dx.doi.org/10.1137/22M1493197},
}

@incollection{CalozNumericalAnalysisNonlinear1997,
  author =        {Caloz, Gabriel and Rappaz, Jacques},
  booktitle =     {Handbook of Numerical Analysis},
  pages =         {487--637},
  publisher =     {{Elsevier}},
  series =        {Techniques of Scientific Computing (Part 2)},
  title =         {Numerical Analysis for Nonlinear and Bifurcation
                   Problems},
  volume =        {5},
  year =          {1997},
  doi =           {10.1016/S1570-8659(97)80004-X},
}

@article{Cliffe2000,
  author =        {Cliffe, K. A. and Spence, A. and Tavener, S. J.},
  journal =       {Acta Numerica},
  month =         {January},
  pages =         {39–131},
  publisher =     {Cambridge University Press (CUP)},
  title =         {The Numerical Analysis of Bifurcation Problems with
                   application to Fluid Mechanics},
  volume =        {9},
  year =          {2000},
  doi =           {10.1017/S0962492900000398},
  issn =          {1474-0508},
}

@article{Pintore2020,
  author =        {Pintore, Moreno and Pichi, Federico and Hess, Martin and
                   Rozza, Gianluigi and Canuto, Claudio},
  journal =       {Advances in Computational Mathematics},
  publisher =     {Springer Science and Business Media LLC},
  title =         {Efficient computation of Bifurcation Diagrams with a
                   Deflated Approach to Reduced Basis Spectral Element
                   Method},
  volume =        {47},
  year =          {2020},
  doi =           {10.1007/s10444-020-09827-6},
  issn =          {1572-9044},
}

@article{FarrellDeflationTechniquesFinding2015,
  author =        {Farrell, P. E. and Birkisson, {\'A}. and
                   Funke, S. W.},
  journal =       {SIAM Journal on Scientific Computing},
  publisher =     {{Society for Industrial and Applied Mathematics}},
  title =         {Deflation Techniques for Finding Distinct Solutions
                   of Nonlinear Partial Differential Equations},
  volume =        {37},
  year =          {2015},
  doi =           {10.1137/140984798},
}

@article{PichiStrazzullo2022,
  author =        {Pichi, Federico and Strazzullo, Maria and
                   Ballarin, Francesco and Rozza, Gianluigi},
  journal =       {ESAIM: Mathematical Modelling and Numerical Analysis},
  pages =         {1361--1400},
  publisher =     {{EDP Sciences}},
  title =         {Driving Bifurcating Parametrized Nonlinear {PDEs} by
                   Optimal Control Strategies: Application to
                   {Navier-Stokes} Equations with Model Order Reduction},
  volume =        {56},
  year =          {2022},
  doi =           {10.1051/m2an/2022044},
}

@article{PichiDeflationbasedCertifiedGreedy2025a,
  title = {Deflation-Based Certified Greedy Algorithm and Adaptivity for Bifurcating Nonlinear {{PDEs}}},
  author = {Pichi, Federico and Strazzullo, Maria},
  year = 2025,
  journal = {Communications in Nonlinear Science and Numerical Simulation},
  volume = {149},
  pages = {108941},
  doi = {10.1016/j.cnsns.2025.108941}
}

@book{hesthaven2015certified,
  author =        {J. S. Hesthaven and G. Rozza and B. Stamm},
  publisher =     {Springer International Publishing},
  series =        {SpringerBriefs in Mathematics},
  title =         {Certified Reduced Basis Methods for Parametrized
                   Partial Differential Equations},
  year =          {2015},
  doi =           {10.1007/978-3-319-22470-1},
  isbn =          {9783319224701},
}

@book{benner2017model,
  author =        {Benner, P. and Cohen, A. and Ohlberger, M. and
                   Willcox, K.},
  publisher =     {SIAM, Society for Industrial and Applied Mathematics},
  title =         {Model Reduction and Approximation: Theory and
                   Algorithms},
  year =          {2017},
  doi =           {10.1137/1.9781611974829},
  isbn =          {9781611974812},
}

@article{Pitton_2017,
  author =        {Pitton, Giuseppe and Rozza, Gianluigi},
  journal =       {Journal of Scientific Computing},
  number =        {1},
  pages =         {157--177},
  publisher =     {Springer},
  title =         {On the application of reduced basis methods to
                   bifurcation problems in incompressible fluid
                   dynamics},
  volume =        {73},
  year =          {2017},
  doi =           {10.1007/s10915-017-0419-6},
}

@article{Khamlich2022,
  author =        {Khamlich, Moaad and Pichi, Federico and
                   Rozza, Gianluigi},
  journal =       {International Journal for Numerical Methods in
                   Fluids},
  month =         {June},
  publisher =     {Wiley},
  title =         {Model Order Reduction for Bifurcating Phenomena in
                   Fluid‐Structure Interaction Problems},
  volume =        {94},
  year =          {2022},
  issn =          {1097-0363},
  doi =           {https://doi.org/10.1002/fld.5118},
}

@article{BruntonDiscoveringGoverningEquations2016,
  author =        {Brunton, Steven L and Proctor, Joshua L and
                   Kutz, J Nathan},
  journal =       {Proceedings of the national academy of sciences},
  number =        {15},
  pages =         {3932--3937},
  publisher =     {National Acad Sciences},
  title =         {Discovering governing equations from data by sparse
                   identification of nonlinear dynamical systems},
  volume =        {113},
  year =          {2016},
  doi =           {10.1073/pnas.1517384113},
}

@article{PichiArtificialNeuralNetwork2023,
  author =        {Pichi, Federico and Ballarin, Francesco and
                   Rozza, Gianluigi and Hesthaven, Jan S.},
  journal =       {Computers \& Fluids},
  pages =         {105813},
  title =         {An Artificial Neural Network Approach to Bifurcating
                   Phenomena in Computational Fluid Dynamics},
  volume =        {254},
  year =          {2023},
  doi =           {10.1016/j.compfluid.2023.105813},
}

@article{PichiGraphConvolutionalAutoencoder2024,
  author =        {Pichi, Federico and Moya, Beatriz and
                   Hesthaven, Jan S.},
  journal =       {Journal of Computational Physics},
  pages =         {112762},
  title =         {A Graph Convolutional Autoencoder Approach to Model
                   Order Reduction for Parametrized PDEs},
  volume =        {501},
  year =          {2024},
  doi =           {10.1016/j.jcp.2024.112762},
}

@article{coscia2024generative,
  author =        {Coscia, Dario and Demo, Nicola and Rozza, Gianluigi},
  journal =       {Scientific Reports},
  number =        {1},
  pages =         {3826},
  publisher =     {Nature Publishing Group},
  title =         {Generative Adversarial Reduced Order Modelling},
  volume =        {14},
  year =          {2024},
  doi =           {10.1038/s41598-024-54067-z},
}

@book{uq-book,
  author =        {Smith, Ralph C.},
  month =         {January},
  publisher =     {Society for Industrial and Applied Mathematics},
  title =         {Uncertainty Quantification: Theory, Implementation,
                   and Applications},
  year =          {2013},
  doi =           {10.1137/1.9781611973228},
  isbn =          {9781611973228},
}

@book{arnold2010,
  address =       {Berlin, Heidelberg},
  author =        {Arnold, Ludwig},
  publisher =     {Springer},
  title =         {Random Dynamical Systems},
  year =          {1998},
  doi =           {10.1007/BFb0095238},
  isbn =          {978-3-642-08355-6 978-3-662-12878-7},
}

@article{SriNamachchivaya1990,
  author =        {N. {Sri Namachchivaya}},
  journal =       {Applied Mathematics and Computation},
  number =        {2},
  title =         {Stochastic Bifurcation},
  volume =        {38},
  year =          {1990},
  doi =           {10.1016/0096-3003(90)90003-L},
  issn =          {0096-3003},
}

@article{Blumenthal2023,
  author =        {{Blumenthal}, Alex and {Engel}, Maximilian and
                   {Neamtu}, Alexandra},
  journal =       {Probability Theory and Related Fields},
  month =         {August},
  title =         {On the Pitchfork Bifurcation for the {Chafee-Infante}
                   equation with additive noise},
  year =          {2021},
  doi =           {10.1007/s00440-023-01235-3},
}

@article{Liu2024,
  author =        {Meng Liu},
  journal =       {Journal of Mathematical Analysis and Applications},
  number =        {1},
  pages =         {128096},
  title =         {Stability and Dynamical Bifurcation of a Stochastic
                   regime-switching Predator–Prey model},
  volume =        {535},
  year =          {2024},
  issn =          {0022-247X},
  doi =           {https://doi.org/10.1016/j.jmaa.2024.128096},
}

@article{Venturi2010,
  author =        {Venturi, Daniele and Wan, Xiaoliang and
                   Karniadakis, George Em},
  journal =       {Journal of Fluid Mechanics},
  pages =         {391–413},
  title =         {Stochastic Bifurcation analysis of
                   {Rayleigh–Bénard} convection},
  volume =        {650},
  year =          {2010},
  doi =           {https://doi.org/10.1017/S0022112009993685},
}

@book{Rubinstein2016,
  author =        {Rubinstein, Reuven Y. and Kroese, Dirk P.},
  journal =       {Wiley Series in Probability and Statistics},
  month =         {November},
  publisher =     {Wiley},
  title =         {Simulation and the Monte Carlo Method},
  year =          {2016},
  doi =           {10.1002/9780470285312},
  isbn =          {9781118631980},
  issn =          {1940-6347},
}

@phdthesis{wang2008karhunen,
  author =        {Wang, Limin},
  school =        {London School of Economics and Political Science
                   (United Kingdom)},
  title =         {{Karhunen-Loeve} expansions and their applications},
  year =          {2008},
}

@mastersthesis{arcozzi,
  author =        {Arcozzi, Nicola and Campanino, Massimo and
                   Giambartolomei, Giordano},
  school =        {University of Bologna},
  title =         {The {K}arhunen-{L}oeve Theorem},
  year =          {2015},
}

@article{huang2001convergence,
  author =        {Huang, SP and Quek, ST and Phoon, KK0994},
  journal =       {International Journal for Numerical Methods in
                   Engineering},
  number =        {9},
  pages =         {1029--1043},
  publisher =     {Wiley Online Library},
  title =         {Convergence study of the truncated
                   {Karhunen-Lo{\`e}ve} expansion for simulation of
                   Stochastic Processes},
  volume =        {52},
  year =          {2001},
  doi =           {https://doi.org/10.1002/nme.255},
}

@article{CameronMartin1947,
  author =        {Robert H. Cameron and William Ted Martin},
  journal =       {Annals of Mathematics},
  pages =         {385},
  title =         {The Orthogonal Development of Non-Linear Functionals
                   in Series of {Fourier-Hermite} Functionals},
  volume =        {48},
  year =          {1947},
  doi =           {10.2307/1969178},
}

@book{GhanemSpanos1990,
  address =       {New York, NY},
  author =        {Ghanem, Roger G. and Spanos, Pol D.},
  publisher =     {Springer},
  title =         {Stochastic Finite Elements: A Spectral Approach},
  year =          {1991},
  doi =           {10.1007/978-1-4612-3094-6},
  isbn =          {978-1-4612-7795-8 978-1-4612-3094-6},
}

@article{XIU2002_gpc,
  author =        {Xiu, Dongbin and Karniadakis, George Em},
  journal =       {SIAM Journal on Scientific Computing},
  number =        {2},
  pages =         {619-644},
  title =         {The {Wiener-Askey Polynomial Chaos} for Stochastic
                   Differential Equations},
  volume =        {24},
  year =          {2002},
  doi =           {10.1137/S1064827501387826},
}

@article{ernst2012,
  author =        {Ernst, Oliver G and Mugler, Antje and
                   Starkloff, Hans-J{\"o}rg and Ullmann, Elisabeth},
  journal =       {ESAIM: Mathematical Modelling and Numerical Analysis},
  number =        {2},
  publisher =     {EDP Sciences},
  title =         {On the convergence of generalized {Polynomial Chaos}
                   Expansions},
  volume =        {46},
  year =          {2012},
  doi =           {10.1051/m2an/2011045},
}

@article{babuvska2007stochastic,
  author =        {Babu{\v{s}}ka, Ivo and Nobile, Fabio and
                   Tempone, Ra{\'u}l},
  journal =       {SIAM Journal on Numerical Analysis},
  number =        {3},
  pages =         {1005--1034},
  publisher =     {SIAM},
  title =         {A Stochastic Collocation Method for Elliptic Partial
                   Differential Equations with Random Input Data},
  volume =        {45},
  year =          {2007},
  doi =           {https://doi.org/10.1137/050645142},
}

@book{xiu2010,
  author =        {Xiu, Dongbin},
  publisher =     {Princeton university press},
  title =         {Numerical Methods for Stochastic computations: a
                   Spectral Method Approach},
  year =          {2010},
  doi =           {10.1515/9781400835348},
}

@book{sullivan2015,
  author =        {Sullivan, Timothy John},
  publisher =     {Springer},
  title =         {Introduction to uncertainty quantification},
  volume =        {63},
  year =          {2015},
  doi =           {10.1007/978-3-319-23395-6},
}

@article{Sobey,
  author =        {Sobey, Ian J. and Drazin, Philip G.},
  journal =       {Journal of Fluid Mechanics},
  pages =         {263–287},
  publisher =     {Cambridge University Press},
  title =         {Bifurcations of two-dimensional channel Flows},
  volume =        {171},
  year =          {1986},
  doi =           {10.1017/S0022112086001441},
}

@article{Cherdron,
  author =        {Cherdron, W. and Durst, F. and Whitelaw, J\. H\.},
  journal =       {Journal of Fluid Mechanics},
  number =        {1},
  pages =         {13–31},
  publisher =     {Cambridge University Press},
  title =         {Asymmetric Flows and Instabilities in symmetric ducts
                   with sudden expansions},
  volume =        {84},
  year =          {1978},
  doi =           {10.1017/S0022112078000026},
}

@article{BravoGeometricallyParametrisedReduced2023,
  author =        {Bravo, J. R. and Stabile, G. and Hess, M. and
                   Hernandez, J. A. and Rossi, R. and Rozza, G.},
  journal =       {Journal of Computational Physics},
  pages =         {113058},
  title =         {Geometrically Parametrised Reduced Order Models for
                   Studying the Hysteresis of the {{Coanda}} Effect in
                   Finite Element-Based Incompressible Fluid Dynamics},
  volume =        {509},
  year =          {2024},
  doi =           {10.1016/j.jcp.2024.113058},
}

@article{TonicelloNonintrusiveReducedOrder2022,
  author =        {Tonicello, Niccol{\`o} and Lario, Andrea and
                   Rozza, Gianluigi and Mengaldo, Gianmarco},
  journal =       {Computers \& Fluids},
  pages =         {106307},
  title =         {Non-Intrusive Reduced Order Models for the Accurate
                   Prediction of Bifurcating Phenomena in Compressible
                   Fluid Dynamics},
  volume =        {278},
  year =          {2024},
  doi =           {10.1016/j.compfluid.2024.106307},
}

@article{HessLocalizedReducedorderModeling2019,
  author =        {Hess, Martin and Alla, Alessandro and
                   Quaini, Annalisa and Rozza, Gianluigi and
                   Gunzburger, Max},
  journal =       {Computer Methods in Applied Mechanics and
                   Engineering},
  pages =         {379--403},
  title =         {A Localized Reduced-Order Modeling Approach for
                   {PDEs} with Bifurcating Solutions},
  volume =        {351},
  year =          {2019},
  doi =           {10.1016/j.cma.2019.03.050},
}

@article{sousedik2016,
  author =        {Soused{\'\i}k, Bed{\v{r}}ich and Elman, Howard C},
  journal =       {Journal of Computational Physics},
  pages =         {435--452},
  publisher =     {Elsevier},
  title =         {{Stochastic Galerkin} Methods for the Steady-State
                   {Navier-Stokes} equations},
  volume =        {316},
  year =          {2016},
  doi =           {10.1016/j.jcp.2016.04.013},
}

@article{Lee2019,
  author =        {Lee, Kookjin and Elman, Howard C. and
                   Soused\'{\i}k, Bed\v{r}ich},
  journal =       {SIAM/ASA Journal on Uncertainty Quantification},
  pages =         {1275-1300},
  title =         {A Low-Rank Solver for the {Navier-Stokes} Equations
                   with Uncertain Viscosity},
  volume =        {7},
  year =          {2019},
  doi =           {10.1137/17M1151912},
}

@article{Girolami2021,
  author =        {Girolami, Mark and Febrianto, Eky and Yin, Ge and
                   Cirak, Fehmi},
  journal =       {Computer Methods in Applied Mechanics and
                   Engineering},
  month =         {March},
  pages =         {113533},
  publisher =     {Elsevier BV},
  title =         {The statistical Finite Element Method (statFEM) for
                   coherent synthesis of observation data and model
                   predictions},
  volume =        {375},
  year =          {2021},
  doi =           {10.1016/j.cma.2020.113533},
  issn =          {0045-7825},
}

@book{KuznetsovElementsAppliedBifurcation2023,
  title = {Elements of Applied Bifurcation Theory},
  author = {Kuznetsov, Yuri A.},
  year = 2023,
  series = {Applied Mathematical Sciences},
  volume = {112},
  publisher = {Springer International Publishing},
  address = {Cham},
  doi = {10.1007/978-3-031-22007-4},
  isbn = {978-3-031-22006-7 978-3-031-22007-4}
}

@article{HallerNonlinearNormalModes2016,
  title = {Nonlinear Normal Modes and Spectral Submanifolds: Existence, Uniqueness and Use in Model Reduction},
  shorttitle = {Nonlinear Normal Modes and Spectral Submanifolds},
  author = {Haller, George and Ponsioen, Sten},
  year = 2016,
  journal = {Nonlinear Dynamics},
  volume = {86},
  number = {3},
  pages = {1493--1534},
  doi = {10.1007/s11071-016-2974-z}
}

@article{BettiniDatadrivenNonlinearModel2025,
  title = {Data-Driven Nonlinear Model Reduction to Spectral Submanifolds via Oblique Projection},
  author = {Bettini, Leonardo and Kasz{\'a}s, B{\'a}lint and Zybach, Bernhard and Dual, J{\"u}rg and Haller, George},
  year = 2025,
  journal = {Chaos: An Interdisciplinary Journal of Nonlinear Science},
  volume = {35},
  number = {4},
  pages = {043135},
  doi = {10.1063/5.0243849}
}

\appendix

\section{Karhunen-Lo{\`e}ve expansion} \label{section:KL_expansion}
The Karhunen-Lo{\`e}ve (K-L) expansion \cite{wang2008karhunen, arcozzi} operates similarly to the Fourier expansion but it applies to the space of stochastic processes with finite second-order moment $\set{L}_2(\Omega, \vectspace{F}, \nu; \set{R}^d)$.
The main idea of the method is to first find the orthonormal eigenfunctions $\vect{\pi}_{i}(\vect{x})$ of
the Hilbert-Schmidt integral operator $\vectspace{H}_K$ defined as:
\begin{equation}\label{eq:KL_coefficient}
    \vectspace{H}_K[\vect{\pi}_{i}(\vect{x})](\vect{x}) = \int_{\vectspace{D}}
    K_{\vect{\alpha}}(\vect{x}, \vect{y}) \vect{\pi}_{i}(\vect{y})
    d\vect{y} = \lambda_i\vect{\pi}_{i}(\vect{x}),
\end{equation}
whose kernel $K_{\vect{\alpha}}(\vect{x}, \vect{y})$ is the covariance function
of the process $\vect{\alpha}(\vect{x}, \omega) \in \set{L}_2(\Omega, \vectspace{F}, \nu; \set{R}^d)$,  and
$\lambda_i$ are the eigenvalues associated with $\vect{\pi}_{i}(\vect{x})$.
Subsequently, the process is projected onto such an orthonormal basis, obtaining the random variables:
\begin{equation}\label{eq:KL_rv}
    \xi_i(\omega) = \int_{\vectspace{D}} (\vect{\alpha}(\vect{x}, \omega)-\overline{\vect{\alpha}}(\bm{x}))
    \vect{\pi}_{i}(\vect{x})d\vect{x} \quad \text{s.t.}\quad \set{E}[\xi_i]=0,\ \text{and} \ \set{E}[\xi_i\xi_j]=\lambda_i\delta_{i,j},
\end{equation}
where the conditions on their mean and variance are easily obtained by substituting the $\vect{\pi}_{i}(\vect{x})$ functions with their definition in Equation \eqref{eq:KL_coefficient}. Namely, they are orthogonal in the space $\set{L}_2(\Omega, \vectspace{F}, \nu; \set{R}^d)$, and thus uncorrelated. The K-L expansion can then be expressed as:
\begin{equation}\label{eq:KL_normalized}
    \vect{\alpha}(\vect{x}, \omega) = \overline{\vect{\alpha}}(\vect{x}) +
    \sum_{i=1}^{\infty}\frac{\xi_i(\omega)}{\sqrt{\lambda_i}}
    \vect{\pi}_{i}(\vect{x})\sqrt{\lambda_i} 
    = \sum_{i=0}^{\infty}\overline{\xi_i}(\omega)\vect{\pi}_{i}(\vect{x})\sqrt{\lambda_i},
\end{equation}
where the first two equivalences follow directly from introducing the normalized random variables $\overline{\xi_i}(\omega)$ such that $\set{E}[\overline{\xi_i}]=0$ and
$\set{E}[\overline{\xi_i}\overline{\xi_j}]=\delta_{i,j}$, and assuming without loss of generality that $\overline{\xi_0}(\omega)\sim\delta_{1}$,
$\vect{\pi}_{0}(\vect{x})=\overline{\vect{\alpha}}(\vect{x})$, and $\lambda_0=1$. 

The convergence properties of this expansion strictly depend on the nature of the stochastic process covariance function, namely on the decay of its eigenvalues. For a further discussion on the topic, we refer to \cite{huang2001convergence}. 
In general, note that the K-L expansion is particularly useful to model the parameters of a generic stochastic PDE \eqref{eq:stochasticPDE}, since their covariance function is known and all the computation comes down to its spectral study.

\section{Non-intrusive stochastic perturbation approach} \label{sec:MC}

A sample-based approach to bifurcation problems may be misleading because it assumes access to all solution branches for a given parameter set, a requirement that cannot be verified without prior knowledge of the system.
Indeed, without exploiting ad-hoc techniques such as the deflation method \cite{FarrellDeflationTechniquesFinding2015}, there are no guarantees that the nonlinear solver can find samples of the solution equally distributed on the possible branches. 
For instance, approaching the problem via Monte-Carlo simulations, or similar non-intrusive methodologies, one collects an ensemble of solutions by repeatedly solving the system for different values of the parameters, without obtaining any assurance on the completeness of the analysis. 

We empirically prove this argument to justify the need of intrusive methods when assuming no previous knowledge of the system. 
We treat the viscosity in \eqref{eq:coanda} as a normally distributed low-variance random variable $\mu\sim\vectspace{N}(0.9,0.001)$ centered at a point of the parameter space where the bifurcation has already occurred, and we draw $300$ samples from such parameter distribution, solving the following nonlinear system coming from FEM formulation of the problem and collecting an ensemble of solutions:
\begin{equation} \label{eq:coanda_wf}
\begin{cases}
\mu\int_{\vectspace{D}}\nabla \vect{v}\cdot \nabla\vect{\psi} d\vect{x} +
\int_{\vectspace{D}}(\vect{v}\cdot\nabla \vect{v})\vect{\psi} d\vect{x}
- \int_{\vectspace{D}}p\nabla\cdot\vect{\psi} d\vect{x}, \quad
&\forall \vect{\psi} \in \set{V}_0\\
\int_{\vectspace{D}}\pi\nabla\cdot \vect{v} d\vect{x}=0, \quad
&\forall \pi \in \set{Q},
\end{cases}
\end{equation}
where $\set{V}_{0} =\left\{\boldsymbol{v} \in H^{1}\left(\vectspace{D}; \mathbb{R}^{2}\right):\left.\vect{\psi}\right|_{\Gamma_{\text{wall}} \cup \Gamma_{\text{in}}} =0\right\}$, $ 
\set{V}_\text{in} = \set{V}_{0} \oplus \vect{v}_\text{in}$, and
$\set{Q}=L^{2}(\vectspace{D})$.

\begin{figure}[t]
\centering
        \captionsetup[sub][figure]{oneside,margin={-7cm,0cm},skip=-1.6cm}
        \subcaptionbox{\label{subfig:MC_noprior}}{        \includegraphics[trim={12cm 10.3cm 6.0cm 11.cm}, clip, width=.65\textwidth]{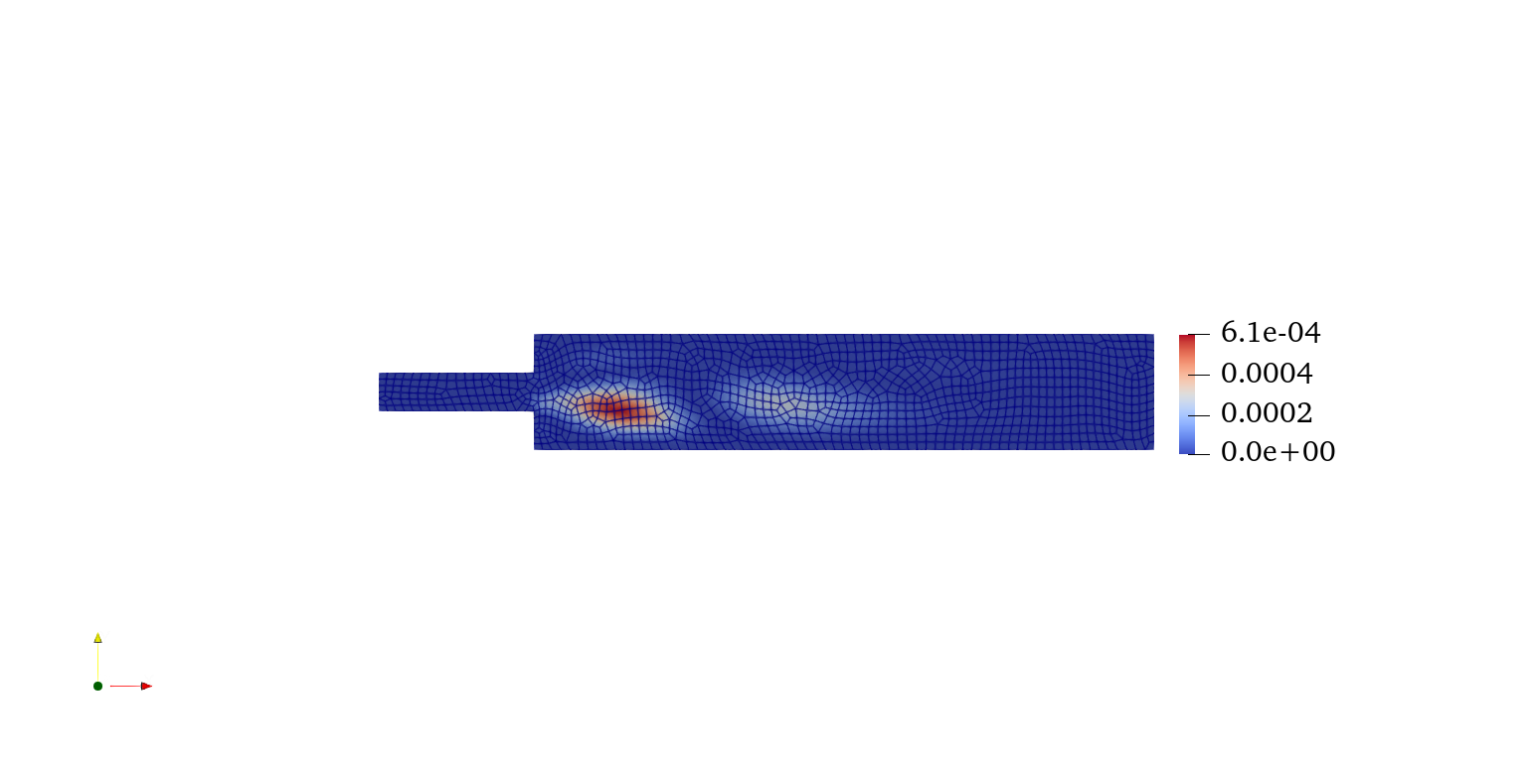}
        }
    \hfill
        \subcaptionbox{\label{subfig:MC_prior}}{\includegraphics[trim={12cm 10.3cm 6.0cm 11.cm}, clip, width=.65\textwidth]{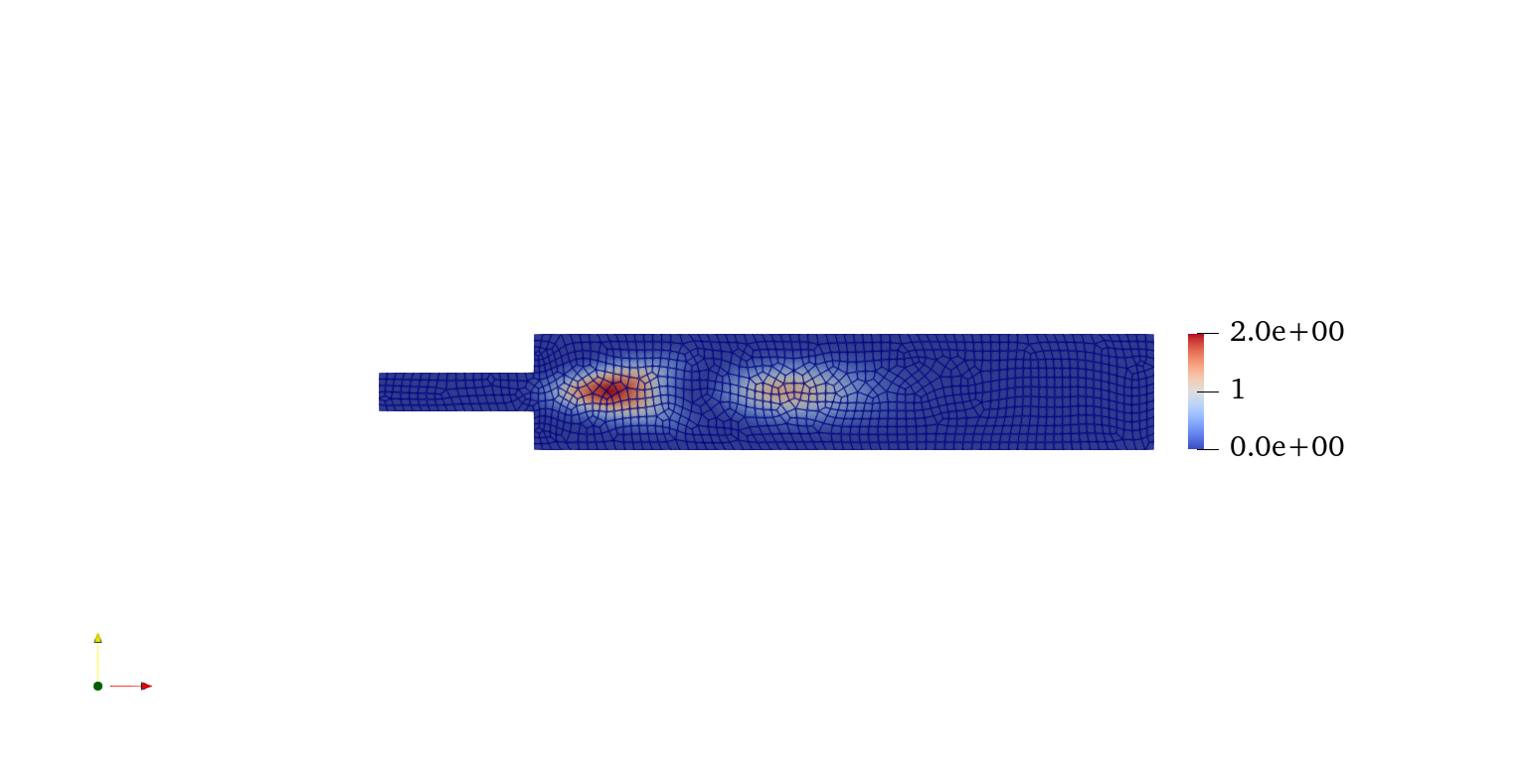}
        \captionsetup[sub][figure]{oneside,margin={-6cm,0cm},skip=-2cm}}
        \captionsetup{subrefformat=parens}
    \caption{Variance of $v_y$ computed via MC simulation of $300$ samples with $\mu\sim\vectspace{N}(0.9,0.001)$ on a quadrilateral grid of $1275$ nodes, \cref{subfig:MC_noprior} fixed and \cref{subfig:MC_prior} ad-hoc knowledge-based initialization.}
    \label{fig:MC}
\end{figure}

Figure \ref{subfig:MC_noprior} shows that, if the nonlinear solver initialization is not exploiting any prior knowledge regarding the possible states, i.e., the initial guess is set to zero, the variance of the collected ensemble has a small magnitude in the whole domain. 
On the other hand, if we incorporate some knowledge by means of a continuation approach for the initialization of the Newton-Krylov solver, the variance clearly shows a different behavior with an increased magnitude, concentrating in the expansion-region of the channel where the bifurcation is mostly affecting the flow.

\begin{figure}[t]
\centering
\captionsetup[sub][figure]{oneside,margin={-5cm,0cm},skip=-3cm}
        \subcaptionbox{\label{subfig:MCsamples_noprior}}{        \includegraphics[trim={2.4cm 1.6cm 2.cm 2cm}, clip, width=0.305\textwidth]{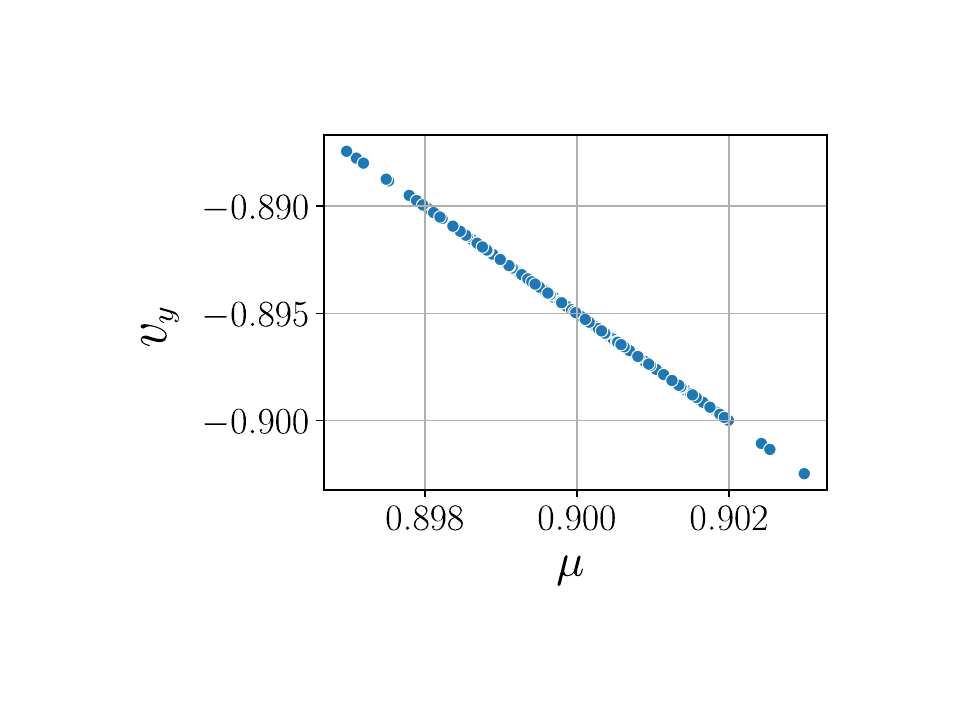}
        }
    \hfill
        \captionsetup[sub][figure]{oneside,margin={-5cm,0cm},skip=-3cm}
        \subcaptionbox{\label{subfig:MCsamples_prior}}{        \includegraphics[trim={2.4cm 1.6cm 2cm 2cm}, clip, width=0.305\textwidth]{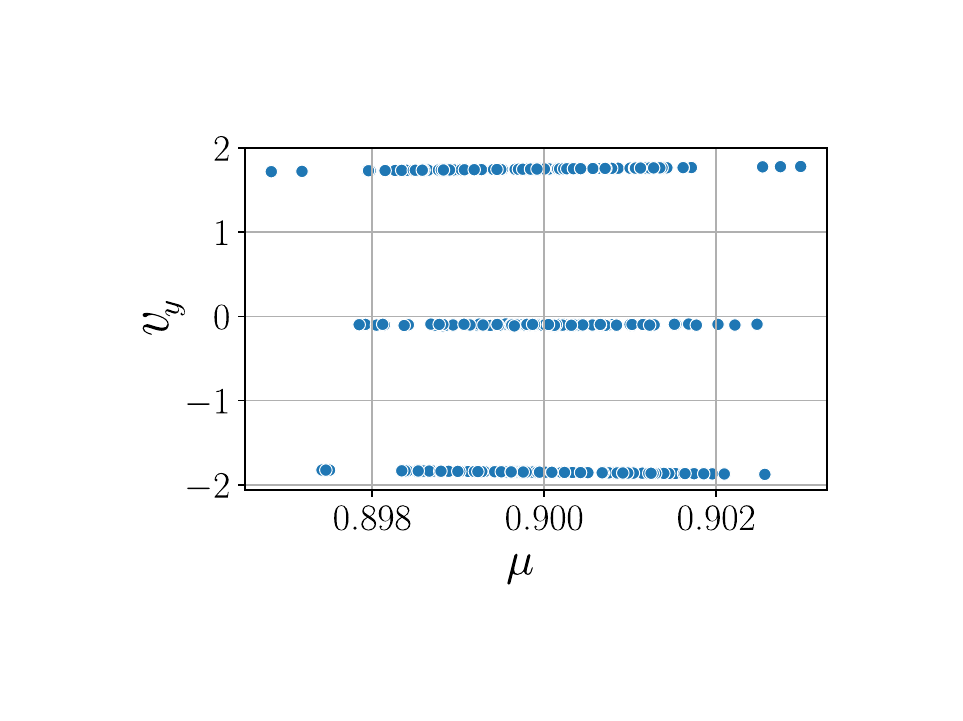}
        }
     \hfill
       \captionsetup[sub][figure]{oneside,margin={-5cm,0cm},skip=-3cm}
        \subcaptionbox{\label{subfig:MCsamples_highvar}}{        \includegraphics[trim={0.5cm 0.cm 0cm 0cm}, clip, width=0.305\textwidth]{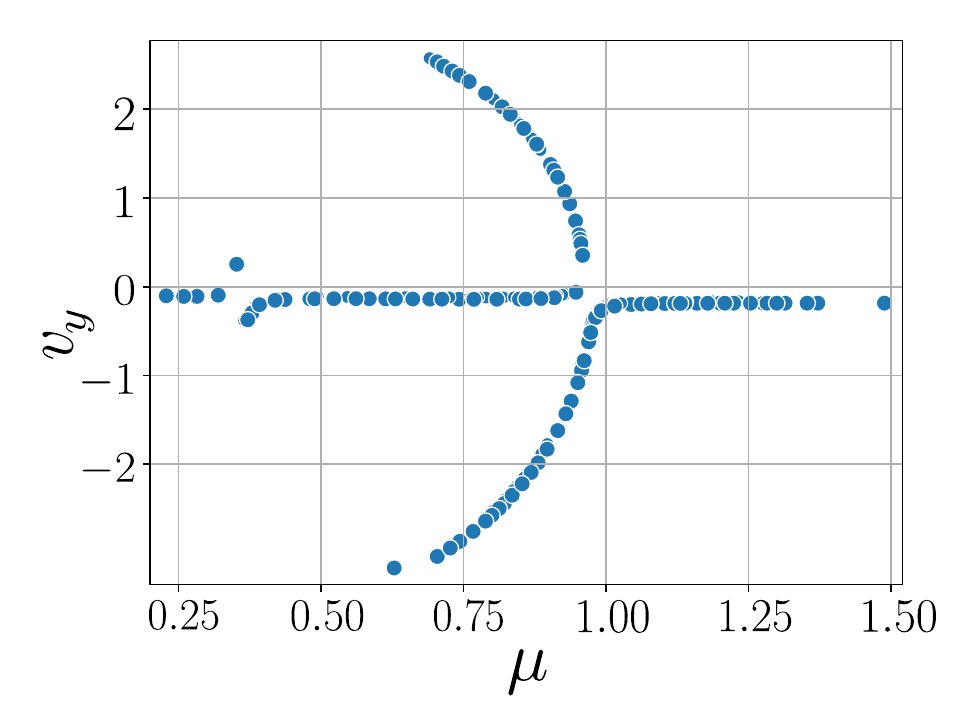}
        }
    \captionsetup{subrefformat=parens}
    \caption{Monte Carlo samples of $v_y$ at the point of maximum variance with $\mu\sim\vectspace{N}(0.9,0.001)$ for Figure \ref{subfig:MC_noprior} in \cref{subfig:MCsamples_noprior}, for Figure \ref{subfig:MC_prior} in \cref{subfig:MCsamples_prior}, and with $\mu\sim\vectspace{N}(0.9,0.2)$ in  \cref{subfig:MCsamples_highvar}.}
    \label{fig:MCsamples}
\end{figure}

%As a further confirmation of the need for additional knowledge in the non-intrusive approach, we can study the actual samples from the grid points realizing the maximum value of the variance in both cases.
Moreover, Figure \ref{subfig:MCsamples_noprior} shows that the increased variance observed in Figure \ref{subfig:MC_noprior} arises solely from a linear trend. In contrast, Figure \ref{subfig:MCsamples_prior}, which corresponds to the high variance in Figure \ref{subfig:MC_prior}, exhibits a multi-valued mapping that correctly identifies the three branches of the deterministic bifurcation diagram shown in Figure \ref{fig:bif-diag-deterministic}.
In Figure \ref{subfig:MCsamples_highvar}, we exploit guided initialization over a broader parameter region by assigning $\mu\sim\vectspace{N}(0.9,0.2)$.
In this case, the topology of the bifurcation diagram is not only reproduced, but an additional bifurcating behavior is also observed around $\mu=0.5$. 
{Thus, although sample-based approaches are often easier to interpret and implement, they lack reliability in regimes of non-uniqueness when no prior knowledge of the system is available. 
Nevertheless, when coupled with prior information or dedicated bifurcation-discovery techniques, they may still serve as a tool for investigation of the parameters space \cite{PichiArtificialNeuralNetwork2023,PichiGraphConvolutionalAutoencoder2024,TonicelloNonintrusiveReducedOrder2022}. 
}

\end{document}